\documentclass{article}

\usepackage{amsmath,amsthm,amsfonts,amssymb}
\usepackage[all]{xy}
\usepackage[margin=1in]{geometry}
\usepackage{dsfont}
\usepackage{mathtools}
\usepackage{hyperref}
\usepackage{tikz}
\usetikzlibrary{matrix}

\newcommand\com{\mathrm{com}}
\newcommand{\Bcom}{B_{\com}}
\newcommand{\Ecom}{E_{\com}}
\newcommand{\RP}{\mathbb{RP}}
\newcommand{\CP}{\mathbb{CP}}
\newcommand{\Z}{\mathbb{Z}}
\newcommand{\R}{\mathbb{R}}
\newcommand{\F}{\mathbb{F}}
\newcommand{\Q}{\mathbb{Q}}
\newcommand{\hofib}{\mathop{\mathrm{hofib}}}
\newcommand{\hocofib}{\mathop{\mathrm{hocofib}}}
\newcommand{\im}{\mathop{\mathrm{im}}}
\newcommand{\id}{\mathrm{id}}
\newcommand{\Sq}{\mathrm{Sq}}
\newcommand{\Gr}{\mathrm{Gr}}
\newcommand{\Map}{\mathrm{Map}}
\newcommand{\Hom}{\mathrm{Hom}}
\newcommand{\one}{\mathds{1}}

\newtheorem{theorem}{Theorem}[section]
\newtheorem*{introthm*}{Theorem \ref{\introthmref}}
\newtheorem*{introcor*}{Corollary \ref{\introthmref}}
\newtheorem{lemma}[theorem]{Lemma}
\newtheorem{proposition}[theorem]{Proposition}
\newtheorem{corollary}[theorem]{Corollary}
\theoremstyle{definition}
\newtheorem{remark}[theorem]{Remark}
\newtheorem{question}[theorem]{Question}

\title{\vspace{-1cm}Classifying spaces for commutativity of low-dimensional Lie groups}
\author{Omar Antol\'{\i}n-Camarena\thanks{Instituto de Matem\'aticas,
    UNAM, Mexico City. \emph{E-mail address}: \href{mailto:omar@matem.unam.mx}{\texttt{omar@matem.unam.mx}}}
  \and
  Simon Philipp Gritschacher\thanks{Centre for Symmetry and Deformation,
    University of Copenhagen. \emph{E-mail address}: \href{mailto:gritschacher@math.ku.dk}{\texttt{gritschacher@math.ku.dk}} This author gratefully acknowledges financial support from the London Mathematical Society through a Postdoctoral Mobility Grant (PMG 16-17 22), and would like to thank the Pacific Institute for the Mathematical Sciences at the University of British Columbia for their hospitality. The author is supported by the Danish National Research Foundation through the Centre for Symmetry and Deformation (DNRF92).}
  \and
  Bernardo Villarreal\thanks{Instituto de Matem\'aticas,
    UNAM, Mexico City. \emph{E-mail address}: \href{mailto:villarreal@matem.unam.mx}{\texttt{villarreal@matem.unam.mx}}}}

\begin{document}

\maketitle

\vspace{-5mm}

\begin{abstract}
  For each of the groups $G = O(2), SU(2), U(2)$, we compute the
  integral and $\F_2$-cohomology rings of $\Bcom G$ (the classifying
  space for commutativity of $G$), the action of the Steenrod algebra
  on the mod 2 cohomology, the homotopy type of $\Ecom G$ (the
  homotopy fiber of the inclusion $\Bcom G \to BG$), and some
  low-dimensional homotopy groups of $\Bcom G$.
\end{abstract}

\tableofcontents

\section{Introduction}

Let $G$ be a topological group for which the identity element $1_G\in G$ is a non-degenerate basepoint. The set of homomorphisms $\Hom(\Z^n,G)$ can be topologized as the subspace of $G^n$ consisting of $n$-tuples $(g_1,...,g_n)$ such that $g_ig_j=g_jg_i$ for all $1\leq i,j\leq n$. The inclusions $i_n\colon \Hom(\Z^n,G)\to G^n$ are compatible with the degeneracy and face maps arising form $NG$, the nerve of $G$ thought of as a (topological) category with one object, which means that the collection $\{\Hom(\Z^n,G)\}_{n\geq0}$ assembles into a simplicial subspace $\Hom(\Z^*,G)$ of $NG$ -- the latter being a simplicial model of the classifying space $BG$. Let
\[\Bcom G:=|\Hom(\Z^*,G)|.\]
This space was originally introduced in \cite{Ad5}, and further studied in \cite{Ad1}. The latter authors among other results, show that $\Bcom G$ classifies principal $G$-bundles whose transition functions $\rho_{\alpha\beta}\colon U_\alpha\cap U_\beta\to G$ are \emph{transitionally commutative}, that is, the functions $\rho_{\alpha\beta}$ commute whenever they are simultaneously defined. In particular, the pullback of the universal bundle $EG\to BG$ along the inclusion $\iota\colon\Bcom G\to BG$ has a transitionally commutative structure. This bundle is denoted $\Ecom G\to \Bcom G$.

So far most of the attention in the literature, both for proving general results and in examples, has been given to compact connected Lie groups. The focus on compact groups is justified in part by the main result of \cite{PettetSouto} which states that if $G$ is the group of complex (or real) points of a (real) reductive algebraic group, and $K \subset G$ is a maximal compact subgroup then the induced map $\Hom(\Z^n, K) \to \Hom(\Z^n, G)$ is a strong deformation retract. As shown in \cite[Theorem 3.1]{Ad5}, this implies that the induced maps $\Bcom K \to \Bcom G$ and $\Ecom K \to \Ecom G$ are homotopy equivalences. This theorem is surprising because the commuting tuples in $G$ can be much more complicated than those in $K$. Notice, for example, that the theorem applies to $G = SL(2, \R)$ whose maximal compact subgroup is $K = SO(2)$, which happens to be abelian. Thus, $\Bcom SL(2, \R) \simeq BSO(2)$ and $\Ecom SL(2, \R)$ is contractible.

\paragraph*{Summary of our main results.} This paper is devoted to concrete calculations with examples of low-dimensional compact Lie groups, calculations which refine the general structural information available from \cite{Ad1} for compact, connected Lie groups. In \cite{Ad1} the rational cohomology rings of $\Bcom G_\one$ (see the start of Section \ref{sec:poBcomG1} for an explanation of the subscript $\one$) and $\Ecom G_\one$ are computed in terms of the Weyl group actions on the rational cohomology of $BT$ and $G/T$ (where $T$, as usual, denotes a maximal torus). Namely \cite[Proposition 7.1]{Ad1} says that \[H^{\ast}(\Bcom G_{\one}; \Q) = (H^{\ast}(BT;\Q) \otimes_{\Q} H^{\ast}(G/T; \Q))^{W}\] and \cite[Corollary 7.4]{Ad1} states that \[H^{\ast}(\Ecom G_{\one}; \Q) = (H^{\ast}(G/T;\Q) \otimes_{\Q} H^{\ast}(G/T; \Q))^{W}.\] The authors note in \cite[Remark 7.3]{Ad1} that integrally those descriptions cannot hold because the example of $G=SU(2)$ shows that $H^*(\Bcom SU(2), \Z)$ contains 2-torsion. One of our main results is the description of the cup products of these interesting torsion classes.

Our example groups are $G = O(2)$, $SU(2)$, $U(2)$ and $SO(3)$ (notice that while our examples are all compact, $O(2)$ is disconnected). In the first two thirds of the paper we compute cohomology of $\Bcom G$ for these groups --- except for $G=SO(3)$, where we focus on $\Bcom SO(3)_{\one}$ instead. In each case we compute the integral cohomology ring. They turn out to all be either torsion free or to have only 2-torsion, so we also compute the $\F_2$-cohomology rings together with the action of the Steenrod algebra. The results of these calculations are collected in Theorems \ref{thm:crsu2}, \ref{thm:cru2}, \ref{thm:cro2} and \ref{thm:crso3_1} (the integral cohomology rings for the example groups $G=U(2),SU(2)$ and $SO(3)$ first appeared in \cite{VillarrealThesis}). In Theorem \ref{thm:cru2} we give a presentation of the ring
\[
H^*(\Bcom U(2);\Z) \cong \Z[c_1,c_2,y_1,y_2]/(2y_2-y_1c_1,y_1^2,y_1y_2,y_2^2)
\]
which contains the subring $H^*(BU(2);\Z)\cong \Z[c_1,c_2]$, where $c_1$ and $c_2$ are the first and second Chern classes, respectively. The embedding $H^\ast(BU(2);\Z)\hookrightarrow H^\ast(\Bcom U(2);\Z)$ is induced by the inclusion map $\iota\colon \Bcom U(2)\to BU(2)$. It turns out that the cohomology of $\Bcom SU(2)$ is easier to describe via the ring $H^*(\Bcom U(2);\Z)$. The following theorem describes the cohomology of $\Bcom SU(2)$.

\begin{samepage}
\begin{theorem}\label{thm:crsu2}\leavevmode
  \begin{enumerate}
  \item The map $\Bcom SU(2) \to \Bcom U(2)$ arising from the
    inclusion induces a surjective homomorphism of integral
    cohomology rings with kernel generated by the class $c_1$. In
    particular
    \[H^*(\Bcom SU(2);\Z)\cong\Z[c_2,y_1,x_2]/(2x_2,y_1^2,x_2y_1,x_2^2)\]
    where $c_2,y_1$ are in degree 4 and $x_2$ in degree 6. The class
    $c_2$ is pulled back from the second Chern class in
    $H^4(BSU(2);\Z)$ under the inclusion $\Bcom SU(2)\to BSU(2)$, and
    the classes $y_1$ and $x_2$ are the images of the classes $y_1$
    and $y_2$ in $H^*(\Bcom U(2); \Z)$.
    
  \item  Let $\bar{z}$ again denote reduction mod 2 of $z$. Then
    \[H^*(\Bcom SU(2);\F_2) \cong \F_2[\bar{c}_2,\bar{y}_1,x_1,\bar{x}_2]/ (\bar{y}_1^2,\bar{y}_1x_1,x_1^2,\bar{x}_2\bar{y}_1,x_1\bar{x}_2,\bar{x}_2^2),\]
    where $x_1$ has degree 5. Moreover, $\beta(x_1)=x_2$ where $\beta$
    denotes the integral Bockstein.

  \item The action of the Steenrod algebra on $H^*(\Bcom SU(2); \F_2)$
    is given by the total Steenrod squares
    \begin{align*}
      \Sq(\bar{c}_2) & = \bar{c}_2 + \bar{c}_2^2\\
      \Sq(\bar{y}_1) & = \bar{y}_1 \\
      \Sq(x_1)       & = x_1       + \bar{x}_2 \\
      \Sq(\bar{x}_2) & = \bar{x}_2 + \bar{c}_2 \bar{y}_1.\\
    \end{align*}
\end{enumerate}
\end{theorem}
\end{samepage}

We next turn to the real matrix groups $SO(2)$ and $O(2)$. For $SO(2)$ there is not much to say, since $\Bcom SO(2)=BSO(2)$. The case of $O(2)$ is more interesting. In fact, we found an intriguing relation intertwining the integral and $\F_2$-cohomology groups of $\Bcom O(2)$: They can be described entirely in terms of the cohomology of $BO(2)$, the monodromy invariants in the cohomology of $\Ecom O(2)$ and the Bockstein homomorphism. To state our result recall that $H^*(BO(2);\F_2)\cong \F_2[w_1,w_2]$ with $w_1$ and $w_2$ the first and second Stiefel-Whitney classes, and $H^*(BO(2);\Z)\cong \Z[W_1,W_2,p_1]/(2W_i,W_2^2-W_1p_1)$, where $W_i=\beta(w_i)$ and $p_1$ is the first Pontryagin class. Abusing notation, we denote by the same letters the pullback of these classes along the inclusion $\iota\colon \Bcom O(2)\to BO(2)$. 

\begin{theorem} \label{thm:cro2}\leavevmode
  \begin{enumerate}
  \item\label{thm:cro2-f2} There is an isomorphism of graded rings
\[H^*(\Bcom O(2);\F_2)\cong \F_2[w_1,w_2,\bar{r},s]/(\bar{r} w_1,\bar{r}^2, \bar{r}s, s^2)\, ,
\]
where $\deg(\bar{r})=2$, $\deg(s)=3$ and $\bar{r}$ is the reduction modulo $2$ of the integral class $r$ (see \ref{thm:cro2-hz}).

\item\label{thm:cro2-hz} There is an isomorphism of graded rings
\[H^*(\Bcom O(2);\Z)\cong \Z[W_1,W_2,p_1,r,b_1,b_2,b_3]/I\,,\]
where $\deg(r)=2$, $\deg(b_i)=i+3$ for $1\leq i\leq 3$, $b_1=\beta(s)$, $b_2=\beta(w_1 s)$, $b_3=\beta(w_2 s)$ and where $I$ is the ideal generated by $W_2^2-p_1W_1$, $r^2-4p_1$, $b_2 p_1-b_3 W_2$, $b_2 W_2-b_3 W_1$, $2W_i$, $rW_i$ and $b_1 W_i$ for $i=1,2$ as well as $2b_i$, $rb_i$ and $b_ib_j$ for $1\leq i,j\leq 3$.

\item The action of the Steenrod algebra on $H^\ast(\Bcom O(2);\F_2)$ is determined by its action on $H^\ast(BO(2);\F_2)$ and the total Steenrod squares
\begin{alignat*}{2}
\Sq(\bar{r}) & =\bar{r}\\
\Sq(s)       & =s+w_2 \bar{r}+w_1^2 s\, .
\end{alignat*}
\end{enumerate}
\end{theorem}

In all three cases $G=SU(2)$, $U(2)$, $O(2)$ our calculations show that the inclusion $\iota\colon \Bcom G\to BG$ induces an embedding $\iota^\ast\colon H^\ast(BG;\Z)\hookrightarrow H^\ast(\Bcom G;\Z)$. In fact, in Remark \ref{rem:so3injective} we argue that this is also true for $G=SO(3)$, without calculating the cohomology of $\Bcom SO(3)$ in this paper. This brings up the following question:

\begin{question}
For which compact Lie groups $G$ does the inclusion $\iota\colon \Bcom G\rightarrow BG$ induce an injective map on cohomology? In particular, does it always happen on integral cohomology when $G$ is compact and connected?
\end{question}

One has to be at least a little cautious with this question. For example, it is also true that $\iota\colon \Bcom G\to BG$ induces an injection on cohomology with $\F_2$-coefficients for $G=SU(2)$, $U(2)$, $O(2)$; and again in Remark \ref{rem:so3injective} we show this is true as well for $G=SO(3)$ (without computing $H^\ast(\Bcom SO(3); \F_2)$). However in \cite[Remark 3.2.2]{VillarrealThesis} it is shown that for the finite quaternion group $Q_8$, the homomorphism $\iota^\ast \colon H^\ast(\Bcom Q_8; \F_2) \to H^\ast(BQ_8; \F_2)$ is \emph{not} injective.

In the last third of our paper we study the homotopy fiber $\hofib(\Bcom G\to BG)=\Ecom G$. For compact and connected Lie groups $G$ a reasonable intuition for $\Ecom G$ is to think of it as how far is $G$ from being abelian. The paper \cite{Ad1} proves that for a compact, connected Lie group $G$, $\Ecom G_\one$ is homotopy equivalent to a finite CW-complex (this direct corollary of \cite[Theorem 6.5]{Ad1} is mentioned in \cite[Remark 6.6]{Ad1}). In this paper we compute the homotopy type of the $\Ecom$'s for our example groups, namely $\Ecom G$ for $G=SU(2), U(2), O(2)$ and $\Ecom SO(3)_\one$. In fact, according to Proposition \ref{prop:equivEcom} to cover most of our cases we only need to consider $\Ecom SU(2)$, because $\Ecom SU(2)$, $\Ecom U(2)$ and $\Ecom SO(3)_\one$ are all homotopy equivalent. 

\begin{theorem}\label{thm:ecomsu2}
  There is a homotopy equivalence $\Ecom SU(2) \simeq S^4 \vee \Sigma^4 \RP^2$.
\end{theorem}

The remaining case is $\Ecom O(2)$, which we compute separately.

\begin{theorem} \label{thm:ecomo2}
There is a homotopy equivalence $\Ecom O(2)\simeq S^2\vee S^2\vee S^3$.
\end{theorem}

\paragraph*{Application to bundle theory: transitionally commutative bundle structures.} We say that a principal $G$-bundle over a space $X$ with classifying map $f\colon X\to BG$ has a \emph{transitionally commutative structure} if there is a lift of $f$ to $\Bcom G$ along the inclusion $\iota$. A choice of a transitionally commutative structure is a homotopy class in $[X,\Bcom G]$. The exact sequence
\begin{equation*} \label{eq:sespin}
1\to \pi_n(\Ecom G)\to\pi_n(\Bcom G)\xrightarrow{\iota_\#} \pi_n(BG)\to1 
\end{equation*}
proved in \cite{Ad1} implies that any bundle over a sphere admits a transitionally commutative structure. It also follows that $\ker \iota_\#\cong\pi_n(\Ecom G)$, which is the set of transitionally commutative structures on the trivial bundle over $S^n$. In Remarks \ref{rem:pi_nBcomO2} and \ref{rem:pi_nBcomSU2} we compute some low dimensional homotopy groups of $\Bcom O(2)$ and $\Bcom SU(2)$. It is very surprising that even in the simplest case of $G=O(2)$, where $\pi_n(BO(2))=0$ for $n>2$, the trivial bundle over $S^n$ has so many distinct transitionally commutative structures. Namely, for $n>2$ each element of the groups listed in Remark \ref{rem:pi_nBcomO2} represents a transitionally commutative structure on the trivial $O(2)$-bundle over $S^n$!

On the other hand, it seems considerably harder to find examples of bundles, which do not admit any transitionally commutative structure. Predictably, these examples should arise from universal bundles, and indeed it was proved in \cite{Grthesis} that for a large class of connected Lie groups $G$ the universal bundle $EG\rightarrow BG$ cannot be made commutative. With the results of the present paper, however, we can give the first example of a vector bundle defined on a complex of dimension as small as four, which does not admit a transitionally commutative structure: the tautological bundle $\gamma_{2,4}\to \Gr_2(\R^4)$ over the Grassmannian of $2$-planes in $\R^4$.

\begin{corollary} \label{cor:tautological}
For any $n\geq 4$ the tautological $2$-plane bundle $\gamma_{2,n}\to \Gr_2(\R^n)$ does not admit a transitionally commutative structure. In particular, the universal bundle over $\Gr_2(\R^\infty)=BO(2)$ does not have a transitionally commutative structure.
\end{corollary}

\section{A homotopy pushout square for $\Bcom G$}\label{sec:poBcomG1}

When $G$ is compact and connected, a very useful tool towards the
understanding of $\Bcom G$ are the maximal tori $T\subset G$. Let
$\Hom(\Z^n,G)_\one$ denote the connected component of the trivial
homomorphism $\one \colon \Z^n\to G$ in $\Hom(\Z^n,G)$. By \cite[Lemma
4.2]{Baird}, the map
\[
\tilde{\phi_n}\colon G\times T^n\to \Hom(\Z^n,G)_\one
\]
given by
$(g,\vec{t})\mapsto g\vec{t}g^{-1}$ is surjective. Moreover,
$\tilde{\phi_n}$ is invariant under the diagonal action of $N(T)$ the
normalizer of $T$, which is given by right translation on $G$ and by
conjugation on $T^n$. Thus the induced map
\[\phi_n\colon G/T\times_WT^n=G\times_{N(T)}T^n\to\Hom(\Z^n,G)_\one\]
is also surjective. Since conjugation by elements of $G$ is a group
homomorphism, it is compatible with the simplicial structure.
Therefore, after taking geometric realizations, we get a map
\[\phi\colon G/T\times_{W}BT\to\Bcom G_\one,\]
where $\Bcom G_\one:=|\Hom(\Z^*,G)_\one|$. This map, which we'll call
the \emph{conjugation map}, is the key to computations for
$\Bcom G_\one$. For $G$ compact and connected, it is a rational
homology isomorphism (see the first page of \cite[Section 7]{Ad1}). In
many interesting cases, $\Bcom G_\one=\Bcom G$; for example, in
\cite{Ad0} it is shown that whenever any abelian subgroup of $G$ is contained in
a path-connected abelian subgroup, $\Hom(\Z^n,G)$ is connected. This
is the case for $G=U(k),SU(k)$ or $Sp(k)$, but not for $SO(3)$.

\medskip

The basis for all of our calculations are certain homotopy pushout
squares for $\Bcom G$ involving the conjugation map.

\begin{lemma}\label{lem:poBcomG1}
  Let $G$ denote $SO(3),SU(2)$ or $U(2)$, $T\subset G$ a maximal
  torus, $N=N(T)$, and $Z$ the center of $G$. Then we have a homotopy
  pushout square
  \begin{equation}\label{dgr:Bcomps}
    \xymatrix{G/N\times BZ\ar[r]^-{i}\ar[d]_{\pi_2}&G/T\times_{W}BT\ar[d]^-{\phi}\\
      BZ\ar@{^(->}[r]^-{j}&\Bcom G_\one,
    }
  \end{equation}
  where $i$ is the map induced by the inclusion $Z\hookrightarrow T$
  and $\pi_2$ is the projection. (As remarked above, the $\one$ only
  makes a difference for $SO(3)$.)
\end{lemma}

This lemma is a corollary of Theorem 6.3 in \cite{Ad1}, which
expresses $\Bcom G_\one$ for a compact connected Lie group as a
homotopy colimit. In the case of $SO(3), SU(2)$ or $U(2)$, their
diagram reduces to precisely the homotopy pushout above. (For $SU(2)$
this is Example 6.4 in \cite{Ad1}, and it is easy to check for the
other two cases as well.) We give an independent proof here, partly
because our approach is more direct and also establishes a homotopy
pushout for $\Hom(\Z^n, G)_\one$ (see Corollary \ref{cor:poHomG1}
below), but mostly because we also need a similar (but not completely
analogous!) pushout for $O(2)$, which is not covered by their theorem.

\begin{lemma} \label{lem:poBcomO2}
  There is a homotopy pushout
  \begin{equation} \label{dgr:poBcomO2}
    \xymatrix{
      O(2)/N(D_4)\times BZ \ar[r]^-{i} \ar[d]^q
      & O(2)\times_{N(D_4)}BD_4 \ar[d]^-{\phi} \\
      BSO(2) \ar[r]^-{j} & \Bcom O(2),
    }
  \end{equation}
  where $Z = \{\pm I\}$ is the center of $O(2)$, $D_4$ is the dihedral
  subgroup of $O(2)$ generated by two reflections in perpendicular
  lines, $i$ is the map induced by the inclusion
  $Z\hookrightarrow D_4$ and $q$ is the composite of the projection
  onto $BZ$ and the map induced by $Z\hookrightarrow SO(2)$, and
  $\phi$ is a conjugation map for $D_4 \le O(2)$ (defined analogously
  to the one above for $T \le G$).
\end{lemma}

We will deduce both from the following technical lemma.

\begin{lemma}\label{lem:abspo}
  Let $G$ be a compact Hausdorff group, $Z$ its center, $A$ be a
  closed abelian subgroup of $G$ with normalizer $N$, and $H$ be a
  closed normal abelian subgroup of $G$ such that (a) any two
  different conjugates of $A$ have intersection $Z$, and (b) $A$
  intersects $H$ in $Z$. Then we have a square which is both a pushout
  and a homotopy pushout:
  \begin{equation}
    \label{dgr:abspo}
    \xymatrix{G/N\times Z^n\ar[r]^-{i_n}\ar[d]_-{q_n} &
      G \times_{N}A^n\ar[d]^-{\phi_n}\\
      H^n\ar@{^(->}[r]^-{j_n}&\Hom(\Z^n,G)_?,}
  \end{equation}
  where $i_n$ is induced by the inclusion $Z \hookrightarrow A$, $q_n$
  is the composite of the projection onto $Z^n$ and the map induced by
  $Z^n \hookrightarrow H^n$, $\phi_n$ is the conjugation map for
  $A \le G$ and $\Hom(\Z^n,G)_?$ denotes the image of the map
  $(j_n, \phi_n) \colon H^n \sqcup G \times_{N}A^n \to \Hom(\Z^n, G)$.
\end{lemma}

In the cases we will use this lemma, we'll have either
$\Hom(\Z^n,G)_? = \Hom(\Z^n,G)$ or $\Hom(\Z^n,G)_\one$, which explains
the odd notation for the image.

\begin{proof}
  It is enough to prove the square is a pushout: the diagonal action
  of $N$ on $G\times A$ induces a $N$-CW-complex structure such that
  $G\times Z^n$ is an $N$-subcomplex, so passing to $N$-orbits we see
  that $i_n$ is a cofibration.

  The map $(j_n, \phi_n)$ is a surjective continuous map between
  compact Hausdorff spaces, and as such is closed and therefore a
  quotient map. This quotient can also be expressed as the quotient by
  an equivalence relation, namely, by the relation $\sim_H$ of having
  the same image under $(j_n, \phi_n)$.
  
  Now, the actual pushout, say $P$, is also a quotient of
  $H^n \sqcup G\times_N A^n$, namely it is the quotient by the
  equivalence relation \emph{generated} by
  $q_n(\zeta) \sim_P i_n(\zeta)$ for any $\zeta \in G/N\times Z^n$.
  Thus, to see that (\ref{dgr:abspo}) is a pushout, it is enough to
  show the equivalence relation $\sim_P$ is the same as $\sim_H$
  above.

  The fact that (\ref{dgr:abspo}) commutes shows that
  $q_n(\zeta) \sim_H i_n(\zeta)$ for any $\zeta \in G/N\times Z^n$.
  Since $\sim_P$ was defined to be generated by those pairs, this
  shows that the relation $\sim_P$ implies the relation $\sim_H$. Now
  we must show that any points related by $\sim_H$ are also related by
  $\sim_P$. There are three cases for a pair of points
  $\alpha, \beta \in H^n \sqcup G\times_N A^n$ related by $\sim_H$:

\begin{enumerate}
\item\emph{Both points are from $H^n$:} Since $j_n$ is injective, we
  have $\alpha = \beta$ and thus $\alpha \sim_P \beta$.
\item\emph{Both points are from $G\times_N A^n$:} Say
  $\alpha = [g, \vec{u}]$ and $\beta = [h, \vec{v}]$ with
  $\phi_n(\alpha) = \phi_n(\beta)$, which means
  $g \vec{u} g^{-1} = h \vec{v} h^{-1}$. Here we have again two cases:
  the conjugates $gAg^{-1}$ and $hAh^{-1}$ can either be the same
  or different. If they are the same, then $g^{-1}h \in N(A)$ and is
  witness to the fact that $[g, \vec{u}] = [h, \vec{v}]$ in
  $G\times_N A^n$, so again $\alpha = \beta$.

  If the conjugates are different, because their intersection is just
  the center $Z$, it follows that $\vec{u} = \vec{v} \in Z^n$. In this
  case, letting $\zeta_\alpha := (gN(A), \vec{u})$ and
  $\zeta_\beta := (h N(A), \vec{v})$ we see that
  \[\alpha = i_n(\zeta_\alpha) \sim_P q_n(\zeta_\alpha) = \vec{u} =
    \vec{v} = q_n(\zeta_\beta) \sim_P i_n(\zeta_\beta) = \beta,\] so
  that $\alpha \sim_p \beta$, as desired.
  
\item\emph{One point is from $H^n$ and one from $G\times_N A^n$:} Say
  $\alpha = \vec{s} \in H^n$ and $\beta = [g, \vec{u}]$ with
  $j_n(\alpha) = \phi_n(\beta)$. Since the intersection of $H$ with
  any conjugate of $A$ is the center, it follows that
  $\vec{s} = \vec{u} \in Z^n$ and thus letting
  $\zeta := (g N(A), \vec{s})$ we have
  $\alpha = q_n(\zeta) \sim_P i_n(\zeta) = \beta$.
\end{enumerate}
\end{proof}

Now we specialize this lemma to the groups we will study.

\begin{corollary}\label{cor:poHomG1}
  Let $G$ denote $SO(3),SU(2)$ or $U(2)$, $T\subset G$ a maximal
  torus, $N=N(T)$, and $Z$ the center of $G$. The following square is
  both a pushout and a homotopy pushout:
  \begin{equation}\label{dgr:poHomG1}
    \xymatrix{G/N\times Z^n\ar[r]^-{j_n}\ar[d]_-{\pi_2}&G/T\times_{W}T^n\ar[d]^-{\phi_n}\\
      Z^n\ar@{^(->}[r]&\Hom(\Z^n,G)_\one.}
  \end{equation}
  The $\one$ modifier is only necessary for $SO(3)$: for
  $G=SU(2)$ or $U(2)$, $\Hom(\Z^n,G)_\one = \Hom(\Z^n,G)$. 
\end{corollary}

\begin{proof}
  We apply Lemma \ref{lem:abspo} with $A=T$ and $H=Z$. Taking $H=Z$
  automatically ensures the second hypothesis, that $H$ intersects any
  conjugate of $A$ in $Z$. We must also show that the intersection of
  any two different maximal tori in $G$ is the center $Z$.

  For $G=SU(2)$, a short proof can be seen in \cite[Example
  2.22]{Villarreal}. The double cover
  $SU(2)\to SU(2)/\{\pm I\}\cong SO(3)$ implies that any two different
  maximal tori in $SO(3)$ intersect at the identity matrix, which is
  the center of $SO(3)$. Consider the 2-fold covering map
  $p\colon S^1\times SU(2)\to U(2)$ (we think of $S^1$ as the scalar
  matrices in $U(2)$ which is the center in this case). Let
  $T_2\subset U(2)$ denote the maximal torus consisting of the
  diagonal matrices in $U(2)$. Then $T_2$ is the image of $S^1\times T$
  under $p$, where $T\subset SU(2)$ is the subspace of diagonal
  matrices in $SU(2)$. Notice that $U(2)$ acts on $S^1\times SU(2)$,
  trivially on $S^1$ and by conjugation on $SU(2)$, making $p$ a
  $U(2)$-equivariant map. Let $g\in U(2)$ such that
  $T_2\ne g^{-1}T_2g$ and suppose $x\in T_2\cap g^{-1}T_2g$. Let
  $(\lambda, \bar{x})$ be an element in $p^{-1}(x)$. Then
  $\bar{x}\in T\cap h^{-1}Th=\{\pm I\}$ (where
  $h\sqrt{\det(g)} =g$) which implies that $x=\pm\lambda I$, a scalar
  matrix. Since any conjugate of $T$ contains the scalar matrices,
  $T_2\cap g^{-1}T_2g=S^1$.

  Finally, recall that \cite[Lemma 4.2]{Baird} says the map
  $\tilde{\phi_n}\colon G\times T^n\to \Hom(\Z^n,G)_\one$ is
  surjective on its own, without any help from $j_n$.
\end{proof}

\begin{corollary}\label{cor:poHomO2}
  The following square is both a pushout and a homotopy pushout:  
  \begin{equation} \label{dgr:poHomO2}
    \xymatrix{
      O(2)/N(D_4)\times Z^n \ar[r]^-{i_n} \ar[d]^-{q_n} &
      O(2)\times_{N(D_4)}D_4^n \ar[d]^-{\phi_n} \\
      SO(2)^n \ar[r]^-{j_n} & \Hom(\Z^n,O(2)).
    }
  \end{equation}
\end{corollary}

\begin{proof}
  We apply Lemma \ref{lem:abspo} with $G=O(2)$, $A=D_4$ and $H=SO(2)$.
  Any conjugate of $A = D_4$ consists of the elements of $Z$ and two
  reflections in perpendicular lines. This shows that the hypothesis
  of Lemma \ref{lem:abspo} are satisfied: any two conjugates of $A$
  have intersection $Z$, and any conjugate of $A$ intersects $H$ in
  precisely $Z$.

  It is an elementary observation that (1) a rotation $R\in O(2)$ and
  a reflection $r\in O(2)$ commute if and only if $R\in Z$ and (2) two
  reflections $r,r'\in O(2)$ commute if and only if they are
  reflections across perpendicular axes. Therefore, any $n$-tuple of
  commuting elements in $O(2)$ is either an $n$-tuple of elements in
  $SO(2)$, or it is conjugate to an $n$-tuple of elements in $D_4$.
  This proves the map
  \[(j_n, \phi_n) \colon SO(2)^n \sqcup O(2)\times_{N(D_4)}D_4^n \to
    \Hom(\Z^n,O(2))\] is surjective, so that in this case the image of
  $(j_n, \phi_n)$ is all of $\Hom(\Z^n, O(2))$.
\end{proof}

Now Lemmas \ref{lem:poBcomG1} and \ref{lem:poBcomO2} follow from
Corollaries \ref{cor:poHomG1} and \ref{cor:poHomO2}.

\begin{proof}[Proof of Lemmas \textup{\ref{lem:poBcomG1}} and \textup{\ref{lem:poBcomO2}}]
  One can easily check that the maps in diagrams (\ref{dgr:poHomG1})
  and (\ref{dgr:poHomO2}) commute with the simplicial structure maps,
  making those diagrams pushouts and homotopy pushouts of simplicial
  spaces. Since geometric realization commutes with both colimits and
  homotopy colimits, the conclusion follows.
\end{proof}

\section{Cohomology of $\Bcom SO(3)_\one$}

In this section we compute the integral and $\F_2$-cohomology rings of $\Bcom SO(3)_\one$, as well as the action of the Steenrod algebra on the $\F_2$-cohomology. We begin by recalling the cohomology of $BSO(3)$. We have
\[H^*(BSO(3);\Z)\cong \Z[\chi,p_1]/(2\chi),\]
where $p_1$ is the first Pontryagin class and $\chi$ is the Euler class of the universal oriented bundle over $BSO(3)$ (in degree 3) and
\[H^*(BSO(3);\F_2)\cong \F_2[w_2,w_3],\]
where $w_i$ is the $i$-th Stiefel-Whitney class. Consider the inclusion $\iota\colon \Bcom SO(3)_\one\to BSO(3)$. Abusing notation, we denote the pullback of $p_1$ and $w_i$ under $\iota$ by the same letters.

\begin{theorem}\label{thm:crso3_1}\leavevmode
\begin{enumerate}
\item\label{thm:crso3_1-hz} There is an isomorphism of graded rings
\[H^*(\Bcom SO(3)_\one;\Z)\cong\Z[p_1,w,y_1]/(2w,y_1^2,wy_1,w^3),\]
where $w=\iota^*(\chi)$ and $\deg y_1=4$.
\item There is an isomorphism of graded rings
\[H^*(\Bcom SO(3)_\one;\F_2)\cong\F[w_2,\bar{w},\bar{y}_1]/(\bar{y}_1^2,\bar{w}\bar{y}_1,\bar{w}^2+w_2\bar{y}_1),\]
where $\bar{w}$ and $\bar{y}_1$ are the reduction mod 2 of $w$ and $y_1$ in part \ref{thm:crso3_1-hz}, respectively, and $\iota^*(w_3)=\bar{w}$.
\item The action of the Steenrod algebra on $H^*(\Bcom SO(3)_\one;\F_2)$ is determined by the total Steenrod squares
\begin{alignat*}{2}
&\Sq(w_2)&&=w_2+\bar{w}\\
&\Sq(\bar{w})&&=\bar{w}+w_2\bar{w}+w_2\bar{y_1}\\
&\Sq(\bar{y}_1)&&=\bar{y}_1+w_2\bar{y}_1.
\end{alignat*}
\end{enumerate}
\end{theorem}

The proof of this theorem occupies the rest of this section. We will tackle each part in turn.

\subsection{The integral cohomology ring}

As we said before, the center of $SO(3)$ is trivial. A maximal torus is the canonically embedded $SO(2)\subset SO(3)$ with normalizer $O(2)\subset SO(3)$. Thinking of the elements of $SO(2)$ as 2 by 2 matrices, the action of $W=O(2)/SO(2)=\Z/2$ on $SO(2)$ is given by swapping columns and then swapping rows. That is, if a rotation in $SO(2)$ is defined by an angle $\theta$, the generating element in $O(2)/SO(2)$ acts by negating the angle of the rotation. This corresponds to complex conjugation in $S^1\cong SO(2)$. The homotopy pushout of Lemma \ref{lem:poBcomG1} for $G=SO(3)$ then reads

\begin{align}
\xymatrix{SO(3)/O(2)\ar[r]^-{i}\ar[d]&SO(3)/SO(2)\times_{W}BSO(2)\ar[d]^-{\phi}\\ {*}\ar[r]& \Bcom SO(3)_\one}\label{dgr:poso3}
\end{align}

To compute the cohomology ring of $\Bcom SO(3)_\one$, it is enough to give a presentation of the cohomology ring of $SO(3)/SO(2)\times_WBSO(2)$, and then describe the top horizontal map of (\ref{dgr:poso3}) in cohomology.

 Consider the inclusion $\iota\colon \Bcom SO(3)_\one\to BSO(3)$. Define $\tilde{p}_1:=(\iota\phi)^*(p_1)$. Let $U\in H^2(\RP^2;\Z)$ be the generator, and $\pi\colon SO(3)/SO(2)\times_WBSO(2)\to SO(3)/SO(2)/W\cong \RP^2$ be the projection. Abusing notation, write $U$ for $\pi^*(U)$. 

\begin{lemma}\label{lem:icr1}
There is an isomorphism of graded rings
\[H^*(SO(3)/SO(2)\times_WBSO(2);\Z)\cong\Z[U,x,\tilde{p}_1,y]/(2U,2x,U^2,xy,Uy,Ux,y^2,x^2-\tilde{p}_1U) ,\]
where $\deg x=3$ and $\deg y=4$. 
\end{lemma}

\begin{proof}
Identifying $SO(3)/SO(2)$ with $S^2$, and writing $SO(2)=T$, we compute the integral cohomology of $S^2\times_WBT$ using the Serre spectral sequence associated to the fibration $BT\to S^2\times_WBT\to S^2/W\cong\RP^2$. To do this we need to determine the action of $\pi_1(\RP^2)$ on $H^*(BT;\Z)\cong \Z[a]$, where $a\in H^2(BT;\Z)$ is a generator. Lifting a representative of the generator of $\pi_1(\RP^2)\cong W$ under the fibration $S^2\times_WBT\to \RP^2$ shows that the monodromy action on $BT$ is the same as the action of the Weyl group, which is induced by complex conjugation on $T$. On the cohomology group $H^1(T;\Z)$ this is multiplication by $-1$, and the same holds for $H^2(\Sigma T;\Z)$ via the suspension isomorphism. This implies that the action on $H^*(BT;\Z)$ is generated by $a^n\mapsto (-a)^n$. We conclude that $H^q(BT;\Z)$ is the trivial representation $\Z$ when $q\equiv 0$ (mod $4$) and the sign representation $\Z_w$ when $q\equiv 2$ (mod $4$). 

Now we can compute the $E_2$-page of the Serre spectral sequence in cohomology. We have that $E_2^{p,q}=H^p(\RP^2;\underline{H}^q(BT;\Z))$, where 
\[\underline{H}^q(BT;\Z)=\left\{
\begin{matrix}
\Z&q\equiv 0\;(\text{mod 4})\\
\Z_w&q\equiv 2\;(\text{mod 4})\\
0&\text{else.}
\end{matrix}\right.\]
Since $H^*(\RP^2;\Z)\cong \Z[U]/(2U,U^2)$, where $\deg U=2$, and $H^p(\RP^2;\Z_w)=H_{2-p}(\RP^2;\Z)$ by Poincar\'e Duality, we get that the $E_2$-page looks like
\[
\begin{tikzpicture}
  \matrix (m) [matrix of math nodes,
    nodes in empty cells,nodes={minimum width=5ex,
    minimum height=5ex,outer sep=-5pt},
    column sep=1ex,row sep=-2ex]{
            8&b^2\Z  & 0      &  U_4\Z/2 & \\             
            7&0   &  0     & 0     & \\      
            6&0   & x_2\Z/2   &  y_2\Z   & \\ 
            5&0   &  0     & 0     & \\         
            4&b\Z  & 0      &  U_2\Z/2 & \\             
            3&0   &  0     & 0     & \\
           2 &0   & x\Z/2   &  y\Z   & \\
           1 &0   &  0     & 0     & \\
        0 &\Z  & 0      &  U\Z/2 &  \\
\quad\strut         &0  & 1      &  2& \strut \\};
  %\draw[-stealth] (m-3-3.north west) -- (m-2-2.south east);
\draw[thick] (m-1-1.east) -- (m-10-1.east) ;
\draw[thick] (m-10-1.north) -- (m-10-5.north) ;
\end{tikzpicture}
\]
and for $p>2$, $E_{2}^{p,q}= 0$. From this it follows that $xy=Uy=Ux=y^2=x^3=U^2=0$. We claim that $b^{n}y=y_{2n},b^{n}U=U_{2n},b^{n}x=x_{2n}$, for all $n\geq1$, and the class $x^2$ in ${H}^2(\RP^2;\underline{H}^4(BT;\Z))$ equals $U_2$. The cup product in $H^*(BT;\Z)$ induces a product in the local coefficient system $\underline{H}^{4n}(BT;\Z)\times \underline{H}^s(BT;\Z)\to \underline{H}^{4n+s}(BT;\Z)$, given by the bilinear maps $\Z\times \Z_w\to \Z_w$ or $\Z\times \Z\to \Z$. In both cases, choosing the positive generator of $\Z$ in the first argument, we get the identity map. Thus, $b^{n}U=U_{2n}$ and $b^{n}y=y_{2n}$.  Consider the diagram
\[\xymatrix{{H}^1(\RP^2;\underline{H}^2(BT;\Z))\times{H}^{1}(\RP^2;\underline{H}^2(BT;\Z))\ar[r]\ar[d]&{H}^{2}(\RP^2;\underline{H}^{4}(BT;\Z))\ar[d]\\
{H}^1(\RP^2;\Z/2)\times{H}^1(\RP^2;\Z/2)\ar[r]^-{\smile}&{H}^{2}(\RP^2;\Z/2)}
\]
where the vertical arrows are induced by the reduction modulo 2
homomorphisms $\Z_w\to\Z/2$ and $\Z\to\Z/2$, respectively. These are
all isomorphisms by Lemma \ref{lem:rp2coh}. Let
$\widehat{x}\in {H}^1(\RP^2;\Z/2)$ be the reduction modulo 2 of $x$ in
${H}^1(\RP^2;\underline{H}^2(BT;\Z))$. Notice that the bottom arrow maps the class $\widehat{x}\times\widehat{x}$ to the generator in ${H}^2(\RP^2;\Z/2)$, that is, $\widehat{x}^2$. By commutativity of the diagram, $x^2$ is also a generator and thus non-zero. Now the relation $x^2-bU=0$ follows.

All further differentials are zero and the spectral sequence collapses at the $E_2$-page. Therefore, the $E_\infty$-page has the multiplicative structure 
\begin{align*}
E_\infty^{*,*}\cong\Z[U,x,b,y]/(2U,2x,U^2,xy,Uy,Ux,y^2,x^2-bU)\label{eq:cr1}
\end{align*}
where $\deg U=2$, $\deg x=3$ and $\deg b=\deg y=4$. Notice there are no additive extensions, since the only groups appearing as quotients are the $\Z$'s in the first column. Next we will show that all relations between the generators still hold in $H^*(S^2\times_W BT;\Z)$, but to see this, we need to make a choice of the class $b$ in $H^4(S^2\times_WBT;\Z)$, that corresponds to a representative of the class $a^2$ in $H^4(BT;\Z)$. Consider the composite
\begin{align}
\xymatrix{BT\ar[r] & S^2\times_{W}BT\ar[r]^{\iota\phi} &BSO(3)}\
\end{align}
The first Pontryagin class $p_1\in H^4(BSO(3);\Z)$ maps to $a^2$, and so we choose $b:=(\iota\phi)^*(p_1)$. 

Any ``zero'' relation in $E_\infty^{*,*}$ takes place in $F_i^{p+q}$, the $i$-th term of the associated filtration of $H^{p+q}(S^2\times_W BT;\Z)$ with $i\geq 2$. For $i\geq 3$, they are all zero and $F^{p+q}_2$ is a subgroup. Therefore we have an isomorphism of graded rings $H^*(S^2\times_W BT;\Z)\cong E_{\infty}^{*,*}$.
\end{proof}

\begin{proposition}\label{prop:icrso3}
  There is an isomorphism of graded rings \[H^*(\Bcom SO(3)_\one;\Z)\cong\Z[p_1,w,y_1]/(2w,y_1^2,wy_1,w^3),\] where $w$ and $p_1$ are the pullbacks of $\chi \in H^3(BSO(3); \Z)$ and $p_1\in H^4(BSO(3);\Z)$ under the inclusion $\Bcom SO(3)_\one\hookrightarrow BSO(3)$, and $y_1$ satisfies $\phi^*(y_1)=y$ with $y$ as in Lemma \ref{lem:icr1}.
\end{proposition}

\begin{proof}
Notice that the composition
\begin{equation} \label{eq:retract}
\RP^2\stackrel{j}{\longrightarrow}S^2\times_WBSO(2)\stackrel{\pi}{\longrightarrow} \RP^2
\end{equation}
is the identity map, where $j$ is the inclusion at the basepoint of $BSO(2)$ and $\pi$ is the projection map in our fibration. By inspection of the above spectral sequence, we see that the class $U$ in Lemma \ref{lem:icr1} maps under $j^*$ to the generator of $H^2(\RP^2;\Z)$. Applying Mayer-Vietoris to (\ref{dgr:poso3}), we see that $\tilde{H}^*(\Bcom SO(3)_\one;\Z)=\ker j^*$. An easy check shows that the kernel of $j^*$ is the ideal generated by $x$, $y$ and $\tilde{p}_1$. Since $Ux = Uy = 0$ and $\tilde{p}_1 U = x^2$, we see that the ideal generated by $x$, $y$ and $\tilde{p}_1$ is equal to the positive degree part of the \emph{subring} of $H^*(S^2 \times_W BSO(2);\Z)$ generated by $x$, $y$ and $\tilde{p}_1$, and thus the cohomology ring $H^*(\Bcom SO(3)_\one; \Z)$ is isomorphic to said subring.

Now we try to identify the classes in $H^*(\Bcom SO(3)_\one;\Z)$ that map to $x$, $y$ and $\tilde{p}_1$ under $\phi^*$. First let's deal with $x$. Using the fact that $\Ecom SO(3)_\one$ is 3-connected and Serre's exact sequence for the homotopy fibration $\Ecom SO(3)_\one\to \Bcom SO(3)_\one\to BSO(3)$, it is readily seen that the pullback of $\chi$ under the inclusion $\Bcom SO(3)_\one\to BSO(3)$ is non-zero, so that we must have $\phi^*(w) = x$. As for $\tilde{p}_1$, it was defined as $(\iota \phi)^*(p_1)$, and therefore $\phi^*(p_1)=\tilde{p}_1$ if we regard $p_1$ as a class in $H^\ast(\Bcom SO(3)_\one;\Z)$. Finally, we have nothing to say about the pre-image of $y$ which we simply call $y_1$.

The relations claimed for $p_1$, $w$ and $y$ follow easily from the presentation in Lemma \ref{lem:icr1}. The only non-immediate one is $w^3 = 0$, which follows from $x^3 = x(x^2-\tilde{p}_1U) + (Ux)\tilde{p}_1$. To show that no other relations are needed, in Section \ref{sec:singular3} we used {\sc Singular} to verify that the kernel of the ring homomorphism $\Z[p_1, w, y_1] \to H^*(S^2 \times_W BSO(2); \Z)$, $(p_1, w, y_1) \mapsto (\tilde{p}_1, x, y)$ is indeed $(2w,y_1^2,wy_1,w^3)$.
\end{proof}

\begin{remark}
When $G$ is a compact Lie group, the inclusion $BT \to BG$ is injective in rational cohomology. The fact that this inclusion factors through $\Bcom G_\one\to BG$ implies that the induced map $H^*(BG;\Q)\to H^*(\Bcom G_\one;\Q)$ is also injective. Proposition \ref{prop:icrso3} shows that this is not true with integral coefficients for $G=SO(3)$. Indeed, $\chi^3$ is pulled back to $w^3=0$.  
\end{remark}

\begin{remark} \label{rem:so3injective}
  In light of the above remark, it is natural to ask about the differences between $\Bcom SO(3)_\one$ and the ``full'' classifying space for commutativity $\Bcom SO(3)$. The main difference between these two spaces is that $BD_4$ appears as a natural subcomplex of $\Bcom SO(3)$, because of the subgroup of $SO(3)$ isomorphic to $D_4$ consisting of diagonal matrices with $\pm 1$ diagonal entries (and determinant $1$), which is not contained in any maximal torus and thus does not contribute to $\Bcom SO(3)_\one$. With help from that $BD_4$ we can show that the inclusion $\iota\colon \Bcom SO(3)\to BSO(3)$ induces an injective homomorphism $H^*(BSO(3))\to H^*(\Bcom SO(3))$ with either $\Z$ or $\F_2$ coefficients.

Consider the inclusions
\begin{align}
\xymatrix{BSO(2)\ar[r] \ar@/_1pc/[rd]_{j_1} & \Bcom SO(3)\ar[d]^{\iota}\\&BSO(3)}&&\xymatrix{BD_4\ar[r] \ar@/_1pc/[rd]_{j_2} & \Bcom SO(3)\ar[d]^{\iota}\\&BSO(3)}
\end{align}
In integral cohomology, $\ker j_1^*=(\chi)$ and it can be shown that $\ker j_2^*=(2p_1)$. Since the above triangles show that $\ker \iota^* \subset \ker j_1^*$ and $\ker \iota^* \subset \ker j_2^*$, we must have $\ker \iota^* \subset \ker j_1^* \cap \ker j_2^* = 0$.

For $\F_2$-coefficients we only need the second diagram. We will show that $j_2^\ast$ is injective on $\F_2$-cohomology, which implies that $\iota^\ast$ must be as well. Let $x, y$ be the generators in degree $1$ of the cohomology ring of $BD_4 = B\Z/2 \times B\Z/2$. Thinking of $D_4$ as $\{\pm 1\} \times \{\pm 1\}$ for convenience, the inclusion $D_4 \to SO(3)$ is given by $(\epsilon_1, \epsilon_2) \mapsto \mathrm{diag}(\epsilon_1, \epsilon_2, \epsilon_1 \epsilon_2)$. Now, that same formula embeds $D_4$ in $O(1)^3$ and the induced map on cohomology,  $\F_2[t_1, t_2, t_3] \cong H^\ast(BO(1)^3; \F_2) \to H^\ast(BD_4; \F_2)$ is given by $t_1 \mapsto x$, $t_2 \mapsto y$ and $t_3 \mapsto x+y$. Recall that the $\F_2$-cohomology of $O(3)$ is the subring of symmetric polynomials inside the cohomology  of its subgroup $O(1)^3$ and the classes $w_2$ and $w_3$ map to $t_1 t_2 + t_1 t_3 + t_2 t_3$ and $t_1 t_2 t_3$, respectively. Therefore the inclusion $D_4 \to SO(3)$ induces the homomorphism on cohomology of classifying spaces given by \[w_2 \mapsto t_1 t_2 + t_1 t_3 + t_2 t_3 \mapsto xy + x(x+y) + y(x+y) = x^2 + xy + y^2\] and \[w_3 \mapsto t_1 t_2 t_3 \mapsto xy(x+y) = x^2y+xy^2.\]
This homomorphism $\F_2[w_2, w_3] \to \F_2[x,y]$ can easily be seen to be injective, for example, by computing the kernel in \textsc{Singular}.
\end{remark} 

\subsection{The $\F_2$-cohomology ring and Steenrod squares}

Let $\tilde{w}_2=(\iota\phi)^*(w_2)$, and write $H^*(\RP^2;\F_2)=\F_2[u]/(u^3)$.

\begin{lemma}\label{lem:f2crS2xBT}
There is an isomorphism of graded rings
\[H^*(S^2\times_WBSO(2);\F_2)\cong \F_2[\tilde{w}_2,u]/(u^3).\]
\end{lemma}

\begin{proof}
Since $\F_2$ has a unique automorphism, the monodromy action of $\pi_1(\RP^2)$ on $H^*(BSO(2);\F_2)$ must be trivial. The $E_2$-page of the Serre spectral sequence associated to $BSO(2)\to S^2\times_WBSO(2)\to \RP^2$ is then $H^*(BSO(2);\F_2)\otimes H^*(\RP^2;\F_2)\cong \F_2[w_2]\otimes \F_2[u]/(u^3)$ and all differentials are zero. In fact, the spectral sequence also collapses at the $E_2$-page. To lift the multiplicative structure from the $E_{\infty}$--page to the cohomology of $S^2\times_WBSO(2)$ we check that $\tilde{w}_2$ is a lift of $w_2$ in the $E_{\infty}$--page. Since $w_2$ is stable, the class $w_2\in H^2(BSO(3);\F_2)$ maps to $w_2\in H^2(BSO(2);\F_2)$ under the inclusion $BSO(2)\hookrightarrow BSO(3)$. The latter factors through the conjugation map $\phi$,
\[
\xymatrix{
BSO(2) \ar[r] & S^2\times_WBSO(2) \ar[r]^-{\iota\phi} & BSO(3)}
\]
and, therefore, $\tilde{w}_2$ is a lift of $w_2$ in $E_\infty^{0,2}$.
\end{proof}

\begin{proposition}\label{prop:f2crso3}
 Let $w_2\in H^2(\Bcom SO(3)_\one;\F_2)$ also denote the pullback of $w_2\in H^2(BSO(3);\F_2)$ under the inclusion $\Bcom SO(3)_\one\to BSO(3)$. Then
\[H^*(\Bcom SO(3)_\one;\F_2)\cong\F[w_2,\bar{w},\bar{y}_1]/(\bar{y}_1^2,\bar{w}\bar{y}_1,\bar{w}^2+w_2\bar{y}_1),\]
where $\bar{w},\bar{y}_1$ are the reduction mod 2 of $w,y_1$ in Proposition \ref{prop:icrso3}, respectively. 
\end{proposition}

\begin{proof}
Similarly to the computation with integer coefficients,here we also have that $\tilde{H}^*(\Bcom SO(3)_\one;\F_2)=\ker j^*$, where $j$ is the first map in (\ref{eq:retract}). In the case of $\F_2$--coefficients the kernel of $j^*$ is the ideal generated by $\tilde{w}_2$. We claim that $\tilde{w}_2u$ and $\tilde{w}_2u^2$ are the reduction mod 2 of the classes $x$ and $y$ in Lemma \ref{lem:icr1}, respectively.  We compare the $E_2$-pages of the integral and the mod 2 Serre spectral sequences via the mod 2 reduction homomorphism $H^2(BSO(2);\Z)\to H^2(BSO(2);\F_2)$. As a morphism of $\pi_1(\RP^2)$-representations, this is the non-trivial morphism $\Z_w \to \F_2$. From Lemma \ref{lem:rp2coh} we see that $x$ and $y$ indeed map to the generators $\tilde{w}_2u$ and $\tilde{w}_2u^2$ on the $E_2$-page, and thus this is also true on the $E_\infty$-page. Since the classes $x,y$ and $\tilde{w}_2u,\tilde{w}_2u^2$ are in honest subgroups of the cohomology groups of the corresponding total spaces, this continues to hold in these cohomology groups. The stated presentation follows from the naturality of reduction modulo $2$.
\end{proof}

\paragraph*{The Steenrod algebra action.} Recall that the Steenrod squares in $H^*(BSO(3);\F_2)\cong\F_2[w_2,w_3]$ are given by $\Sq^1(w_2)=w_3$, $\Sq^1(w_3)=0$ and $\Sq^2(w_3)=w_2w_3$. Now under the inclusion $\iota\colon \Bcom SO(3)_\one\to BSO(3)$ we have $\iota^*(w_3)=\bar{w}$. We can immediately conclude 
\begin{align*}
\Sq^1(w_2)=\bar{w}& &\Sq^1(\bar{w})&=0\\
&  &\Sq^2(\bar{w})&=\bar{w}w_2.
\end{align*}
To determine $\Sq^n(\bar{y}_1)$, recall from the proof of Proposition \ref{prop:f2crso3} that $\phi^*(\bar{y}_1)=\tilde{w}_2u^2$ and $\phi^*(\bar{w})=\tilde{w}_2u$. Since $\phi^*$ is injective, it is enough to calculate $\Sq^n(\tilde{w}_2u^2)$. By construction $(\iota\phi)^*(w_2)=\tilde{w}_2$, and hence 
\[\Sq^1(\tilde{w}_2)=(\iota\phi)^*(\Sq^1(w_2))=(\iota\phi)^*(w_3)=\phi^*(\bar{w})=\tilde{w}_2u\]
An easy check shows that $\Sq^1(\tilde{w}_2u^2)=0,\Sq^2(\tilde{w}_2u^2)=\tilde{w}_2^2u^2,
$ and $\Sq^3(\tilde{w}_2u^2)=0$. Thus,
\begin{align*}
\Sq^1(\bar{y}_1)&=0\\
\Sq^2(\bar{y}_1)&=w_2\bar{y}_1\\
\Sq^3(\bar{y}_1)&=0.
\end{align*}

\section{Cohomology of $\Bcom U(2)$ and $\Bcom SU(2)$}

The goal of this section is to compute the cohomology rings of the spaces  $\Bcom U(2)$ and $\Bcom SU(2)$ with both integral and $\F_2$-coefficients, as well as the action of the Steenrod algebra on the $\F_2$-cohomology. Our calculations for the two groups are inextricably intertwined, which is why we present them together. Our calculations of the ring structure and the Steenrod squares for $\F_2$-cohomology rings are also intertwined.

Recall that the integral cohomology rings of $BU(2)$ and $BSU(2)$ are described in terms of Chern classes, namely, $H^*(BU(2); \Z) \cong \Z[c_1, c_2]$, where $c_i$ has degree $2i$, and $H^*(BSU(2); \Z) \cong \Z[c_2]$. The inclusion $SU(2) \hookrightarrow U(2)$ does the obvious thing on cohomology: $c_1 \mapsto 0$, $c_2 \mapsto c_2$.

Our results for $SU(2)$ were already stated in Theorem \ref{thm:crsu2}. Here is the corresponding omnibus theorem for $U(2)$.

\begin{theorem}\label{thm:cru2}\leavevmode  
  \begin{enumerate}
  \item There is an isomorphism of graded rings
    \[H^*(\Bcom U(2);\Z) \cong
      \Z[c_1,c_2,y_1,y_2]/(2y_2-y_1c_1,y_1^2,y_1y_2,y_2^2)\] where
    $c_i\in H^{2i}(\Bcom U(2);\Z)$ and
    $y_i\in H^{2i+2}(\Bcom U(2);\Z)$. Moreover, $c_i$ is the pullback
    of the $i$-th Chern class along the inclusion
    $\Bcom U(2)\to BU(2)$.

  \item Let $\bar{z}$ denote reduction mod 2 of $z$.
    Then
    \[H^*(\Bcom U(2);\F_2) \cong \F_2[\bar{c}_1,\bar{c}_2,\bar{y}_1,\bar{y}_2]/(\bar{y}_1\bar{c}_1,\bar{y}_1^2,\bar{y}_1\bar{y}_2,\bar{y}_2^2).\]
  \item The Steenrod algebra action is given by the following total Steenrod squares:
    \begin{align*}
      \Sq(\bar{c}_1) & = \bar{c}_1 + \bar{c}_1^2 \\
      \Sq(\bar{c}_2) & = \bar{c}_2 + \bar{c}_1 \bar{c}_2 + \bar{c}_2^2\\
      \Sq(\bar{y}_1) & = \bar{y}_1 \\
      \Sq(\bar{y}_2) & = \bar{y}_2 + \bar{c}_2 \bar{y}_1 + \bar{c}_1^2 \bar{y}_2. \\
    \end{align*}
\end{enumerate}
\end{theorem}

Even though $x_2 \in H^6(\Bcom SU(2); \Z)$ is the image of $y_2 \in H^6(\Bcom U(2); \Z)$, we've changed the name to heighten awareness of the fact that $y_2$ is non-torsion while $x_2$ is 2-torsion.

\subsection{The integral cohomology rings}

The homotopy pushout for $G=U(2)$ (Lemma \ref{lem:poBcomG1}) reads 
\begin{align}
\xymatrix{\RP^2\times BS^1\ar[r]^-{i}\ar[d]_-{\pi_2}&S^2\times_{W}BT_2\ar[d]^-{\phi}\\ BS^1\ar[r]^-{j}&\Bcom U(2)}\label{dgr:pou2}
\end{align}

To compute the integral cohomology ring of $\Bcom U(2)$, we need to know what the map $i^*$ does in cohomology.

Let $\tilde{c}_i:=(\iota\phi)^*(c_i)$, where $c_i\in H^{2i}(BU(2);\Z)$
is the $i$-th Chern class.

\begin{samepage}
\begin{lemma}\label{lem:crS2xBT2}
There is an isomorphism of graded rings
\[H^*(S^2\times_WBT_2;\Z)\cong \Z[\tilde{c}_1,\tilde{c}_2,c,U]/(2U,U^2,c^2,\tilde{c}_1U,cU)\]
where $\deg{c}=4$ and $\deg U=2$.
\end{lemma}
\end{samepage}

\begin{proof}
We have a fibration $BT_2\to S^2\times_W BT_2\to \RP^2$. The action of the Weyl group $W$ on $T_2$ is by swapping the entries of the diagonal. Thus the induced action of $\pi_1(\RP^2)$ on $H^2(BT_2;\Z)=\langle a_1,a_2\rangle$ is generated by $(a_1,a_2)\mapsto (a_2,a_1)$, i.e. $\underline{H}^2(BT;\Z)$ is the regular representation $\Z[\Z/2]=\{m+\tau n\;|\;m,n\in\Z,\;\tau^2=1\}$. Since $H^q(BT_2;\Z)$ is zero for $q$ odd and is generated by the elements $a_1^{i_1}a_2^{i_2}$, where $i_1,i_2\geq 0$ and $i_1+i_2=\frac{q}{2}$ when $q$ is even,  we obtain
\[\underline{H}^q(BT_2;\Z)=\left\{
\begin{matrix}
\Z[\Z/2]^{n+1}&q=4n+2&n\geq0\\
\Z[\Z/2]^{n}\oplus\Z&q=4n&n\geq 0\\
0&\text{else.}
\end{matrix}\right.\] 
Recall that $H^p(\RP^2;\Z[\Z/2])\cong H^p(S^2;\Z)$. Thus the $E_2$-page of the Serre spectral sequence associated to the above fibration is given by \[E^{p,q}_2=H^p(\RP^2;\underline{H}^q(BT_2;\Z))\cong\left\{
\begin{matrix}
\Z^{n+1}&p=0,2&q=4n+2&n\geq0\\
\Z^{n+1}&p=0&q=4n&n\geq 0\\
\Z^{n}\oplus\Z/2&p=2&q=4n&n\geq 0\\
0&\text{else.}
\end{matrix}\right.\]

It is easy to determine the multiplicative structure of the $E^{*,*}_2$-page:

\begin{equation} \label{eq:ssbt2}
  \begin{minipage}{0.5\linewidth}
  \begin{tikzpicture}
  \matrix (m) [matrix of math nodes,
    nodes in empty cells,nodes={minimum width=5ex,
    minimum height=5ex,outer sep=-5pt},
    column sep=1ex,row sep=-2ex]{
4n+2          &\Z^{n+1}            & 0       &  \Z^{n+1}               & \\             
4n+1          &0                   &  0      & 0                       & \\      
           4n &\Z^{n+1}            &0        &  \Z^n\oplus\Z/2         & \\ 
       \vdots &\vdots              &  \vdots & \vdots                  & \\         
            4 &b_1^2\Z\oplus b_2\Z & 0       &  b_1c \Z\oplus b_2u\Z/2 & \\             
            3 &0                   &  0      & 0                       & \\
           2  &b_1\Z               & 0       &  c\Z                    & \\
           1  &0                   &  0      & 0                       & \\
        0     &\Z                  & 0       &  U\Z/2                  &  \\
\quad\strut   &0                   & 1       &  2                      & \\};
  %\draw[-stealth] (m-3-3.north west) -- (m-2-2.south east);
\draw[thick] ([xshift=7mm]m-1-1.north) -- ([xshift=7mm]m-10-1.south) ;
\draw[thick] ([yshift=1.5mm]m-10-1.north) -- ([yshift=1.5mm]m-10-5.north) ;
\end{tikzpicture}
\end{minipage}
\end{equation}

Here $c^2=U^2=b_1U=cU=0$. The only non-trivial products occur in the zero column, as well as between the zero and second column, where we have $b_2^n U$ for $n\geq 0$, each giving a copy of $\Z/2$, and $p(b_1,b_2) c$ (where $p$ is any polynomial in $b_1,b_2$) giving the free summands in the second column. These products are easy to determine using Lemma \ref{lem:rp2coh}. All differentials are zero, in fact $E_2^{p,q}\cong E_\infty^{p,q}$ and
\begin{align}
E^{*,*}_\infty\cong \Z[b_1,b_2,c,U]/(2U,U^2,c^2,b_1U,cU)\,,\label{eq:crS2xBT2}
\end{align}
where $\deg b_1=\deg U=2$ and $\deg b_2=\deg c=4$. Since $E_\infty^{p,q}=0$ for $p\geq 3$, the classes in the second column lift uniquely to $H^\ast(S^2\times_W BT;\Z)$ and any product involving these classes is uniquely determined. On the other hand, there may be multiplicative extensions for the products in the zero column. However, these products are determined by those in $H^\ast(BU(2);\Z)$ once we've shown that $\tilde{c}_1=(\iota\phi)^*(c_1)$ and $\tilde{c}_2=(\iota\phi)^*(c_2)$ represent $b_1$ and $b_2$, respectively, in the $E_\infty$--page. The latter follows from the factorization
\begin{align}
\xymatrix{BT_2\ar[r] & S^2\times_{W}BT_2\ar[r]^{\iota\phi} &BU(2)}
\end{align}
and the fact that the Chern classes $c_1$, $c_2$ pull back to the generators of the invariant ring $H^0(\RP^2;\underline{H}^*(BT_2;\Z))$ isomorphic to $H^*(BT_2;\Z)^W$. Note that any other representative of $b_1$, say $\tilde{c}_1+kU$, still satisfies $(\tilde{c}_1+kU)U=0$ so that the presentation of the ring (\ref{eq:crS2xBT2}) remains unaltered. Upon replacing $b_1$, $b_2$ by $\tilde{c}_1$, $\tilde{c}_2$, this gives the cohomology ring claimed in the lemma.
\end{proof}

Consider the trivial fibration $BS^1\to\RP^2\times BS^1\to\RP^2$ and let ${E^\prime}_r^{p,q}$ denote the associated Serre spectral sequence. The map $i$ induces a map of fibrations and hence a morphism of spectral sequences and $E_\infty$--pages $i_\infty\colon E_{\infty}^{p,q}\to {E^\prime}^{p,q}_\infty$. We use $i_\infty$ to describe the map \[i^*\colon H^*(S^2\times_WBT_2;\Z)\to H^*(\RP^2\times BS^1;\Z)\cong \Z[U,t]/(2U,U^2)\,,\]
where $\deg U=\deg t=2$.
\begin{lemma} \label{lem:jbt2}
In the presentation of Lemma \ref{lem:crS2xBT2}, $i^*(U)=U$, $i^*(\tilde{c}_1)=2t$, $i^*(c)=tU$, and $i^*(\tilde{c}_2)=t^2$.
\end{lemma}
\begin{proof}
The map of fibrations is the identity map on base spaces, thus $i^*(U)=i_\infty(U)=U$. We analyze the class  $i_\infty(c)\in H^2(\RP^2;H^2(BS^1;\Z))$. The map $i_\infty$ in bi-degree (2,2) is induced by the morphism of $\Z/2$-representations $\Z[\Z/2]\to \Z$ sending $m+\tau n\mapsto m+n$. The induced map $H^2(\RP^2;\Z[\Z/2])\to H^2(\RP^2;\Z)$ is shown to be reduction modulo 2 in Lemma \ref{lem:rp2coh}. Both spectral sequences are concentrated in degrees $p=0,2$, which means that $E_\infty^{2,2}$ and ${E^\prime}^{2,2}_\infty$ are actual subgroups of $H^4(S^2\times_WBT_2;\Z)$ and $H^4(\RP^2\times BS^1;\Z)$, respectively. Thus $i^*(c)=tU$. Consider the commutative diagram
\[
\xymatrix{
\RP^2\times BS^1\ar[r]^-{i}\ar[d]_-{\pi_2}&S^2\times_{W}BT_2\ar[r]^-{\phi}&\Bcom U(2)\ar[d]^-{\iota}\\ BS^1\ar[r]^-{j_1}&BT_2\ar[r]^-{j_2}&BU(2)\,,
}
\]
where $j_1$ is the diagonal inclusion. Suppose $a_1,a_2$ are the generators of degree 2 of the cohomology ring $H^*(BT_2;\Z)$. Recall that the inclusion $T_2\hookrightarrow U(2)$ induces an isomorphism $H^*(BU(2);\Z) \cong H^*(BT_2;\Z)^W$ given by $j_2^*(c_1)=a_1+a_2$ and $j_2^*(c_2)=a_1a_2$. Then $j_1^*(a_k)=t$ implies that that $(j_2\circ j_1\circ \pi_2)^*(c_1)=2t$ and $(j_2\circ j_1\circ \pi_2)^*(c_2)=t^2$. 
Now we have a complete description of $i^*$, since by the definition of $\tilde{c}_1$ and $\tilde{c}_2$ we have $i^*(\tilde{c}_1)=2t$ and $i^*(\tilde{c}_2)=t^2$. 
\end{proof}

\begin{lemma}\label{lem:icru2}
  The integral cohomology ring of $\Bcom U(2)$ is as claimed in Theorem \ref{thm:cru2}.
\end{lemma}

\begin{proof}
From the previous lemma it follows that $(\pi_2^\ast-i^\ast)_n\colon H^n(BS^1;\Z)\oplus H^n(S^2\times_WBT_2;\Z)\to H^n(\RP^2\times BS^1;\Z)$ is surjective for every $n\geq 0$. From the Mayer-Vietoris sequence associated to (\ref{dgr:pou2}) we see that the ring homomorphism
\[
\begin{split}
H^\ast(\Bcom U(2);\Z) \xrightarrow{(j^\ast, \phi^\ast)} & \,H^\ast(BS^1\vee (S^2\times_W BT_2);\Z) \\& \quad \cong \Z[t,\tilde{c}_1,\tilde{c}_2,c,U]/(2U,U^2,c^2,\tilde{c}_1U,cU,tz;\;z\in\{\tilde{c}_1,\tilde{c}_2,c,U\})
\end{split}
\]
is injective. Therefore, $H^*(\Bcom U(2);\Z)$ is isomorphic to the
subring $\im(j^*,\phi^*)$, where
$\im(j^*,\phi^*)_n=\ker(\pi_2^\ast-i^\ast)_n$ for each $n$.

The above reasoning implies that $\im(j^*,\phi^*)$ is generated by the kernel of $i^{*}$, given by $\ker i^*=(2c,\tilde{c}_1c,\tilde{c}_1^2-4\tilde{c}_2)$, together with the two classes $2t+\tilde{c}_1$ and $t^2+\tilde{c}_2$. Define a map $f\colon \Z[c_1,c_2,y_1,y_2]\to \im(j^*,\phi^*)$ by $f(c_1)=2t+\tilde{c}_1$, $f(c_2)=t^2+\tilde{c}_2$, $f(y_1)=2c$ and $f(y_2)=\tilde{c}_1c$. We claim that $f$ is surjective. Indeed, all there remains to check is that $f$ surjects onto $\ker i^\ast$. This follows from $\tilde{c}_1^2-4\tilde{c}_2=f(c_1^2-4c_2)$ and the relations $\tilde{c}_1(\tilde{c}_1^2-4\tilde{c}_2)=f(c_1^3-4c_1 c_2)$, $\tilde{c}_2(\tilde{c}_1^2-4\tilde{c}_2)=f(c_1^2c_2-4c_2^2)$ and $c(\tilde{c}_1^2-4\tilde{c}_2)=f(y_2c_1-2y_1c_2)$, together with the fact that $U$ annihilates $\ker i^\ast$. Using {\sc Singular} (see Section \ref{sec:singular1}) we can verify that the kernel of $f$ is $(2y_2-y_1c_1,y_1^2,y_1y_2,y_2^2)$. The fact that $c_1$ and $c_2$ (as generators of the polynomial ring) can indeed be identified with the pullback of the Chern classes follows from that fact that $(j^\ast,\phi^\ast)(\iota^\ast(c_1))=2t+\tilde{c}_1$ and $(j^\ast,\phi^\ast)(\iota^\ast(c_2))=t^2+\tilde{c}_2$.
\end{proof}

\begin{lemma} \label{lem:crsu2f2}
  The integral cohomology ring of $\Bcom SU(2)$ is as stated in Theorem \ref{thm:crsu2}.
\end{lemma}

\begin{proof}
  In \cite[Lemma 3.5]{Gr}, the author verifies the existence of a homotopy fiber sequence $\Bcom SU(2)\to\Bcom U(2)\to BS^1$, and hence a homotopy fiber sequence $S^1\to\Bcom SU(2)\to\Bcom U(2)$. The $E_2^{p,q}$-page of the associated Serre spectral sequence is concentrated in $q=0,1$. Let $e$ be the generator of $H^1(S^1;\Z)$, and consider $d_2\colon E_2^{0,1}\to E_{2}^{2,0}$. Then $d_2(e)=c_1\in H^2(\Bcom U(2);\Z)$, since
$\Bcom SU(2)$ is 3-connected (see \cite{Ad1}) and in particular its cohomology in degrees 1 and 2 is zero. Therefore $d_2\colon E_2^{4,1}\to E_2^{6,0}$ satisfies $d_2(ey_1)=2y_2$. Computing the remaining differentials $d_2$, we see that the $E_3^{p,q}$-page is concentrated in degree $q=0$, and is generated by $c_2,y_1$ and $x_2$, the generator of $y_2\Z/2y_2\Z$ in bi-degree $(6,0)$.
\end{proof}

\begin{remark}
  The integral cohomology groups of $\Bcom SU(2)$ were computed in \cite{Ad1}, and our computations agree with their results.
\end{remark}

\begin{corollary}\label{cor:cohEcomSU2}
$H^*(\Ecom SU(2);\Z)\cong \Z[y_1,x_2]/(2x_2,y_1^2,y_1x_2,x_2^2)$,
where $\deg y_1=4$ and $\deg x_2=6$.
\end{corollary}

\begin{proof}
We analyze the Serre spectral sequence $(E^{*,*},d)$ associated to the fiber sequence $SU(2)\to\Ecom SU(2)\to\Bcom SU(2)$ by comparing it to the one associated to the universal bundle $SU(2)\to ESU(2)\to BSU(2)$. Denote the latter by $({E^\prime} ^{*,*},d')$. By Theorem \ref{thm:crsu2}, $\Bcom SU(2)\to BSU(2)$ induces an injection $H^*(BSU(2);\Z)\to H^*(\Bcom SU(2);\Z)$ sending the second Chern class $c_2$ to the identically named class $c_2$. Let $e$ be the generator in $H^3(SU(2);\Z)$. Since $d_4^\prime (e)=c_2$, then also $d_4(e)=c_2$. Therefore, the only non-vanishing generators in the $E_5$-page are $y_1$ and $x_2$ in bi-degrees $(4,0)$ and $(6,0)$, respectively.
\end{proof}

\subsection{The $\F_2$-cohomology rings and Steenrod squares}

Now we move on to the $\F_2$-cohomology rings. For $U(2)$ there's not much left to do.

\begin{lemma}
  The $\F_2$-cohomology ring of $\Bcom U(2)$ is as claimed in Theorem \ref{thm:cru2}.
\end{lemma}

\begin{proof}
  From the presentation we've obtained for it, it is clear that the
  integral cohomology ring is torsion free. Then the presentation for
  the $\F_2$-cohomology ring follows from the universal coefficients
  theorem.
\end{proof}

For $SU(2)$, on the other hand, we still have some work to do. This includes figuring out a small part of the Steenrod algebra action.

\begin{lemma}\label{lem:f2crsu2}
  The $\F_2$-cohomology ring of $\Bcom SU(2)$ is as claimed in Theorem \ref{thm:crsu2}.
\end{lemma}

\begin{proof}
We use the mod 2 Serre spectral sequence of the homotopy fibration used in Lemma \ref{lem:crsu2f2}. Now let this time $e$ denote the generator of $H^1(S^1;\F_2)$. The same argument made above shows that $d_2(e)=\bar{c}_1$. In this case, the class $e\bar{y}_1$ in bi-degree $(4,1)$ maps to zero under $d_2$. The zeroth row of the $E_3$-page is generated by $\bar{c}_2$, $\bar{y}_1$ and $\bar{x}_2$. In the first row we have the class $e\bar{y}_1$, and the remaining classes are determined by the non-zero multiples of $\bar{y}_1$ in $H^*(\Bcom U(2);\F_2)$ (non-zero entries in bi-degree $(p,1)$ lie in honest subgroups of $H^{p+1}(\Bcom SU(2);\F_2)$, since the cohomology of the base is concentrated in even degrees). The spectral sequences collapses at the $E_3$-page and the only remaining product to check is $(e\bar{y}_1)^2$ (bi-degree $(8,2)$). This requires a different approach: because of possible multiplicative extensions we cannot decide if the product is zero or equal to $\bar{c}_2\bar{x}_2$ (bi-degree (10,0)).

Instead, we use the homotopy pushout for $G=SU(2)$ (Lemma \ref{lem:poBcomG1}) that reads
\begin{align}
\xymatrix{\RP^2\times \RP^\infty\ar[r]^-{i}\ar[d]_-{\pi_2}&S^2\times_{W}BT\ar[d]^-{\phi}\\ \RP^\infty\ar[r]^-{j}&\Bcom SU(2).}\label{dgr:posu2}
\end{align}
Let $x_1$ denote the (unique) lift of $e\bar{y}_1$ to $H^5(\Bcom SU(2);\F_2)$. We claim that $x_1$ is in the image of $\delta$, the connecting homomorphism in the Mayer-Vietoris sequence of (\ref{dgr:posu2}). To see this we use Lemma \ref{lem:f2crS2xBT}. Notice that $\Bcom SU(2)$ being 3-connected and considerations with the rank imply that $(\pi_2^*-i^*)\colon H^3(\RP^\infty;\F_2)\oplus H^3(S^2\times_WBT;\F_2)\to H^3(\RP^2\times\RP^\infty;\F_2)$ is injective and of rank 2. Therefore, the kernel of the connecting homomorphism $\delta \colon H^3(\RP^2\times\RP^\infty;\F_2)\to H^4(\Bcom SU(2);\F_2)$ has rank 2. Also, $H^4(S^2\times_{W}BT;\F_2)$ has rank 2, so that the relevant piece of the Mayer-Vietoris sequence in degree $4$ reads
\[0\to \F_2\to H^4(\Bcom SU(2);\F_2)\to \F_2^3\to \F_2^3\xrightarrow{\delta} H^5(\Bcom SU(2);\F_2)\to \cdots\] 
Since $H^4(\Bcom SU(2);\F_2)$ has rank 2 (by the first paragraph of this proof), exactness implies $\ker \delta=\F_2^2$ and therefore $\im \delta = \langle x_1\rangle$, which proves our claim. 

Now we can easily compute $x_1^2$ via $\Sq^5(x_1)$. Since $\Sq^n$ commutes with the connecting homomorphism, $\Sq^5(x_1)=\delta(\Sq^5(l))$ for some $l$, but $\Sq^5(l)=0$, because $l$ has degree 4. Therefore $x_1^2=\Sq^5(x_1)=0$. The presentation now follows.

Finally we must have $\beta(x_1)=x_2$, simply because $x_2$ is $2$-torsion meaning it must be in the image of $\beta$ and $x_1$ is the only non-zero class in degree $5$.
\end{proof}

\medskip

\paragraph*{The Steenrod algebra action.} We've already computed some of the Steenrod squares; let's take stock of what we have, what is obvious and what is left to do.

For $\Bcom U(2)$, the claimed total Steenrod squares for $\bar{c}_1$ and $\bar{c}_2$ are simply what they are well-known to be for the Chern classes in $H^*(BU(2);\F_2)$. Also, since all the cohomology of $\Bcom U(2)$ is in even degree, all Steenrod squares $\Sq^n$ with $n$ odd are automatically zero. The ring structure says that $\Sq^4(\bar{y}_1) = \bar{y}_1^2 = 0$, and the same for $\Sq^6(\bar{y}_2)$. That only leaves three squares to calculate: $\Sq^2(\bar{y}_1)$, $\Sq^2(\bar{y}_2)$ and $\Sq^4(\bar{y}_2)$.

For $\Bcom SU(2)$, the total Steenrod square of $\bar{c}_2$ again ultimately comes from $BU(2)$. Since $\bar{y}_1$ and $\bar{x}_2$ are the images of $\bar{y}_1$ and $\bar{y}_2$ in $H^*(\Bcom U(2); \F_2)$, we can already say that $\Sq^n(\bar{y}_1) = 0$ when $n>0$ except possibly for $n=2$, and that $\Sq^n(\bar{x}_2) = 0$ when $n>0$ except possibly for $n=2$ or $n=4$.

Of course, there is also $x_1$ which is not in the image of the map from $H^*(\Bcom U(2); \F_2)$. We already know that $\beta(x_1) = x_2$, which implies upon reducing mod 2, that $\Sq^1{x_1} = \bar{x}_2$. We now show that all higher Steenrod squares of $x_1$ vanish.

\begin{lemma}\label{lem:sqx1}
  The total Steenrod square of the class $x_1 \in H^5(\Bcom SU(2); \F_2)$ is given by $\Sq(x_1) = x_1 + \bar{x}_2$.
\end{lemma}

\begin{proof}
  We have $\Sq^2(x_1) \in H^7(\Bcom SU(2); \F_2) = 0$; and then, the Adem relation $\Sq^3 = \Sq^1 \Sq^2$ implies that $\Sq^3(x_1) = 0$ as well. From our calculation of the ring structure we know that $\Sq^5(x_1) = x_1^2 = 0$.

  Now we want to find $\Sq^4(x_1) \in H^9(\Bcom SU(2); \F_2) = \langle x_1 \bar{c}_2 \rangle$. If it happened that $\Sq^4 (x_1) = x_1 \bar{c}_2$, by the Adem relation $\Sq^5 = \Sq^1 \Sq^4$, we'd have $0 = \Sq^5(x_1) = \Sq^1(x_1 \bar{c}_2) = \Sq^1(x_1) \bar{c}_2 + x_1 \Sq^1(\bar{c}_2) = \bar{x}_2 \bar{c}_2$, contradicting the ring structure we obtained. Therefore, $\Sq^4(x_1) = 0$, as desired.
\end{proof}

Before we tackle the six remaining squares, we need to recall exactly how we defined the classes $c_1, c_2, y_1, y_2$ in $H^*(\Bcom U(2); \Z)$. This was done in the proof of Lemma \ref{lem:icru2}, which calculates the integral cohomology ring of $\Bcom U(2)$. Recall the inclusion $j \colon BS^1 \to \Bcom U(2)$ and conjugation map $\phi \colon S^2 \times_W BT_2 \to \Bcom U(2)$ appearing in the crucial homotopy pushout (\ref{dgr:pou2}). Together they specify a map $f := (j, \phi) \colon BS^1 \vee (S^2 \times_W BT_2) \to \Bcom U(2)$, which in the course of the proof of Lemma \ref{lem:icru2} was shown to induce an injective homomorphism on \emph{integral} cohomology. The integral cohomology classes in $\Bcom U(2)$ were then defined to satisfy $f^*(c_1)=2t+\tilde{c}_1$, $f^*(c_2)=t^2+\tilde{c}_2$, $f^*(y_1)=2c$ and $f^*(y_2)=\tilde{c}_1c$.

Notice that for $\F_2$-cohomology we have $f^*(\bar{y}_1) = 2\bar{c} = 0$, using $\bar{z}$ to denote the mod 2 reduction of an integral class $z$. So unlike what happens integrally, with $\F_2$-coefficients $f^*$ is not injective. The kernel of the $\F_2$-version of $f^*$ is the ideal generated by $\bar{y}_1$, but we'll only need a few low-degree instances of this fact.

\begin{lemma}
  The class $\bar{y}_1 \in H^4(\Bcom U(2); \F_2)$ satisfies $\Sq^2(\bar{y}_1)=0$, and thus the same is true for its image, the identically named $\bar{y}_1 \in H^4(\Bcom SU(2); \F_2)$.
\end{lemma}

\begin{proof}
  We have $f^*(\Sq^2(\bar{y}_1)) = \Sq^2(f^*(\bar{y}_1)) = 0$. While $f^*$ is not injective in general, on $H^6(\Bcom U(2); \F_2) = \langle\bar{c}_1^3,\bar{c}_1\bar{c}_2,\bar{y}_2\rangle$ it clearly is, by the above formulas, so we conclude that $\Sq^2(\bar{y}_1) = 0$.
\end{proof}

\begin{lemma} \label{lem:sq4y2}
  The class $\bar{y}_2 \in H^6(\Bcom U(2); \F_2)$ satisfies $\Sq^4(\bar{y}_2) = \bar{c}_1^2\bar{y}_2$.
\end{lemma}

\begin{proof}
  We have $f^*(\Sq^4(\bar{y}_2)) = \Sq^4(f^*(\bar{y}_2)) = \Sq^4(\bar{\tilde{c}}_1 \bar{c})$. The class $\tilde{c}_1$ was pulled back from the Chern class $c_1$, so $\bar{\tilde{c}}_1$ only has non-zero $\Sq^0$ and $\Sq^2$. The Cartan formula then says that $\Sq^4(\bar{\tilde{c}}_1 \bar{c}) = \bar{\tilde{c}}_1 \Sq^4(\bar{c}) + \bar{\tilde{c}}_1^2 \Sq^2(\bar{c})$. Since $c^2 = 0$, the mod 2 reduction satisfies $\Sq^4(\bar{c}) = \bar{c}^2 = 0$, too. Putting this together we have $f^*(\Sq^4(\bar{y}_2)) = \bar{\tilde{c}}_1^2 \Sq^2(\bar{c})$.
  
We next compute $\Sq^2(\bar{c})$. This square must be a linear combination of $\bar{\tilde{c}}_1 \bar{c}$ and $\bar{\tilde{c}}_2 \bar{U}$. This can be easily seen from the spectral sequence for the fibration $S^2\times_W BT_2\rightarrow \RP^2$ with $\F_2$--coefficients. The relevant part of the $E_2$-page (total degree $\leq 6$) is the reduction modulo $2$ of (\ref{eq:ssbt2}) with additional classes in bi-degrees $(1,0)$ and $(1,4)$. For degree reasons, however, the only classes contributing to $\Sq^2 (\bar{c})$ are those in bi-degree $(2,4)$, which are named $\bar{b}_1 \bar{c}$ and $\bar{b}_2 \bar{U}$ in the $E_2$--page, or $\bar{\tilde{c}}_1 \bar{c}$ and $\bar{\tilde{c}}_2 \bar{U}$ in the cohomology of $S^2\times_W BT_2$. We next observe that the contribution from $\bar{\tilde{c}}_2 \bar{U}$ is non-zero. Indeed, by Lemma \ref{lem:jbt2} and naturality of reduction modulo $2$, we have $i^\ast(\Sq^2(\bar{c}))=\Sq^2(i^\ast(\bar{c}))=\Sq^2(\bar{t}\bar{U})=\bar{t}^2\bar{U}$, as well as $i^\ast(\bar{\tilde{c}}_1\bar{c})=0$ and $i^\ast(\bar{\tilde{c}}_2 \bar{U})=\bar{t}^2\bar{U}$. To determine the contribution from $\bar{\tilde{c}}_1\bar{c}$ to $\Sq^2(\bar{c})$ we consider a map
  \[
  \rho: S^2\times_W (\CP^1\times\CP^1)\rightarrow S^2\times_{W} BT_2
  \]
induced by the obvious inclusion $\CP^1\times\CP^1\rightarrow\CP^\infty\times\CP^\infty\simeq BT_2$. We claim that $\rho^\ast(\bar{\tilde{c}}_1\bar{c})=\rho^\ast(\bar{\tilde{c}}_2\bar{U})$ and $\rho^\ast(\bar{\tilde{c}}_1\bar{c})\neq 0$. By naturality, we may alternatively check $\rho^\ast(\tilde{c}_1 c)=\rho^\ast(\tilde{c}_2U)$. To prove this we compare the spectral sequences for the fibrations over $\RP^2$. The relevant bi-degree is $(2,4)$, where $\rho^\ast$ corresponds to a map
  \[
  \sigma: H^2(\RP^2; \Z[\Z/2]\oplus \Z) \rightarrow   H^2(\RP^2; \Z)
  \]
  induced by the projection of coefficient modules $\Z[\Z/2]\oplus \Z\rightarrow \Z$ (see the proof of Lemma \ref{lem:crS2xBT2}). It can be checked that under the isomorphism $H^2(\RP^2; \Z[\Z/2]\oplus \Z)\cong H^2(\RP^2;\Z[\Z/2])\oplus H^2(\RP^2;\Z)\cong \Z\oplus \Z/2$ the class $b_1 c$ corresponds to the element $(1,1)\in \Z\oplus \Z/2$ and $b_2 U$ corresponds to $(0,1)$. Therefore, $\sigma(b_1 c)=\sigma(b_2 U)\neq 0$ in $\Z/2$ and, consequently, $\rho^\ast(\bar{\tilde{c}}_1 \bar{c})=\rho^\ast(\bar{\tilde{c}}_2 \bar{U})\neq 0$.
  
 We are now left with two possibilities for $\Sq^2(\bar{c})$ depending on whether $\rho^\ast(\Sq^2(\bar{c}))$ is zero or not. In the first case $\Sq^2(\bar{c})=\bar{\tilde{c}}_1\bar{c}+\bar{\tilde{c}}_2 \bar{U}$, while in the second case $\Sq^2(\bar{c})=\bar{\tilde{c}}_2 \bar{U}$. We claim that $\rho^\ast(\Sq^2(\bar{c}))=0$. Since $W$ acts freely on $S^2$, the space $M:= S^2\times_W(\CP^1\times\CP^1)$ is a $6$-manifold. The class $\rho^\ast(\bar{c})$ has degree four, so $\Sq^2(\rho^\ast(\bar{c}))$ lands in the top dimensional cohomology group of $M$. Thus, the Steenrod square can be computed by taking cup product with the second Wu class $\nu_2(M)$, which can be expressed in terms of Stiefel-Whitney classes as follows: $\nu_2(M)=w_2(M)+w_1(M)^2$. In the $E_2$--page of the spectral sequence for the fibration $M\rightarrow \RP^2$ there are no non-trivial cup products between the first and second columns, so $w_1(M)\rho^\ast(\bar{c})=0$. For the same reason, the product $w_2(M)\rho^\ast(\bar{c})$ only depends on the restriction of $w_2(M)$ to the fiber in the fibration sequence $\CP^1\times \CP^1\rightarrow M\rightarrow \RP^2$. Let $\iota$ denote the inclusion of the fiber. Since $M\rightarrow \RP^2$ is a locally trivial fiber bundle, the normal bundle of $\iota$ is trivial, hence $\iota^\ast(w_2(M))=w_2(\CP^1\times \CP^1)=0$. Therefore, $w_2(M)\rho^\ast(\bar{c})=0$ and consequently $\Sq^2(\rho^\ast(\bar{c}))=\nu_2(M)\rho^\ast(\bar{c})=0$, too. This proves $\Sq^2(\bar{c})=\bar{\tilde{c}}_1\bar{c}+\bar{\tilde{c}}_2\bar{U}$.
 
Finally, the result for $\Sq^2(\bar{c})$ implies that $f^\ast(\Sq^4(\bar{y}_2))=\bar{\tilde{c}}_1^3\bar{c}$. The kernel of $f^\ast$ is the ideal generated by $\bar{y}_1$, hence $f^\ast$ is injective in degree $10$. Since $f^\ast(\bar{c}_1^2 \bar{y}_2)=\bar{\tilde{c}}_1^3\bar{c}$, we conclude that $\Sq^4(\bar{y}_2)=\bar{c}_1^2\bar{y}_2$.
\end{proof}

\begin{lemma}
   We have $\Sq^4(\bar{x}_2) = 0$ for the class $\bar{x}_2 \in H^6(\Bcom SU(2); \F_2)$.
 \end{lemma}

\begin{proof} 
  The claim for $\bar{x}_2$ follows because $\bar{y}_2 \in H^6(\Bcom U(2); \F_2)$ maps to it and by the previous lemma $\Sq^4(\bar{y}_2)$ is a multiple of $\bar{c}_1$ which maps to $0$.

  Alternatively one can show $\Sq^4(\bar{x}_2) = 0$ directly by applying the Adem relation $\Sq^4 \Sq^1 = \Sq^5 + \Sq^2 \Sq^3$ to $x_1$. The left-hand side is $\Sq^4(\bar{x}_2)$ and the right-hand side is zero by Lemma \ref{lem:sqx1}.
\end{proof}

\begin{lemma} \label{lem:sq2x2}
  We have $\Sq^2(\bar{x}_2) = \bar{c}_2 \bar{y}_1$ for the class $\bar{x}_2 \in H^6(\Bcom SU(2); \F_2)$.
\end{lemma}

\begin{proof}
  This is similar to the argument for $x_1^2=0$ in the proof of Lemma \ref{lem:f2crsu2}. Consider again the Mayer-Vietoris sequence for the homotopy pushout square (\ref{dgr:posu2}).
  Let $H^*(\RP^2 \times \RP^\infty; \F_2) = \F_2[u,v]/(u^3)$, where $u$ and $v$ are the generators in degree $1$ for $H^*(\RP^2)$ and $H^*(\RP^\infty)$. We will show that $\bar{x}_2$ is in the image of connecting homomorphism, specifically that $\bar{x}_2 = \partial(u^2v^3)$. To do that we look at the portion of the Mayer-Vietoris around degrees 5 and 6 (let's omit the coefficients for brevity, they are $\F_2$ throughout this proof):
  \[ H^5(\RP^\infty) \oplus H^5(S^2 \times_W BT) \xrightarrow{\pi_2^* - i^*}
    H^5(\RP^2 \times \RP^\infty) \xrightarrow{\partial} H^6(\Bcom SU(2)).\]

  Now, $H^5(\RP^2 \times \RP^\infty) = \langle v^5, uv^4, u^2v^3 \rangle$. The image of $\pi_2^*$ is clearly $\langle v^5 \rangle$. We claim that the image of $i^*$ is rank one and is \emph{either} $\langle uv^4 \rangle$ \emph{or} $\langle uv^4 + u^2v^3\rangle$. Assuming that for the moment, we see that in either case $u^2 v^3$ is guaranteed to \emph{not} be in the image of $\pi_2^* - i^*$, and thus $\partial(u^2 v^3) \neq 0$, which establishes that $\partial(u^2 v^3) = \bar{x}_2$, as there is only one non-zero class in degree $6$ in $H^*(\Bcom SU(2))$.

  First we need to compute $H^*(S^2 \times_W BT; \F_2)$, but by a charming coincidence we already did that in the course of computing the cohomology of $\Bcom SO(3)_{\one}$. According to Lemma \ref{lem:f2crS2xBT}, the cohomology ring is $\F_2[u,z]/(u^3)$ (where we've used $z$ to name the generator of $H^{\ast}(BT; \F_2)$ instead of $\tilde{w}_2$ as in the lemma).  

  The map $i \colon \RP^2 \times \RP^\infty \to S^2 \times_W BT$ fits into a map of fibrations from the trivial fibration $\RP^\infty \to \RP^2 \times \RP^\infty \to \RP^2$ to the fibration considered above. This map of fibrations is the identity on the base and the inclusion $\RP^\infty = BZ \hookrightarrow BT$ on the fibre. So we can easily identify the map on $E_\infty$-pages: it is given by $u \mapsto u$ and $z \mapsto v^2$. In degree $5$ the only class in the $E_\infty$-page of the $S^2 \times_W BT$ fibration is $uz^2$ and this maps to $uv^4$ in the $E_\infty$-page of the trivial fibration. On cohomology, rather than $E_\infty$-pages, there can be correction terms of higher filtration, so we've only shown that $i^*(uz^2)$ is $uv^4$ or $uv^4+u^2v^3$, as mentioned earlier.

  All of this was to show that $\bar{x}_2 = \partial(u^2v^3)$. Since the connecting homomorphism commutes with Steenrod squares, we see that $\Sq^2(\bar{x}_2) = \partial(\Sq^2(u^2 v^3)) = \partial(u^2 v^5)$. The same argument as above, but two degrees higher, shows that $u^2 v^5$ is not in the image of $\pi_2^* - i^*$, and thus that $\partial(u^2v^5) \neq 0$. The only non-zero classes in degree $8$ in $H^*(\Bcom SU(2))$ are $\bar{c}_2 \bar{y}_1$ and $\bar{c}_2^2$. Since the pullback of $\bar{x}_2$ to $BT\subset \Bcom SU(2)$ is zero, $\bar{c}_2^2$ cannot contribute to the Steenrod square. Hence, $\Sq^2(\bar{x}_2)=\bar{c}_2 \bar{y}_1$ as claimed.
\end{proof}

\begin{lemma}
  We have $\Sq^2(\bar{y}_2) = \bar{c}_2 \bar{y}_1$ for the class $\bar{y}_2 \in H^6(\Bcom U(2); \F_2)$.
\end{lemma}

Note that the previous lemma, $\Sq^2(\bar{x}_2) = \bar{c}_2 \bar{y}_1$, is a direct consequence of this lemma, but that fact will be needed in the proof.

\begin{proof}
We first observe that $f^\ast(\Sq^2(\bar{y}_2))=\Sq^2(\bar{\tilde{c}}_1\bar{c})=0$, using knowledge of $\Sq^2(\bar{c})$ from the proof of Lemma \ref{lem:sq4y2}. Therefore, $\Sq^2(\bar{y}_2)\in \ker(f^\ast)=\langle \bar{c}_2 \bar{y}_1\rangle$. By Lemma \ref{lem:sq2x2} the image of $\Sq^2(\bar{y}_2)$ in $H^\ast(\Bcom SU(2);\F_2)$ is non-zero, hence $\Sq^2(\bar{y}_2)=\bar{c}_2\bar{y}_1$.
\end{proof}

That's the last of the Steenrod squares we needed and this concludes the proof of Theorems \ref{thm:cru2} and \ref{thm:crsu2}.

\begin{remark}
  Notice that it is not that easy to tell $\Bcom SU(2)$ and $\Ecom SU(2) \times BSU(2)$ apart: they become homotopy equivalent after looping once according to \cite[Theorem 6.3]{Ad5}, and they have isomorphic integral cohomology rings and $\F_2$-cohomology rings! But our calculation of the Steenrod algebra action on $H^*(\Bcom SU(2); \F_2)$ shows they are not equivalent, indeed $\Sq^{2}(\bar{x}_2) = \bar{c}_2 \bar{y}_1$, mixes the two factors.

  The non-splitting was already proven for an arbitrary compact and connected Lie group in \cite[Theorem 1.2.2]{Grthesis}, but this proof for $SU(2)$ via Steenrod squares is essentially different.
\end{remark}

\section{Cohomology of $\Bcom O(2)$} \label{sec:cro2}

In this section we compute the cohomology ring of $\Bcom O(2)$ with both integral and $\F_2$ coefficients, as well as the action of the Steenrod algebra on the $\F_2$-cohomology. These results were stated as Theorem \ref{thm:cro2}.

Let $\beta$ denote the Bockstein homomorphism for the coefficient sequence  $\Z\xrightarrow{2} \Z\rightarrow \Z/2$. Recall that
\[
H^\ast(BO(2);\F_2)\cong \F_2[w_1,w_2]\,,
\]
where $w_1$ and $w_2$ are the first two Stiefel-Whitney classes, and that
\[
H^\ast(BO(2);\Z)\cong \Z[W_1,W_2,p_1]/(2W_1,2W_2, W_2^2-p_1W_1)\,,
\]
where $W_i=\beta(w_i)\in H^{i+1}(BO(2);\Z)$ and where $p_1\in H^4(BO(2);\Z)$ is the first Pontryagin class of the universal $2$-plane bundle over $BO(2)$ (see \cite{BrownBOn}). In order to describe the corresponding cohomology rings of $\Bcom O(2)$ we will only need to define two additional classes. To define these we consider the following three maps:

\begin{itemize}
\item The map $j: BSO(2)\rightarrow \Bcom O(2)$ induced by the inclusion of the identity component $SO(2)\hookrightarrow O(2)$.
\item The map $k: BD_4\rightarrow \Bcom O(2)$ induced by the inclusion of the dihedral group $D_4\hookrightarrow O(2)$ as the subgroup generated by the reflections across the first and second axis.
\item A map $l: T_\tau \to \Bcom O(2)$, where $\tau$ is the self-map of $S^1\times S^1$ interchanging the two factors and $T_{\tau}$ is the corresponding mapping torus. This map is obtained as follows: Let the normalizer $N(D_4)$ act on $BD_4$ by conjugation and on $O(2)$ by translation from the right. One can show that $N(D_4)\cong D_8$ and that the action on $S^1\times S^1$ factors through $D_8\rightarrow D_8/D_4\cong \Z/2$ and swaps the two factors. The quotient space $O(2)\times_{N(D_4)}BD_4$ has the natural subcomplex $O(2)\times_{N(D_4)} (S^1\times S^1)\simeq T_{\tau}$ and $l$ is the composition
\[
O(2)\times_{N(D_4)}(S^1\times S^1)\xrightarrow{\textnormal{incl}.} O(2)\times_{N(D_4)} BD_4\xrightarrow{\phi} \Bcom O(2),
\]
where $\phi$ is the conjugation map. Note that $T_{\tau}$ is a (non-orientable) $3$-manifold, thus $H^3(T_{\tau};\F_2)\cong \F_2$.
\end{itemize}

Pick a generator $a\in H^2(BSO(2);\Z)$. Later we will see that there is a unique class $r\in H^2(\Bcom O(2);\Z)$ restricting along $j$ to $2a$ and restricting to zero along both $k$ and $l$. There is also a unique class $s\in H^3(\Bcom O(2);\F_2)$ which restricts to zero along both $j$ and $k$ and goes to the generator of $H^3(T_{\tau};\F_2)$ under $l$. In addition, we have the images of $w_1$, $w_2$, $W_1$, $W_2$ and $p_1$ in the cohomology of $\Bcom O(2)$, which we denote by the same letters.

Theorem \ref{thm:cro2} yields the rational cohomology ring of $\Bcom O(2)$ as a corollary. Note that this result cannot be directly deduced from the work of Adem, Cohen and Torres-Giese \cite[Theorem 6.1]{Ad5}, since their theorem only applies to the decorated version of the classifying space, i.e., to $\Bcom G_\one$.

\begin{corollary} \label{cor:qcohbcomo2}
We have an isomorphism of $\Q$--algebras $H^\ast(\Bcom O(2);\Q)\cong \Q[r]$ where $\deg (r)=2$.
\end{corollary}

We can now prove the application mentioned in the introduction:

\newcommand{\introthmref}{cor:tautological}
\begin{introcor*}
For any $n\geq 4$ the tautological $2$-plane bundle $\gamma_{2,n}\to \Gr_2(\R^n)$ does not admit a transitionally commutative structure. In particular, the universal bundle over $\Gr_{2}(\R^\infty)=BO(2)$ does not have a transitionally commutative structure.
\end{introcor*}

\begin{proof}
  It is enough to prove the case $n=4$ since given any transitionally commutative structure on $\gamma_{2,n}$ we could pull it back to one on $\gamma_{2,4}$ via the usual inclusion $\Gr_2(\R^4) \hookrightarrow \Gr_2(\R^n)$. Assume now that $\gamma_{2,4}$ had a transitionally commutative structure corresponding to some factorization $\Gr_2(\R^4) \xrightarrow{g} \Bcom O(2) \xrightarrow{\iota} BO(2)$ of the classifying map for $\gamma_{2,4}$. According to Borel's classic computation \cite{Borel}, $H^\ast(\Gr_2(\R^4); \Q) = \Q[p_1]/(p_1^2)$, where $p_1$ is indeed the pullback of $p_1 \in H^{\ast}(BO(2); \Q)$ under the classifying map of $\gamma_{2,4}$. Then the class $g^{\ast}(r) \in H^2(\Gr_2(\R^4); \Z)$ satisfies $\left( g^{\ast}(r) \right)^2 = 4p_1$, but this is impossible because, passing to rational cohomology, $H^2(\Gr_2(\R^4); \Q) = 0$ and thus contains no square root of $4p_1$.
\end{proof}

In Section \ref{sec:integralbcomo2} we begin the proof of Theorem \ref{thm:cro2} by calculating the integral cohomology ring, i.e. part \ref{thm:cro2-hz} of the theorem, but without proving the statement about the Bockstein homomorphism. In Section \ref{sec:f2bcomo2} we establish the result for $\F_2$-coefficients, i.e. part \ref{thm:cro2-f2} of the theorem, and compute the image of the Bockstein homomorphism. To determine the ring structure of $H^\ast(\Bcom O(2);\F_2)$ we use Steenrod squares, which we determine along the way. Putting everything together gives us the proof of Theorem \ref{thm:cro2}.

\subsection{The integral cohomology ring} \label{sec:integralbcomo2}

The basis for all our cohomology computations is the homotopy pushout
square given in Lemma \ref{lem:poBcomO2}. We start by describing the
cohomology rings of the spaces in the top left and right corner of
diagram (\ref{dgr:poBcomO2}) and the map between them. \bigskip

Let $\pi_1: O(2)\times_{N(D_4)}BD_4\rightarrow O(2)/N(D_4)\cong S^1$ be the projection and let $e\in H^1(S^1;\Z)$ be a generator. One can show that $N(D_4)=D_8$, which we will use from now on. Let $t=\pi_1^\ast(e)$ and let $\tilde{W}_1$, $\tilde{W}_2$ and $\tilde{p}_1$ denote the images of $W_1$, $W_2$ and $p_1$, respectively, under the map induced by $\iota\phi: O(2)\times_{D_8} BD_4\rightarrow BO(2)$.

\begin{lemma}\label{lem:cro2xbd54}
There is an isomorphism of graded rings
\[
H^*(O(2)\times_{D_8}BD_4;\Z)\cong \Z[t,\tilde{W}_1,\tilde{W}_2,\{x\},\tilde{p}_1]/I\,,
\]
where $\deg(\{x\})=3$, and $I$ is the ideal generated by $2\tilde{W}_1$, $2\tilde{W}_2$, $2\{x\}$, $2\tilde{p}_1$, $\tilde{W}_2^2+\tilde{p}_1\tilde{W}_1$, $t^2$, $t\tilde{W}_1$, $t\{x\}$ as well as $\{x\}^2$.
\end{lemma}

%The notation $\{x\}$ for a cohomology class may seem cumbersome at first, but it turns out to be very helpful and we will use it several times in the following two sections. We use curly brackets to denote the equivalence class of a cohomology class in a module of coinvariants.

\begin{proof}
We describe the $E_2$-page of the Serre spectral sequence associated to the homotopy fibration $BD_4\to O(2)\times_{D_8}BD_4\to O(2)/D_8$. Since $O(2)/D_8\cong S^1$ the $p$-th column of the $E_2$-page is zero for $p>1$, and thus $E_2^{*,*}\cong E_\infty^{*,*}$. Recall that
\[
H^*(BD_4;\Z)\cong \Z[x,y,w]/(2x,2y,2w, w^2+x^2y+xy^2)\,,
\]
where $x=\beta(u)$, $y=\beta(v)$, $w=\beta(uv)$, and $u$, $v$ are the basis elements for $H^1(BD_4;\Z/2)$ determined by reflection in the first and second axis. A generator of $\pi_1(O(2)/D_8)\cong \Z$ acts on the generators of $H^*(BD_4;\Z)$ by swapping $x$ and $y$, and fixing $w$.

The $0$-column and the $1$-column of the $E_2$-page are given by $H^0(S^1;\underline{H}^\ast(BD_4;\Z))$ and $H^1(S^1,\underline{H}^\ast(BD_4;\Z))$, respectively. The 0-column can be identified with the ring of invariants $H^*(BD_4;\Z)^{\Z}$, which is generated by $x+y,w$ and $xy$ with bi-degrees $(0,2)$, $(0,3)$ and $(0,4)$, respectively. By Poincar\'e Duality with twisted coefficients, the 1-column of the $E_2$-page is given by
\[
H^1(S^1;\underline{H}^\ast(BD_4;\Z))\cong H_0(S^1;\underline{H}^\ast(BD_4;\Z))\,,
\]
which can be identified with the module of coinvariants $H^\ast(BD_4;\Z)_{\Z}$. The structure of the $1$-column as a module over the $0$-column via cup product then corresponds to the usual structure of the coinvariants as a module over the ring of invariants.

For a class $z\in H^\ast(BD_4;\Z)$ let $\{z\}$ denote its image in $H^\ast(BD_4;\Z)_{\Z}$. The 1-column consists of linear combinations of $\{w^{\epsilon}x^ky^l\}$ for $\epsilon\in \{0,1\}$ and $0\leq l\leq k$. Since $\{w^{\epsilon}x^ky^l\}=w^\epsilon (xy)^l\{x^{k-l}\}$ and $\{x^{k+1}\}=(x+y)\{x^{k}\}+(xy)\{x^{k-1}\}$ the $1$-column is generated as a module over the $0$-column by the classes $\{1\}$ and $\{x\}$ of bi-degrees $(1,0)$ and $(1,2)$, respectively. This shows that $E_2^{*,*}$ is generated multiplicatively by $x+y$, $w$, $xy$, $\{1\}$ and $\{x\}$, hence so is $E_\infty^{*,*}$.
  
We need to describe the products of these generators. The only products missing are the products of two coinvariant classes, but they are all zero for dimension reasons. To finish the proof we must lift the multiplicative structure from the $E_\infty$-page to $H^\ast(O(2)\times_{D_8}BD_4;\Z)$. We first choose lifts of $x+y$, $w$ and $xy$. Recall that the map in cohomology induced by the inclusion $BD_4\hookrightarrow BO(2)$ is given on generators by $p_1\mapsto xy$, $W_1\mapsto x+y$ and $W_2\mapsto w$ (this is easily seen from the description via the Bockstein homomorphism given earlier, the fact that the Stiefel-Whitney classes $w_1$ and $w_2$ map to the symmetric polynomials $u+v$ and $uv$ in $H^\ast(BD_4;\Z/2)$ and $p_1$ maps to $w_2^2$ under the coefficient homomorphism $\Z\rightarrow \Z/2$). Since the inclusion factors as
\begin{align}
\xymatrix{BD_4\ar[r] & O(2)\times_{D_8}BD_4\ar[r]^{\iota\phi}&BO(2)}\label{dgr:fBD4}
\end{align}
the classes $\tilde{W}_1$, $\tilde{W}_2$ and $\tilde{p}_1$ can be taken as lifts of $x+y$, $w$ and $xy$, respectively. The products of these lifts are then determined by the products in $H^\ast(BO(2);\Z)$. Since $E^{p,q}=0$ for $p>1$, the classes $\{1\}$ and $\{x\}$ lift uniquely and there are no multiplicative extensions. Renaming $t:=\{1\}$ and extracting the relations from $H^*(BO(2);\Z)$ and $E_\infty^{*,*}$ gives the desired presentation. Also recall that $\{1\}$ is the Poincar\'e dual of a generator $e\in H^1(S^1;\Z)$, so $t$ can be identified with $\pi_1^\ast(e)$.
\end{proof}

Next consider the diagram
\begin{equation} \label{dgr:iBD4}
\xymatrix{
O(2)/D_8\times BZ\; \ar[r]^-{i} \ar[d]^-{}	& O(2)\times_{D_8} BD_4 \ar[d]^-{\pi_1} \\
O(2)/D_8 \ar@{=}[r]								& O(2)/D_8
}
\end{equation}
induced by the inclusion $Z\hookrightarrow D_4$, where both vertical maps are the projection onto the first factor. We want to describe the map $i$ in integral cohomology. By slight abuse of notation, we write
\[
H^*(O(2)/D_8\times BZ;\Z)\cong \Z[t,z]/(2z,t^2)\,,
\]
where $z\in H^2(BZ;\Z)\cong \Z/2$ is a generator and $t$ is the pullback of $e\in H^1(O(2)/D_8;\Z)$ under the projection, thus also the image of $t\in H^\ast(O(2)\times_{D_8}BD_4;\Z)$ under $i^\ast$.

\begin{lemma}\label{lem:istar}
In the presentation of Lemma \ref{lem:cro2xbd54} we have $i^*(\tilde{W}_1)=0$, $i^*(\{x\})=tz$, $i^*(\tilde{W}_2)=0$ and $i^*(\tilde{p}_1)=z^2$. 
\end{lemma}
\begin{proof}
Consider the inclusion $BZ\hookrightarrow BD_4$. On cohomology this map is determined by sending $x\mapsto z$ and $y\mapsto z$. Thus it sends $x+y\mapsto0$, $x^2y+xy^2\mapsto 0$ and $xy\mapsto z^2$. Under the inclusion $BZ\hookrightarrow BO(2)$ we therefore have $W_1\mapsto 0$, $W_2\mapsto 0$ and $p_1\mapsto z^2$. From the commutativity of
\[
\xymatrix{
O(2)/D_8\times BZ \ar[r]^-{i} \ar[d]^-{\pi_2}	& O(2)\times_{D_8} BD_4 \ar[d]^-{\iota\phi} \\
BZ\, \ar@{^(->}[r]						& BO(2)
}
\]
we now obtain $i^\ast(\tilde{W}_1)=0$, $i^\ast(\tilde{W}_2)=0$ and $i^\ast(\tilde{p}_1)=z^2$. To show that $i^\ast(\{x\})=tz$ we compare the Serre spectral sequences associated to the two fibrations in diagram (\ref{dgr:iBD4}). On the $E_2$-page the map of $1$-columns is the map $H^1(S^1;\underline{H}^\ast(BD_4;\Z))\rightarrow H^1(S^1;\underline{H}^\ast(BZ;\Z))$ induced by the morphism of local coefficient systems coming from $BZ\hookrightarrow BD_4$. Note that the local coefficient system $\underline{H}^\ast(BZ;\Z)$ is constant. By naturality of Poincar\'e duality, we can identify this map with the obvious map
\[
H^\ast(BD_4;\Z)_{\Z}\rightarrow H^\ast(BZ;\Z)_{\Z}\cong H^\ast(BZ;\Z)\, ,
\]
which is an isomorphism in degree $2$ sending $\{x\}\mapsto \{z\}$. The coinvariant class $\{z\}$ lifts to the class $tz$ and since both spectral sequences have $E_2^{p,q}=0$ for $p>1$ there are no ambiguities in choosing this lift.
\end{proof}

Now we return to the homotopy pushout (\ref{dgr:poBcomO2}) in Lemma \ref{lem:poBcomO2}. Fix a generator $a\in H^2(BSO(2);\Z)$. The map $q$ factors through the projection $O(2)/D_8\times BZ\to BZ$, thus $q^*(a)= z$. By Lemma \ref{lem:istar}, the homomorphism
\[H^n(BSO(2);\Z)\oplus H^n(O(2)\times_{D_8}BD_4;\Z)\xrightarrow{q^*-i^*} H^n(O(2)/D_8\times BZ;\Z)\]
is surjective for every $n\geq 0$. Thus, from the Mayer-Vietoris sequence for the homotopy pushout we get an injective ring map
\begin{equation} \label{eq:bvsod4}
H^*(\Bcom O(2);\Z)\to H^*(BSO(2)\vee O(2)\times_{D_8}BD_4;\Z)\,,
\end{equation}
whose image is isomorphic to $\ker (q^*-i^*)$. As a group we can describe the kernel by $\ker (i^*)\oplus J$, where $\ker (i^\ast)$ can be identified with the ideal $((0,\tilde{W}_1),(0,\tilde{W}_2))$ in the cohomology ring of the wedge sum, and $J$ is the $\Z$-linear span of $\{(a^{2l},\tilde{p}_1^l)\,|\, l\geq 0\}$. Then the ring structure on $\ker (q^\ast-i^\ast)$ (and thus on the image of (\ref{eq:bvsod4})) corresponds to componentwise products. A presentation for this ring will be a presentation for the cohomology ring $H^\ast(\Bcom O(2);\Z)$.\bigskip

The discussion in the preceding paragraph shows that there is a unique class $r\in H^2(\Bcom O(2);\Z)$ satisfying
\begin{align*}
j^\ast(r)&=2a;\\
\phi^\ast(r)&=0\, .
\end{align*} 
It is immediate that this class also satisfies $k^\ast(r)=0$ and $l^\ast(r)=0$ (where $k$ and $l$ were defined in the introduction to this chapter). Conversely, if $r$ is a class of degree $2$ in the kernel of $k^\ast$, then Lemma \ref{lem:cro2xbd54} implies that $\phi^\ast(r)=0$, since $\phi^\ast(r)$ is a multiple of $\tilde{W}_1$ and $BD_4\rightarrow O(2)\times_{D_8} BD_4$ detects this class. This shows that the conditions above specify the same class $r\in H^2(\Bcom O(2);\Z)$ as the one defined in the introduction to this chapter. Similarly, let $b_i\in H^{i+3}(\Bcom O(2);\Z)$ for $1\leq i\leq 3$ be the unique classes satisfying 
\begin{align*}
j^\ast(b_i)&=0\textnormal{ for all }i;\\
\phi^\ast(b_1)&= t\tilde{W}_2,\;  \phi^\ast(b_2)=\{x\} \tilde{W}_1\textnormal{ and }  \phi^\ast(b_3)= \{x\}\tilde{W}_2\,.
\end{align*}

The following intermediate result is part (2) of Theorem \ref{thm:cro2} modulo the statement about the Bockstein homomorphism.

\begin{proposition}\label{prop:cohbcomo2part1}
There is an isomorphism of graded rings
\[H^*(\Bcom O(2);\Z)\cong \Z[W_1,W_2,p_1,r,b_1,b_2,b_3]/I\,,\]
where $\deg(r)=2$, $\deg(b_i)=i+3$ for $1\leq i\leq 3$, and where $I$ is the ideal generated by $W_2^2-p_1W_1$, $r^2-4p_1$, $b_2 p_1-b_3 W_2$, $b_2 W_2-b_3 W_1$, $2W_i$, $rW_i$ and $b_1 W_i$ for $i=1,2$ as well as $2b_i$, $rb_i$ and $b_ib_j$ for $1\leq i,j\leq 3$.
\end{proposition}
\begin{proof}
Consider the ring homomorphism
\[
f: \Z[W_1,W_2,p_1,r,b_1,b_2,b_3]\rightarrow H^\ast(BSO(2)\vee O(2)\times_{N(D_4)}BD_4;\Z)
\]
sending $W_1\mapsto (0,\tilde{W}_1)$, $W_2\mapsto (0,\tilde{W}_2)$, $p_1\mapsto (a^2,\tilde{p}_1)$, $r\mapsto (2a,0)$, $b_1\mapsto (0,t\tilde{W}_2)$, $b_2\mapsto (0,\{x\}\tilde{W}_1)$ and $b_3\mapsto (0,\{x\}\tilde{W}_2)$. It is readily checked that the image of $f$ is isomorphic to $\ker(q^\ast-i^\ast)$, and we show in Section \ref{sec:singular2} that the kernel of $f$ is $I$ using {\sc Singular}.
\end{proof}

\subsection{The $\F_2$-cohomology ring and Steenrod squares} \label{sec:f2bcomo2}

The next lemma is the $\F_2$-analogue of Lemma \ref{lem:cro2xbd54}. If $e\in H^1(S^1;\F_2)$ is the generator and $\pi_1: O(2)\times_{D_8}BD_4\rightarrow O(2)/D_8\cong S^1$ is the projection, let $t=\pi_1^*(e)$. Thus, in this section $t$ is a mod $2$ cohomology class. Let $\tilde{w}_1$ and $\tilde{w}_2$ denote the images of $w_1$ respectively $w_2$ under the map induced by $\iota\phi: O(2)\times_{D_8}BD_4\rightarrow BO(2)$.

\begin{lemma}\label{lem:f2cro2xbd4}
There is an isomorphism of graded rings
\[H^*(O(2)\times_{D_8}BD_4;\F_2)\cong \F_2[t,\tilde{w}_1,\tilde{w}_2,\{u\}]/(t^2,t\tilde{w}_1,t\{u\},\{u\}^2),\] where $\deg(\{u\})=2$.
\end{lemma}
\begin{proof}
The proof is parallel to that of Lemma \ref{lem:cro2xbd54} and we leave the details to the reader. The new class $\{u\}$ is the (unique) lift of the coinvariant $\{u\}\in H^1(BD_4;\F_2)_{\Z}$ of bi-degree $(1,1)$, where the action of $\Z$ on $H^\ast(BD_4;\F_2)\cong \F_2[u,v]$ is by swapping $u$ and $v$.
\end{proof}

Write $H^*(O(2)/D_8\times BZ;\F_2)\cong \F_2[t,z]/(t^2)$, where $t\in H^1(O(2)/D_8;\F_2)$ and $z\in H^1(BZ;\F_2)$ are generators. Consider the inclusion $O(2)/D_8\times BZ\xhookrightarrow{i}O(2)\times_{D_8}BD_4$. 

\begin{lemma}\label{lem:f2istar}
In the presentation of Lemma \ref{lem:f2cro2xbd4} we have $i^*(\tilde{w}_1)=0$, $i^*(\{u\})=tz$ and $i^*(\tilde{w}_2)=z^2$. 
\end{lemma}
\begin{proof}
Again, the proof is parallel to that of Lemma \ref{lem:istar} and will be omitted.
\end{proof}

The Mayer-Vietoris sequence for cohomology with $\F_2$-coefficients applied to the homotopy pushout square (\ref{dgr:poBcomO2}) reads

\[\cdots \to H^{n-1}(O(2)/D_8\times BZ;\F_2)\xrightarrow{\delta} H^n(\Bcom O(2);\F_2)\to\]\[ H^n(BSO(2);\F_2)\oplus H^n(O(2)\times_{D_8}BD_4;\F_2)\xrightarrow{q^*-i^*}H^{n}(O(2)/D_8\times BZ;\F_2)\to\cdots \]
Let $H^*(BSO(2);\F_2)=\F_2[a]$ with $\deg (a)=2$. Then $q^*(a)=z^2$. Lemma \ref{lem:f2istar} shows that the image of $q^\ast-i^\ast$ is spanned by the elements $z^{2k}$ and $tz^{k}$ for $k\geq 0$. In particular, $q^*-i^*$ is not surjective, so we cannot regard $H^\ast(\Bcom O(2);\F_2)$ as a subring of $H^\ast(BSO(2)\vee O(2)\times_{D_8}BD_4;\F_2)$ as we have done before with integer coefficients. However, $q^*-i^*$ is surjective in even degrees and we find that
\begin{equation}
H^*(\Bcom O(2);\F_2)\cong \ker (q^*-i^*)\oplus \im(\delta)\label{eq:cgBcomO2}
\end{equation}
as a vector space, where the image of $\delta$ is the $\F_2$-linear span of all $\delta(z^{2i+1})$ for $i\geq 0$ and
\[
\ker (q^*-i^*)=((0,\tilde{w}_1))\oplus \F_2\langle (a^l,\tilde{w}_2^l)\;|\;l\geq 0\rangle\,,
\]
where $((0,\tilde{w}_1))\subset H^\ast(BSO(2)\vee O(2)\times_{D_8}BD_4;\F_2)$ is the ideal generated by $(0,\tilde{w}_1)$ and $\F_2\langle - \rangle$ means $\F_2$-linear span. This completely describes the cohomology groups $H^n(\Bcom O(2);\F_2)$. It remains to determine the products.\bigskip

It follows from the preceding paragraph that there is a unique $s\in H^3(\Bcom O(2);\F_2)$ satisfying
\begin{align*}
j^\ast(s)&=0;\\
\phi^\ast(s)&=\{u\}\tilde{w}_1.
\end{align*}
To see that this is the same cohomology class as the one specified at the beginning of this chapter, consider the map $\rho: T_\tau\simeq O(2)\times_{D_8}(S^1\times S^1)\rightarrow O(2)\times_{D_8}BD_4$ induced by the obvious inclusion of $S^1\times S^1\cong \RP^1\times \RP^1$ into $\RP^\infty\times \RP^\infty\simeq BD_4$. The Serre spectral sequence for the fibration of $T_{\tau}$ over the circle takes a similar form as the one for $O(2)\times_{D_8} BD_4$ (cf. the proof of Lemma \ref{lem:f2cro2xbd4}) in terms of invariants and coinvariants of $H^\ast(\RP^1\times \RP^1;\F_2)\cong \F_2[u,v]/(u^2,v^2)$. In particular, the class $\{u\}\in H^2(O(2)\times_{D_8} BD_4;\F_2)$ maps to the corresponding class $\{u\}\in H^2(T_{\tau};\F_2)$ (abusing notation), and similarly for $t$. In the spectral sequence the generator of $H^3(T_{\tau};\F_2)$ corresponds to the coinvariant class $\{uv\}$ of bi-degree $(1,2)$. Therefore,
\begin{equation} \label{eq:ttauuw1}
l^\ast(s)=\rho^\ast(\{u\}\tilde{w}_1)=\{u\}(u+v)=\{u^2\}+\{uv\}=\{uv\}
\end{equation}
is the generator for $H^3(T_{\tau};\F_2)$. We also have that $k^\ast(s)=0$, since the restriction along $BD_4\rightarrow O(2)\times_{D_8} BD_4$ sends $\{u\}$ to zero. Conversely, suppose that $s\in H^3(\Bcom O(2);\F_2)$ satisfies $j^\ast(s)=0$, $k^\ast(s)=0$ and $l^\ast(s)\neq 0$. By Lemma \ref{lem:f2cro2xbd4}, the class $\phi^\ast(s)$ is a linear combination of $\{u\}\tilde{w}_1$, $t\tilde{w}_2$ and terms involving only $\tilde{w}_1$ and $\tilde{w}_2$. By (\ref{eq:cgBcomO2}), $\phi^\ast(s)$ must not contain the term $t \tilde{w}_2$. Furthermore, since $k^\ast(s)=0$, it cannot contain any monomials in $\tilde{w}_1$, $\tilde{w}_2$. Finally, since $l^\ast(s)$ is non-zero, we must have $\phi^\ast(s)=\{u\}\tilde{w}_1$. Together this shows that the class $s$ specified via $j^\ast$ and $\phi^\ast$ is the same class as the one specified at the beginning of Section \ref{sec:cro2}.\bigskip

Next we compute the total Steenrod square $\Sq(s)$. To do this, we first consider $\{u\}\in H^2(O(2)\times_{D_8}BD_4;\F_2)$.

\begin{lemma} \label{lem:squ}
We have $\Sq^1 \{u\}=\{u\} \tilde{w}_1+t\tilde{w}_2$.
\end{lemma}
\begin{proof}
In general we have that $\Sq^1 \{u\}=\epsilon_1 \{u\} \tilde{w}_1+\epsilon_2 t\tilde{w}_2$ for some $\epsilon_1,\epsilon_2\in \F_2$. By Lemma \ref{lem:f2istar},
\[
\epsilon_2 tz^2=i^\ast(\Sq^1\{u\})=\Sq^1(tz)=tz^2\,,
\]
hence $\epsilon_2=1$. We claim that also $\epsilon_1=1$. Since $\rho^\ast(\{u\}\tilde{w}_1)=\{uv\}$ (see (\ref{eq:ttauuw1})) and similarly $\rho^\ast(t\tilde{w}_2)=\{uv\}$, it is enough to check that $\Sq^1\{u\}$ computed in $H^\ast(T_{\tau};\F_2)$ is zero. Since $T_{\tau}$ is a $3$-manifold, $\Sq^1\{u\}$ is given by cup product with the first Wu class $\nu_1(T_\tau)$. The latter is the same as the first Stiefel-Whitney class $w_1(T_{\tau})$ of the manifold $T_{\tau}$ (by the formula $\Sq(\nu)=w$). Using the fact that the fiber of $T_{\tau}\rightarrow S^1$ is orientable, one can see that $w_1(T_{\tau})$ restricts to zero on the fibers. From the spectral sequence for the fibration $T_{\tau}\rightarrow S^1$ we then see that $\Sq^1\{u\}=w_1(T_{\tau})\{u\}$ must be zero in $H^\ast(T_{\tau};\F_2)$.
\end{proof}

Let $\bar{r}\in H^2(\Bcom O(2);\F_2)$ be the reduction modulo $2$ of the class $r\in H^2(\Bcom O(2);\Z)$. Then clearly $\Sq^1(\bar{r})=0$, and $\Sq^2(\bar{r})=\bar{r}^2=0$, since $r^2=4p_1$ by Proposition \ref{prop:cohbcomo2part1}. Thus, the following lemma completes the proof of part (3) in Theorem \ref{thm:cro2}.

\begin{lemma} \label{lem:sqd2}
We have $\Sq^1s=w_2 \bar{r}$ and $\Sq^2s=w_1^2s$.
\end{lemma}
\begin{proof}
We begin with $\Sq^2 s$. Since the degree of $\Sq^2 s$ is odd, it is uniquely characterized by its image under $j^\ast$ and $\phi^\ast$. We have $j^\ast(\Sq^2 s)=0$ (by definition of $s$) and
\begin{equation*}
\begin{split}
  \phi^\ast(\Sq^2 s)=\Sq^2(\{u\}\tilde{w}_1)&= \Sq^2(\{u\})\tilde{w}_1+ \Sq^1(\{u\}) \tilde{w}_1^2 \\& =(\{u\}\tilde{w}_1 +t\tilde{w}_2)\tilde{w}_1^2 \quad [\Sq^2\{u\}=\{u\}^2=0\textnormal{ and Lemma \ref{lem:squ}}]\\&=\{u\}\tilde{w}_1^3\quad [\textnormal{since } t\tilde{w}_1=0]\\&= \phi^\ast(w_1^2s)\,,
\end{split}
\end{equation*}
and therefore $\Sq^2s=w_1^2s$.\bigskip

Now $j^\ast(\Sq^1 s)=\Sq^1(j^\ast s)=0$ and
\begin{equation*}
\begin{split}
\phi^\ast(\Sq^1s)=\Sq^1(\{u\} \tilde{w}_1)&=t\tilde{w}_1 \tilde{w}_2  \quad [\textnormal{by Lemma \ref{lem:squ}}] \\&= 0 \quad [\textnormal{since }t\tilde{w}_1=0]\, .
\end{split}
\end{equation*}
Therefore $\Sq^1s \in \im(\delta)$, so either $\Sq^1s=0$ or $\Sq^1s=\delta(z^3)$. We claim that $\Sq^1s\neq 0$.

Assume, for contradiction, $\Sq^1s=0$. Recall that $\Sq^1$ can be expressed as the composite of the Bockstein homomorphism $\beta$ and the natural transformation $R: H^\ast(-;\Z)\rightarrow H^\ast(-;\Z/2)$ induced by reduction modulo $2$. By assumption, $R(\beta(s))=0$, hence
\[
\beta(s)\in \ker(R)=\im(H^4(\Bcom O(2);\Z)\stackrel{2}{\longrightarrow}H^4(\Bcom O(2);\Z))\cong \Z\langle 2p_1\rangle\,,
\]
where the isomorphism follows from Proposition \ref{prop:cohbcomo2part1}. At the same time $\beta(s)$ is torsion, hence $\beta(s)=0$. Therefore, $s\in \ker(\beta)=\im(R)$. However, there is no element which reduces modulo $2$ to $s$. The only candidate is $b_1$, but
\begin{equation} \label{eq:rb1}
\phi^\ast(R(b_1))=R(\phi^\ast b_1)=R(t\tilde{W}_2)=t \Sq^1\tilde{w}_2=t\tilde{w}_1\tilde{w}_2=0\,,
\end{equation}
while $\phi^\ast(s)\neq 0$. We arrive at a contradiction. Hence $\Sq^1s=\delta(z^3)$. The assertion of the lemma follows from the identity $\delta(z^3)=w_2\bar{r}$ which will be established in the proof of Proposition \ref{prop:cohbcomo2part2} (without making our argument circular).
\end{proof}

The following proposition corresponds to part \ref{thm:cro2-f2} of Theorem \ref{thm:cro2} as well as the statement about the Bockstein homomorphism in part \ref{thm:cro2-hz}. It thus concludes the proof of Theorem \ref{thm:cro2}.

\begin{proposition} \label{prop:cohbcomo2part2}
There is an isomorphism of graded rings
\[H^*(\Bcom O(2);\F_2)\cong \F_2[w_1,w_2,\bar{r},s]/(\bar{r} w_1,\bar{r}^2, \bar{r}s, s^2)\, ,
\]
where $\deg(\bar{r})=2$ and $\deg(s)=3$. Moreover, $\beta(s)=b_1$, $\beta(w_1 s)=b_2$ and $\beta(w_2 s)=b_3$ in $H^\ast(\Bcom O(2);\Z)$.
\end{proposition}
\begin{proof}
The isomorphism (\ref{eq:cgBcomO2}) and Lemma \ref{lem:f2cro2xbd4} show that $H^\ast(\Bcom O(2);\F_2)$ is generated as a ring by $w_1$, $w_2$, $s$ and $\delta(z^{2i+1})$ for $i\geq 0$. More precisely, the classes $\delta(z^{2i+1})$ account for $\im(\delta)$, $w_2$ accounts for $\F_2\langle (a^l,\tilde{w}_2^l)\,|\, l\geq 0\rangle$, since $j^\ast(w_2^l)=a^l$ and $\phi^\ast(w_2^l)=\tilde{w}_2^l$, and finally $\tilde{w}_1$, $\tilde{w}_2$ and $\phi^\ast(s)=\{u\} \tilde{w}_1$ are enough to generate the ideal $((0,\tilde{w}_1))$.\bigskip

\paragraph*{Multiplicative relations:}
\begin{itemize}
\item $w_1 \delta(z)=0$ and $s \delta(z)=0$. These products have odd degree, and since $\delta(z)$ is in $\ker (\phi^\ast)$ and $\phi^\ast$ is injective in these degrees, the products must be zero.

\item $\delta(z)^2=0$. This follows from $\delta(z)^2=\Sq^2\delta(z)=\delta(\Sq^2 z)=0$, since $z$ has degree one.

\item $w_2 \delta(z^{2i-1})=\delta(z^{2i+1})$ for all $i>0$. First, we establish the relation $R(r)=\delta(z)$. We have that $j^\ast(R(r))=R(2a)=0$ and $\phi^\ast(R(r))=0$, so $R(r)$ is either $0$ or $\delta(z)$. Since $r$ is not divisible by $2$, $R(r)$ is non-zero, hence $R(r)=\delta(z)$. For every $i\geq 0$ the monomial $rp_1^i\in H^{4i+2}(\Bcom O(2);\Z)$ is not divisible by $2$, so $R(rp_1^i)=w_2^{2i}\delta(z)$ is non-zero, hence also $w_2^i \delta(z)\neq 0$. Since $w_2^i \delta(z)$ lies in the image of $\delta$, we must have
\[
w_2^i \delta(z)=\delta(z^{2i+1})\,.
\]
In particular, this shows $w_2\bar{r}=\delta(z^3)$ and finishes the proof of Lemma \ref{lem:sqd2}.

\item $s^2=0$. Using the Adem relation $\Sq^3=\Sq^1\Sq^2$, Lemma \ref{lem:sqd2} and $w_1\bar{r}=0$ we get
\[
s^2=\Sq^3s=\Sq^1\Sq^2 s=\Sq^1(w_1^2s)=w_1^2 \Sq^1 s=w_1^2w_2\bar{r}=0\,.
\]
\end{itemize}

Together this shows that we have a surjective ring homomorphism
\[
\F_2[w_1,w_2,\bar{r},s]/(\bar{r}w_1,\bar{r}^2,\bar{r}s,s^2)\rightarrow H^\ast(\Bcom O(2);\F_2)\,.
\]
This is easily seen to be an isomorphism, by comparing the dimension of the domain and target in every degree, using (\ref{eq:cgBcomO2}) and Lemma \ref{lem:f2cro2xbd4}.\bigskip

\paragraph*{The image of the Bockstein homomorphism:}\mbox{}\medskip

Note that the image of multiplication by $2$ on $H^\ast(\Bcom O(2);\Z)$ is spanned by $2r^ip_1^j$ for $i=0,1$ and $j\geq 0$. In particular, $R(b_1)$, $R(b_2)$ and $R(b_3)$ are all non-zero.

\begin{itemize}
\item $\beta(s)=b_1$. We have that $j^\ast(R(b_1))=0$ and $\phi^\ast(R(b_1))=0$ (see (\ref{eq:rb1})), hence $R(b_1)$ is in the image of the connecting homomorphism, i.e. $R(b_1)=w_2 \bar{r}$. On the other hand, we also have that $R(\beta(s))=\Sq^1(s)=w_2\bar{r}$, hence $\beta(s)-b_1\in \ker (R)\cong  \Z\langle 2p_1\rangle$. Since both $\beta(s)$ and $b_1$ are torsion, we must have $\beta(s)-b_1=0$.

\item $\beta(w_1 s)=b_2$. Again, we compute $j^\ast(R(b_2))=0$ and
\[
  \begin{aligned}
    \phi^\ast(R(b_2))&=R(\phi^\ast(b_2))=R(\{x\}\tilde{W}_1)\\
    &=\{u^2\} \tilde{w}_1^2 =(\{u\}\tilde{w}_1+t\tilde{w}_2)\tilde{w}_1^2=\{u\}\tilde{w}_1^3=\phi^\ast(w_1^2s)\, .
\end{aligned}
\]
This implies that $R(b_2)=w_1^2 s=\Sq^1(w_1 s)=R(\beta(w_1 s))$, hence $\beta(w_1 s)-b_2=0$.

\item $\beta(w_2 s)=b_3$. We find $j^\ast(R(b_3))=0$ and $\phi^\ast(R(b_3))=\phi^\ast(w_1w_2 s)$. However, in this degree $j^\ast$ and $\phi^\ast$ only determine $R(b_3)$ modulo the image of the connecting homomorphism, so that $R(b_3)=w_1w_2 s+\epsilon w_2^2 \bar{r}$ where $\epsilon$ may be $0$ or $1$. We also have that $R(\beta(w_2 s))=w_1w_2 s+w_2^2 \bar{r}$, hence $\beta(w_2 s)=b_3+ (1+\epsilon) r p_1$ modulo $\Z\langle 2rp_1\rangle$. Again, since both $\beta(w_2 s)$ and $b_3$ are torsion and $rp_1$ is of infinite order, this implies that $\beta(w_2s)=b_3$.
\end{itemize}
This finishes the proof of the proposition.
\end{proof}

\section{A homotopy pushout square for $\Ecom G$}

Just as all our cohomology calculations for $\Bcom G$ were based on the homotopy pushout squares in Lemmas \ref{lem:poBcomG1} and \ref{lem:poBcomO2}, our determination of the homotopy type of $\Ecom G$ will also be based on describing those spaces as homotopy pushouts. The goal of this section is to obtain such a description.

We start with finding the homotopy fiber of the variant of the conjugation map that lands in $BG$. Let $G$ be a compact Lie group, $H\leqslant G$ a closed subgroup and $N(H)$ the normalizer of $H$ in $G$. For any $N'\leqslant N(H)$ there is a conjugation map $\psi: G\times_{N'} BH\rightarrow BG$. Consider the action of $N(H)$ on $G/H$ given by
\begin{equation} \label{eq:act1}
n\cdot gH=ngn^{-1}H\quad (n\in N(H),\, gH\in G/H)\,,
\end{equation}
and write $G\times_{N'} G/H$ for the quotient of $G\times G/H$ by the induced diagonal $N'$-action.

\begin{lemma} \label{lem:hofibconj}
For any closed subgroup $N'\leqslant N(H)$ there is a homotopy fiber sequence
\[
G\times_{N'} G/H\longrightarrow G\times_{N'} BH\stackrel{\psi}{\longrightarrow} BG\, .
\]
\end{lemma}

\begin{proof}
We replace the conjugation map $\psi$ by an equivalent fibration. Let $B(G,G/H)$ be the bar construction for the action of $G$ on $G/H$ by left-translation. Let us write $\mu: G\times BG\rightarrow BG$ for the $G$-action on $BG$ induced by conjugation in $G$. Similarly, we have an action of $N(H)$ on $B(G,G/H)$ induced by the conjugation action of $N(H)$ on $G$ and the action of $N(H)$ on $G/H$ specified in (\ref{eq:act1}). Let us write $\pi: B(G,G/H)\rightarrow BG$ for the projection map. Then the composite map $\mu\circ (\id\times \pi): G\times B(G,G/H)\rightarrow BG$ factors through a map
\begin{equation*}
\tilde{\psi}: G\times_{N'} B(G,G/H)\longrightarrow BG\, .
\end{equation*}
Now consider the homotopy equivalence $BH\stackrel{\sim}{\longrightarrow} B(G,G/H)$ induced by $H\hookrightarrow G$. This is an $N(H)$-equivariant map. After taking the product with $G$ we obtain a homotopy equivalence of free $N(H)$-spaces $G\times BH\simeq  G\times B(G,G/H)$, thus an $N(H)$-equivalence. In particular, there is an induced homotopy equivalence of orbit spaces
\[
f: G\times_{N'} BH\stackrel{\sim}{\longrightarrow} G\times_{N'} B(G,G/H)\, .
\]
This fits into a commutative diagram
\[
\xymatrix{
G\times_{N'} BH \ar[r]^-{\psi} \ar[d]^-{f}_-{\simeq}	& BG \ar@{=}[d]	\\
G\times_{N'} B(G,G/H) \ar[r]^-{\tilde{\psi}}			& BG\,,
}
\]
thus $\hofib(\psi)\simeq \hofib(\tilde{\psi})$. To compute the homotopy fiber of $\tilde{\psi}$ we check that $\tilde{\psi}$ is a fibration. To see this, it is enough to show that $\mu\circ (\id\times \pi)$ is a fibration and we can check this separately for $\mu$ and for $\pi$. Clearly, $\pi$ is a fibration. Moreover, $\mu$ is isomorphic to the projection onto $BG$ via the shear map on $G\times BG$ and thus is a fibration too. The homotopy fiber of $\tilde{\psi}$ is then simply the fiber over the basepoint, which is $G\times_{N'}G/H$.
\end{proof}

\begin{remark} \label{rem:hofibconj}
If $G$ is connected, $H=T$ is a maximal torus for $G$, $N'=N(T)$ and $W=N(T)/T$ is the Weyl group, the homotopy fiber sequence of the lemma takes the form
\begin{equation} \label{eq:hfib}
G/T\times_W G/T\longrightarrow G/T\times_W BT\stackrel{\iota\phi}{\longrightarrow} BG\, ,
\end{equation}
where an element $w=nT\in W$ acts on $G/T$ by $w\cdot gT=gn^{-1}T$. To see that (\ref{eq:hfib}) is a fiber sequence note that $G\times_{N(T)}BT=G/T\times_W BT$ and that there is a homeomorphism
\[
G\times_{N(T)} G/T \stackrel{\cong}{\longrightarrow} G/T\times_W G/T 
\]
sending an equivalence class $[g,g'T]\mapsto [gT,gg' T]$, whose inverse given by $[gT,g'T]\mapsto [g, g^{-1}g' T]$.
\end{remark}

\begin{remark}
In view of \cite[Theorem 6.1]{Ad5} the homotopy fiber sequence (\ref{eq:hfib}) gives an alternative and short proof of the first statement in \cite[Corollary 7.4]{Ad1} which says that for a compact connected Lie group $G$ there is an isomorphism of algebras
\[
H^\ast(\Ecom G_\one;\Q)\cong (H^\ast(G/T;\Q)\otimes H^\ast(G/T;\Q))^W\, .
\]
This follows from the comparison theorem for the Serre spectral sequence.
\end{remark}

In the following we only compute the homotopy types of $\Ecom O(2)$ and $\Ecom SU(2)$, but according to the following proposition this covers all of our targets.

\begin{proposition}\label{prop:equivEcom}
The are homotopy equivalences \[\Ecom SU(2) \simeq \Ecom U(2) \simeq \Ecom SO(3)_\one.\]
\end{proposition}

This is a direct corollary of Lemma 1.2.8 from \cite{Grthesis}, which we now state:

\begin{lemma}[{\cite[Lemma 1.2.8]{Grthesis}}]\label{lem:covhom}
  If $\tilde{G} \to G$ is a homomorphism of compact connected Lie groups which is also a covering map, then the following induced square is a homotopy pullback square:
\[\xymatrix{\Bcom \tilde{G}_\one\ar[r]\ar[d]&B\tilde{G}\ar[d]\\
  \Bcom G_\one\ar[r]&BG}\]  
\end{lemma}

\begin{proof}[Proof of Proposition \ref{prop:equivEcom}]
Applying the previous lemma to the double cover homomorphisms $SU(2) \to SO(3)$ and $S^1 \times SU(2) \to U(2)$ (given by $(\lambda, M) \mapsto \lambda M$), we get the following homotopy pullbacks squares:
\[\xymatrix{\Bcom SU(2)\ar[r]\ar[d]&BSU(2)\ar[d]\\
  \Bcom SO(3)_\one\ar[r]&BSO(3)}
   \;\text{ and }\;
   \xymatrix{BS^1\times \Bcom SU(2)\ar[r]\ar[d]&BS^1\times BSU(2)\ar[d]\\
\Bcom U(2)\ar[r]&BU(2).}\]
It follows the homotopy fibers of the horizontal arrows are weakly homotopy equivalent, so that
\[\Ecom SU(2)\simeq \Ecom SO(3)_\one\;\text{ and }\Ecom U(2)\simeq ES^1\times\Ecom SU(2)\simeq\Ecom SU(2).\qedhere\]
\end{proof}

Now we will use Lemma \ref{lem:hofibconj} to obtain descriptions for $\Ecom SU(2)$ and $\Ecom O(2)$ as homotopy pushouts. A common language for both was setup in Section \ref{sec:poBcomG1}, which we now recall. Let $G$ be either $SU(2)$ or $O(2)$ and let $Z$ be its center. 
\begin{itemize}
\item For $G = SU(2)$, let $A=T=S^1$ the maximal torus, and $H = Z$.
\item For $G = O(2)$, let $A=D_4$, and $H = SO(2)$.
\end{itemize}
In both cases, let $N=N_G(A)$ and $W = N/A$. With that notation, for both groups we had the following homotopy pushout:
\[\xymatrix{
    G/N \times BZ \ar[d] \ar[r] & G/A \times_W BA \ar[d] \\
    BH \ar[r] & \Bcom G \\
}\]

\begin{lemma}\label{lem:poEcomG}
  With the above notation for $G=SU(2)$ or $G=O(2)$, we have the following homotopy pushout square for $\Ecom G$:
  \[\xymatrix{
      G/N \times G/Z \ar[d] \ar[r] & G/A \times_W G/A \ar[d] \\
      G/H \ar[r] & \Ecom G\,, \\
    }\]
  where the left vertical arrow is the composite of the projection onto $G/Z$ and the canonical projection $G/Z \to G/H$, and the top horizontal map is given by $(gN, xZ) \mapsto [gA, xgA]$.
\end{lemma}

\begin{proof}
  By composing with the inclusion $\iota \colon \Bcom G \to BG$, we can
regard all the spaces as objects in the category of spaces over $BG$
and by the second Mather cube theorem, the homotopy fibers over the
basepoint of $BG$ also form a homotopy pushout square. Let's look at
each corner.

\begin{itemize}
\item The bottom right corner is
  $\hofib(\Bcom G \to BG) \simeq \Ecom G$, which is the reason we are
  doing this.
\item The bottom left corner is $\hofib(BH \to BG)$ where the map is
  the one induced by the inclusion $H \hookrightarrow G$. Thus, this
  corner is $G/H$.
\item Similarly, the top right corner is $G/N \times G/Z$, because the
  map $G/N \times BZ \to BG$ factors through $BZ \to BG$ (whose
  homotopy fiber is $G/Z$).
\item Finally, the top right corner is the most interesting one. By
  Lemma \ref{lem:hofibconj}, it is $G \times_N G/A$ and by Remark
  \ref{rem:hofibconj} this is homeomorphic to $G/A \times_W G/A$.
\end{itemize}

The left vertical map is simply the projection. To understand the top
horizontal map, say $\sigma$, notice that to figure out what the top left
corner of the square is we could have also applied lemma
\ref{lem:hofibconj} to the group $H = Z$ and the subgroup $N$ of its
normalizer (here, of course, $N_G(Z) = G$). We would have obtained the
description $G \times_N G/Z$ for the space in the top left corner. The
$G \times_N (-)$ descriptions have the advantage that the top
horizontal map $G \times_N G/Z \to G \times_N G/A$ is simply the map
induced by the canonical projection $G/Z \to G/A$, namely
$[g, xZ] \mapsto [g, xA]$. Let's see what this map corresponds to for
the homeomorphic descriptions of the spaces we used in the square
above. We have the following commuting square where the vertical maps
are homeomorphisms (the one on the right is the one from Remark
\ref{rem:hofibconj}):

\[\xymatrix{
    G \times_N G/Z \ar[r] \ar[d]_{\cong} & G \times_N G/A \ar[d]^{\cong} \\
    G/N \times G/Z \ar[r]^-{\sigma} & G/A \times_W G/A. \\
  }\]

The left vertical homeomorphism is given by
$[g,xZ] \mapsto (gN, gxg^{-1}Z)$ (to check this is well-defined,
recall the diagonal $N$ action on $G \times G/Z$ used in the top left
corner is given by $n \cdot (g, xZ) = (gn^{-1}, nxn^{-1}Z)$), and its
inverse is then $(gN, xZ) \mapsto [g, g^{-1}xgZ]$. Following the upper
path through the square, we see the bottom horizontal map is given by
$\sigma(gN, xZ) = [gA, xgA]$.
\end{proof}

\section{The homotopy type of $\Ecom O(2)$}

Throughout this section, we write $R_{\alpha}\in O(2)$ for a rotation in the plane by an angle $\alpha\in [0,2\pi)$ and $r_{\beta}\in O(2)$ for a reflection in a line which makes an angle $\beta\in [0,\pi)$ with the first axis. The main goal of this section is to prove the following result stated in the introduction.

\renewcommand{\introthmref}{thm:ecomo2}
\begin{introthm*}
There is a homotopy equivalence $\Ecom O(2) \simeq S^2 \vee S^2 \vee S^3$.
\end{introthm*}

\begin{proof}
Recall that $N(D_4)=D_8$, so that the homotopy pushout diagram for $\Ecom O(2)$ given in Lemma \ref{lem:poEcomG} reads
\begin{equation} \label{dgr:ecomo2}
\xymatrix{
  O(2)/D_8 \times O(2)/\{\pm I\} \ar[r]^-{s} \ar[d]^-{r} & O(2)/D_4 \times_{D_8/D_4} O(2)/D_4 \ar[d]	\\
  O(2)/SO(2) \ar[r]                                      & \Ecom O(2)\,,
}
\end{equation} 
where $s(gD_8, xZ) = [gD_4, xgD_4]$ and $r(gD_8, xZ) = xSO(2)$.

The next step is to identify the spaces in the top half of the diagram with more familiar spaces. There are homeomorphisms
\begin{equation} \label{eq:homeo1}
  O(2)/D_8 \times O(2)/\{\pm I\}\cong S^1\times (S^1 \sqcup S^1)
\end{equation}
and
\begin{equation} \label{eq:homeo2}
O(2)/D_4 \times_{D_8/D_4} O(2)/D_4 \cong S^1\times S^1\,.
\end{equation}

First we explain (\ref{eq:homeo1}). We have $O(2)/D_8= SO(2)/\langle R_{\pi/2}\rangle$, where $\langle R_{\pi/2}\rangle$ denotes the cyclic group generated by the rotation $R_{\pi/2}$. Thus $O(2)/D_8$ is an $S^1$ parametrizing rotations modulo $\pi/2$. The space $O(2)/\{\pm I\}$ is the disjoint union of two circles. The first one is $SO(2)/\langle R_{\pi}\rangle$ and parametrizes rotations modulo $\pi$. The second one parametrizes reflections $r_{\beta}\in O(2)$, where $\beta$ is to be taken modulo $\pi/2$. Together this gives the identification in (\ref{eq:homeo1}).

Now let's deal with (\ref{eq:homeo2}). We have $O(2)/D_4=SO(2)/\langle R_\pi\rangle$. Under this identification the non-trivial element of $D_8/D_4\cong \Z/2$ acts as multiplication by $R_{\pi/2}$ on $SO(2)/\langle R_{\pi}\rangle$. We now have a homeomorphism
\[
SO(2)/\langle R_{\pi}\rangle \times_{\Z/2} SO(2)/\langle R_{\pi}\rangle\stackrel{\cong}{\longrightarrow} SO(2)/\langle R_{\pi}\rangle\times SO(2)/\langle R_{\pi}\rangle
\]
given by $[R_{\alpha}\langle R_{\pi}\rangle, R_{\beta}\langle R_{\pi}\rangle]\mapsto (R_{\alpha+\beta}\langle R_{\pi}\rangle, R_{\alpha-\beta}\langle R_{\pi}\rangle)$. Then both factors of $S^1$ in (\ref{eq:homeo2}) are circles parametrizing rotations modulo $\pi$.

We can now describe the map induced by $s$ from the right hand side of (\ref{eq:homeo1}) to the right hand side of (\ref{eq:homeo2}). Denote this map by $s'$. Let us write $S^0=\{0,1\}$ and $S^1\times (S^1\sqcup S^1)\cong S^0\times S^1\times S^1$. Chasing through the maps, we find that $s': S^0\times S^1\times S^1\rightarrow S^1\times S^1$ is given by
\begin{alignat*}{2}
&(0,R_{\alpha}\langle R_{\pi/2}\rangle, R_{\beta}\langle R_{\pi}\rangle) &&\longmapsto (R_{2\alpha+\beta}\langle R_{\pi}\rangle, R_{-\beta}\langle R_{\pi}\rangle)\\
&(1,R_{\alpha}\langle R_{\pi/2}\rangle, r_{\beta\textnormal{ mod }\pi/2}) &&\longmapsto (R_{2\beta}\langle R_{\pi}\rangle, R_{2(\alpha-\beta)}\langle R_{\pi}\rangle)\,.
\end{alignat*}
Equivalently, if we identify $S^1=\{\lambda \in \mathbb{C}\,|\, |\lambda|=1\}$, then $s'$ sends
\begin{alignat*}{2}
& (0,\lambda,\mu) && \longmapsto (\lambda \mu,\bar{\mu}) \\
& (1,\lambda,\mu) && \longmapsto (\mu,\lambda\bar{\mu})\, .
\end{alignat*}

Both component maps are homeomorphisms of the torus. By precomposing with their inverses we make $s'$ the folding map. That is, let
\begin{equation} \label{eq:homeo4}
\sigma: S^0\times S^1\times S^1\stackrel{\cong}{\longrightarrow} S^0\times S^1\times S^1
\end{equation}
be the map $(0,\lambda,\mu)\mapsto (0,\lambda\mu,\bar{\mu})$ and $(1,\lambda,\mu)\mapsto (1,\lambda\mu, \lambda)$. Then $s'\circ \sigma=\pi_2$ is the projection onto $S^1\times S^1$. Moreover, we can identify $O(2)/SO(2)$ with $S^0$ and the map $S^0\times S^1\times S^1 \rightarrow S^0$ induced by $r$ is simply the projection $\pi_1$ onto the first factor. Therefore, diagram (\ref{dgr:ecomo2}) is equivalent to
\begin{equation} \label{dgr:ecomo22}
\xymatrix{
S^0\times (S^1\times S^1) \ar[d]^-{\pi_1} \ar[r]^-{\pi_2}	& S^1\times S^1 \ar[d] \\
S^0 \ar[r]									& \Ecom O(2)\,,
}
\end{equation}
which gives
\[
\Ecom O(2)\simeq S^0\ast (S^1\times S^1)\simeq \Sigma(S^1\times S^1)\simeq S^2\vee S^2\vee S^3\, ,
\]
as claimed.
\end{proof}

\begin{remark}\label{rem:pi_nBcomO2}
  The theorem of Hilton \cite[Section XI.6]{Whitehead} allows us to compute the homotopy groups of $S^2\vee S^2\vee S^3$ in terms of homotopy groups of spheres. Thus, combining our Theorem \ref{thm:ecomo2} with the isomorphism of graded groups $\pi_\ast(\Bcom O(2))\cong \pi_\ast(\Ecom O(2))\oplus \pi_\ast(BO(2))$ from \cite[Theorem 6.3]{Ad5}, we obtain the following low-dimensional homotopy groups of $\Bcom O(2)$:
  \begin{center} {\def\arraystretch{1.2}\tabcolsep=3pt
      \begin{minipage}[t]{0.25\linewidth}
        \begin{tabular}{ll}
          $n$ & $\pi_n(\Bcom O(2))$     \\ \hline
          $1$ & $\Z/2$                  \\
          $2$ & $\Z^3$                  \\
          $3$ & $\Z^4$                  \\
          $4$ & $\Z^4\oplus (\Z/2)^4$   \\
          $5$ & $\Z^7\oplus (\Z/2)^8$ \\
        \end{tabular}
      \end{minipage}
      \quad
      \begin{minipage}[t]{0.65\linewidth}
        \begin{tabular}{ll}
          $n$ & $\pi_n(\Bcom O(2))$                       \\ \hline
          $6$ & $\Z^{16}\oplus (\Z/2)^{11}\oplus (\Z/{12})^4$ \\
          $7$ & $\Z^{34}\oplus (\Z/2)^{27}\oplus (\Z/{12})^4$   \\ 
          $8$ & $\Z^{68}\oplus (\Z/2)^{58}\oplus (\Z/{24})^7$   \\ 
          $9$ & $\Z^{140}\oplus (\Z/2)^{113}\oplus (\Z/3)^4 \oplus (\Z/{24})^{16}$ \\ 
          $10$& $\Z^{308}\oplus (\Z/2)^{215}\oplus (\Z/3)^4 \oplus (\Z/{15})^4 \oplus (\Z/{24})^{34}$ \\ 
        \end{tabular}
      \end{minipage}
    }
  \end{center}
\end{remark}\bigskip

The homotopy pushout (\ref{dgr:ecomo22}) makes it easy to describe the monodromy action for the homotopy fibration $\Ecom O(2)\rightarrow \Bcom O(2)\rightarrow BO(2)$. First, we describe the monodromy action for the homotopy fiber sequence of Lemma \ref{lem:hofibconj}.

\begin{lemma} \label{lem:mondrom}
Suppose that $g\in G$ represents an element of $\pi_1(BG)\cong \pi_0(G)$. The monodromy action in the homotopy fiber sequence of Lemma \ref{lem:hofibconj} is determined by the $G$-action
\[
[x,yH]\longmapsto [x,x^{-1}gxyH]
\]
for $[x,yH]\in G\times_{N(H)}G/H$.
\end{lemma}
\begin{proof}
A loop in the bar construction $BG$ representing the specified class in $\pi_1(BG)$ is given by the image of the $1$-simplex $\{g\}\times \Delta^1$ in the geometric realization. A lift of this loop under the fibration $\tilde{\psi}: G\times_{N(H)} B(G,G/H)\rightarrow BG$ with initial point $[x, yH]\in G\times_{N(H)}G/H$ is given by the image of $(x, \{(x^{-1}g x, yH)\}\times \Delta^1)$ in $G\times_{N(H)} B(G,G/H)$. To compute the end point of this lift we have to use the last face map in the bar construction $B(G,G/H)$, which multiplies $x^{-1}gx$ onto $yH$ from the left. Thus the end point of the lift is the point $[x, x^{-1}gxyH] \in G\times_{N(H)}G/H$.
\end{proof}

\begin{corollary} \label{cor:mondr}
  For any coefficient ring $R$, the monodromy representation of the fundamental group  $\pi_1(BO(2))=\Z/2$ on $H^2(\Ecom O(2);R)=R\oplus R$ is given by $(r_1,r_2)\mapsto (-r_2,-r_1)$, and it is trivial on $H^3(\Ecom O(2);R)=R$.
\end{corollary}
\begin{proof}
Each corner of diagram (\ref{dgr:ecomo22}) is the homotopy fiber of a map into $BO(2)$. We determine the monodromy action of $\pi_1(BO(2))$ on $S^0$, $S^0\times (S^1\times S^1)$ and $S^1\times S^1$ in order to obtain the action on $\Ecom O(2)$. The non-trivial element of $\pi_1(BO(2))$ is represented by the reflection $r_0\in O(2)$. The monodromy action on $S^0=O(2)/SO(2)$ is given by left-multiplication by $r_0$, which swaps the two points of $S^0$.

To determine the action on $S^0\times (S^1\times S^1)$ we use Lemma \ref{lem:mondrom}, which asserts that the action on $O(2)\times_{D_8} O(2)/\{\pm I\}$ is given by $[x,y\{\pm I\}]\mapsto [x, x^{-1}r_0 x \{\pm I\}]$. Upon identifying $O(2)\times_{D_8} O(2)/\{\pm I\}$ with $O(2)/D_8\times O(2)/\{\pm I\}$ (as in the proof of Lemma \ref{lem:poEcomG}) this corresponds to left multiplication by $r_0$ on $O(2)/\{\pm I\}$ while leaving $O(2)/D_8$ fixed. Thus, under the identification (\ref{eq:homeo1}), the action is given by
\begin{alignat*}{3}
& (0, R_{\alpha}\langle R_{\pi/2}\rangle, R_{\beta}\langle R_{\pi}\rangle) &&\longmapsto (1, R_{\alpha}\langle R_{\pi/2}\rangle, r_0R_{\beta}\langle R_{\pi}\rangle) &&=(1, R_{\alpha}\langle R_{\pi/2}\rangle, r_{-\beta/2\textnormal{ mod }\pi/2})\\
& (1, R_{\alpha}\langle R_{\pi/2}\rangle, r_{\beta\textnormal{ mod }\pi/2}) &&\longmapsto (0, R_{\alpha}\langle R_{\pi/2}\rangle, r_0r_{\beta\textnormal{ mod }\pi/2})&&=(0, R_{\alpha}\langle R_{\pi/2}\rangle, R_{-2\beta}\langle R_{\pi}\rangle)\, .
\end{alignat*}
Equivalently, this is the map $\varepsilon: S^0\times (S^1\times S^1)\stackrel{\cong}{\longrightarrow} S^0\times (S^1\times S^1)$ given by
\begin{alignat*}{2}
&(0,\lambda,\mu)&&\longmapsto (1, \lambda,\bar{\mu}) \\
&(1,\lambda,\mu)&&\longmapsto (0,\lambda,\bar{\mu})\,.
\end{alignat*}
To obtain the action in diagram (\ref{dgr:ecomo22}) we must conjugate $\varepsilon$ by the homeomorphism $\sigma$ (\ref{eq:homeo4}). It is easy to check that $\sigma^{-1}\circ \varepsilon\circ \sigma$ is given by $(0,\lambda,\mu)\mapsto (1,\mu,\lambda)$ and $(1,\lambda,\mu)\mapsto (0,\mu,\lambda)$.

In a similar way, we use Lemma \ref{lem:mondrom} to describe the monodromy action on $S^1\times S^1$. We leave it to the reader to check that the action is given by $(\lambda,\mu)\mapsto (\mu,\lambda)$. Therefore, the monodromy action in diagram (\ref{dgr:ecomo22}) is given by simultaneously swapping the two points in $S^0$ and the two factors in $S^1\times S^1$. On the homotopy pushout $\Ecom O(2)\simeq \Sigma(S^1\times S^1)$ this becomes reflection of the suspension coordinate and simultaneously swapping the factors of $S^1$. On $H^\ast(\Ecom O(2);R)$ this induces the monodromy representation stated in the corollary.
\end{proof}

\begin{remark}
  In \cite{Ad5} the authors show that the homotopy fibre sequence $\Omega \Ecom G\rightarrow \Omega \Bcom G\rightarrow \Omega BG$ is split for any compact Lie group $G$. As a consequence, there is a homotopy equivalence $\Omega \Bcom G\simeq \Omega \Ecom G\times \Omega BG$. They also show that the alternating group $A_5$ is an example for which this splitting does not deloop. More generally, it was proved in \cite[Lemma 1.2.5]{Grthesis} that $BG$ is not a retract up to homotopy of $\Bcom G$ (and therefore the splitting does not deloop) whenever $G$ is a non-abelian discrete group. For compact connected Lie groups $G$ it was shown in \cite[Theorem 1.2.2]{Grthesis} that $BG$ is not a retract up to homotopy of $\Bcom G_\one$. However, none of these results applies to $O(2)$!

  The calculations in this paper do show that the splitting of $\Omega \Bcom G$ does not deloop for $G=O(2)$. Indeed, we can show this in two different ways: it is a consequence of Corollary \ref{cor:mondr} and it is also a consequence of the cohomology calculations for $\Bcom O(2)$.
\end{remark}

We finish this section with a description of the map $\Ecom O(2)\rightarrow \Bcom O(2)$ on cohomology, showing that the generators $r\in H^2(\Bcom O(2);\Z)$ and $s\in H^3(\Bcom O(2);\F_2)$ can also be identified with monodromy invariants in the cohomology of $\Ecom O(2)$.\bigskip

Theorem \ref{thm:ecomo2} shows $H^\ast(\Ecom O(2);\Z)\cong\Z\langle \omega_1,\omega_2,\kappa\rangle$, where $\omega_i\in H^2(S^2;\Z)$ and $\kappa\in H^3(S^3;\Z)$ are generators. By Corollary \ref{cor:mondr}, the monodromy action on these generators is $\omega_1\mapsto -\omega_2$ and $\omega_2\mapsto -\omega_1$ while $\kappa$ remains fixed. Thus, the monodromy invariants are generated by $\omega_1-\omega_2$ and $\kappa$.

\begin{corollary}
In integral cohomology, the map $\Ecom O(2)\rightarrow \Bcom O(2)$ sends $r\mapsto \pm (\omega_1-\omega_2)$, while for $\F_2$-coefficients it sends $\bar{r}\mapsto \bar{\omega}_1+\bar{\omega}_2$ and $s\mapsto \bar{\kappa}$. All other classes are mapped to zero.
\end{corollary}
\begin{proof}
The kernel of the homomorphism $H^\ast(\Bcom O(2);\Z)\rightarrow H^\ast(\Ecom O(2);\Z)$ contains $W_1$, $W_2$ and $p_1$ (as they are pulled back from the base) as well as $b_1$, $b_2$ and $b_3$ (for degree reasons). That $r$ is mapped to the monodromy invariant $\pm(\omega_1-\omega_2)$ is most easily seen from a portion of the Serre spectral sequence associated to the homotopy fibration $\Ecom O(2)\rightarrow \Bcom O(2)\rightarrow BO(2)$. The left hand picture in the following figure depicts the relevant part of the $E_2$-page where $\times$ and $\bullet$ indicate a copy of $\Z$ and $\Z/2$, respectively.
\begin{figure}
  \centering
\begin{tikzpicture}
  \matrix (m) [matrix of math nodes,
    nodes in empty cells,nodes={minimum width=5ex,
    minimum height=5ex,outer sep=-5pt},
    column sep=1ex,row sep=-2ex]{  
           	&		& 		& &		\\
         3	&\times	& 0	& \bullet & \bullet	\\
         2	&\times	& 0	& \times & 0	\\
         1	&0		& 0	& 0	& 0	\\
         0	&\times	& 0	& \bullet & \bullet	\\
\quad\strut &0  & 1      &  2  & 3    \strut \\};
%\draw[-stealth]  (-0.8,0.71) -- (0.8,0.26);
\node (d) at (0.44,0.18) {{\footnotesize $d_2$}};
\draw[-stealth]  (-0.8,0.24) -- (0.8,-0.19);
\draw[thick] (m-1-1.east) -- (m-6-1.east) ;
\draw[thick] (m-6-1.north) -- (m-6-5.north) ;
\end{tikzpicture}
\begin{tikzpicture}
  \matrix (m) [matrix of math nodes,
    nodes in empty cells,nodes={minimum width=5ex,
    minimum height=5ex,outer sep=-5pt},
    column sep=1ex,row sep=-2ex]{  
           	&		& 		& &		\\
         3	&\bullet	& \bullet	& \bullet \bullet & \bullet \bullet	\\
         2	&\bullet	& 0	& \bullet & 0	\\
         1	&0		& 0	& 0	& 0	\\
         0	&\bullet	& \bullet	& \bullet \bullet & \bullet \bullet	\\
\quad\strut &0  & 1      &  2  & 3    \strut \\};
\draw[-stealth]  (-0.8,0.71) -- (0.8,0.26);
\node (d) at (0.45,0.60) {{\footnotesize $d_2$}};
%\draw[-stealth]  (-0.8,0.24) -- (0.8,-0.19);
\draw[thick] (m-1-1.east) -- (m-6-1.east) ;
\draw[thick] (m-6-1.north) -- (m-6-5.north) ;
\end{tikzpicture}
\end{figure}
The differential drawn is $d_2$ originating from $\omega_1-\omega_2$ and is clearly zero. The corresponding differential $d_3$ will vanish too, since the cohomology of $BO(2)$ injects into the cohomology of $\Bcom O(2)$, so $\omega_1-\omega_2$ does not transgress. Thus, $\omega_1-\omega_2$ lies in the image of $H^\ast(\Bcom O(2);\Z)\rightarrow H^\ast(\Ecom O(2);\Z)$. It also follows that the homomorphism $H^\ast(\Bcom O(2);\F_2)\rightarrow H^\ast(\Ecom O(2);\F_2)$ maps $\bar{r}\mapsto \bar{\omega}_1+\bar{\omega}_2$. To see that $s$ is mapped to the top degree class $\bar{\kappa}$ we consider again a portion of the Serre spectral sequence, this time for $\F_2$-coefficients (the right hand picture in the figure above). The differential labelled $d_2$ originates from the class $\bar{\kappa}$. A simple rank argument shows that the differential must be zero. The next differential $d_3$ has trivial target, and $d_4(\bar{\kappa})$ must vanish, since the mod $2$ cohomology of $BO(2)$ injects into the mod $2$ cohomology of $\Bcom O(2)$. Thus, $\bar{\kappa}$ lies in the image of $H^\ast(\Bcom O(2);\F_2)\rightarrow H^\ast(\Ecom O(2);\F_2)$, and $s$ is the only class that can map to it.
\end{proof}

\section{The homotopy type of $\Ecom SU(2)$}

In this section we compute the homotopy type of $\Ecom SU(2)$,
starting from the homotopy pushout square  given in
Lemma \ref{lem:poEcomG}. For convenience, we'll give one letter names
to several groups: in this section, we will let $G = SU(2)$;
$Z = \{\pm I\}$, its center; $T \cong S^1$, the maximal torus in $G$;
$N$ the normalizer of the torus; and $W = N/T \cong \Z/2$, the Weyl
group. The homotopy pushout square for $\Ecom SU(2)$ then reads:

\[\xymatrix{
    G/N \times G/Z \ar[d] \ar[r]^{\sigma} & G/T \times_W G/T \ar[d] \\
    G/Z \ar[r] & \Ecom G \\
  }\]

The left vertical map is simply the projection and the top horizontal
map is given by $\sigma(gN, hZ) = [gT, hgT]$.

Now, using that $G = SU(2) \cong S^3$, $G/T \cong S^2$, our pushout
square for $\Ecom SU(2)$ looks like:

\begin{equation}
  \xymatrix{
    \RP^2 \times \RP^3 \ar[d] \ar[r]^{\sigma} & S^2 \times_W S^2 \ar[d] \\
    \RP^3 \ar[r] & \Ecom SU(2) \\
  }
  \label{dgr:ecomsu2}
\end{equation}

We'll use this square to compute the homotopy type of $\Ecom SU(2)$,
but the strategy will be somewhat indirect. We'll start by computing
the pushout restricted to a certain copy of $\RP^2$ inside of $\RP^3$.

\subsection{An easy piece of $\Ecom SU(2)$}

Fix an ``equator'' of $G=SU(2)\cong S^3$: that is, a subspace
$S \subset G$ (\emph{not} a sub\emph{group}!), homeomorphic to $S^2$
and closed under multiplication by the center $Z$. Then $S/Z$ (which
just means the space of all cosets $sZ \in G/Z$ for $s \in S$) is the
copy of $\RP^2$ inside of $G/Z \cong \RP^3$ we will focus on.

\begin{lemma}
  Let $X$ be obtained as the following homotopy pushout in which the
  top horizontal map is a restriction of the corresponding map in
  (\ref{dgr:ecomsu2}):
  \begin{equation}
    \xymatrix{
      \RP^2 \times \RP^2 \ar[d]_{\pi_2} \ar[r]^{\sigma_|} & S^2
      \times_W S^2 \ar[d] \\
      \RP^2 \ar[r] & X. \\
    }
    \label{dgr:s4copy}
  \end{equation}
  Then there is homotopy cofiber sequence
  \[ X \to \Ecom SU(2) \to \Sigma^4 \RP^2,\]
  and a homotopy equivalence $X \simeq S^4$.\label{lem:s4}
\end{lemma}

\begin{proof}
  Consider the following commutative diagram, whose rows are minor
  modifications (for convenience) of the pushouts we have for $X$ and
  $\Ecom SU(2)$:
  \begin{equation}
    \xymatrix{
      \RP^2_+ \ar[d]^{\iota_+} & \RP^2_+ \wedge \RP^2_+ \ar[l]
      \ar[r]^{\sigma_|} \ar[d]^{\id \wedge \iota_+} & S^2 \times_W S^2
      \ar[d]^{\id} \\
      \RP^3_+ & \RP^2_+ \wedge \RP^3_+ \ar[l] \ar[r]^{\sigma} & S^2 \times_W
      S^2 \\
    }\label{dgr:XEcomG}
  \end{equation}
  Here $\iota \colon \RP^2 \to \RP^3$ is the inclusion $S/Z \to G/Z$
  and $A_+$ denotes $A$ with a disjoint basepoint added.

  Despite the changes, the homotopy pushouts of the rows are
  still $X$ and $\Ecom SU(2)$, respectively. Take the top row, for
  example. The middle space is
  $\RP^2_+ \wedge \RP^2_+ \cong (\RP^2 \times \RP^2)_+$, which makes
  the top row of the form $A_+ \leftarrow B_+ \to C$ and such a span
  has the same homotopy pushout as $A \leftarrow B \to C$.

  Thus the diagram induces a map $j \colon X \to \Ecom SU(2)$ and we
  can compute $\hocofib(j)$ as the homotopy pushout of the homotopy
  cofibers of the vertical maps in (\ref{dgr:XEcomG}). Let's compute
  those cofibers: first, $\hocofib(\iota_+) \simeq S^3$; second, since
  a functor of the form $A \wedge (-)$ preserves homotopy cofibers,
  the middle vertical map has cofiber $\RP^2_+ \wedge S^3$. Also,
  since the left-pointing horizontal maps in (\ref{dgr:XEcomG}) can be
  written as $q_+ \wedge \id_{\RP^i}$ ($i = 2, 3$) where $q$ is the
  unique map $\RP^2 \to \ast$, we see that the induced map between the
  cofibers of the first two vertical maps is $q_+ \wedge \id_{S^3}$.

  Since the last vertical map is the identity, its cofiber is $\ast$.
  This means the homotopy pushout of the homotopy cofibers reduces to
  $\hocofib(q_+ \wedge \id_{S^3} \colon \RP^2_+ \wedge S^3 \to S^3)
  \simeq \hocofib(q_+) \wedge S^3 \simeq \Sigma \RP^2 \wedge S^3
  \simeq \Sigma^4 \RP^2$. This establishes the cofiber sequence in the
  claim.

  \smallskip
  
  To prove that $X \simeq S^4$, we will show that $X$ is a
  simply-connected integral homology sphere. That implies first, that
  $\pi_4 X \cong H_4(X; \Z) = \Z$ and second, that any choice of a
  generator $S^4 \to X$ of $\pi_4 X$ is an integral homology
  equivalence between simply-connected spaces and thus is an
  equivalence.

  To compute the homology of $X$, we use the two homotopy pushouts we
  have involving $X$. First, from Mayer-Vietoris sequence for the
  defining homotopy pushout we can readily compute that $H_k(X) = 0$
  for $k>4$. Indeed, $S^2 \times_W S^2$ and $\RP^2 \times \RP^2$ are
  both $4$-manifolds, and the latter is additionally non-orientable,
  so that $H_k(S^2 \times_W S^2) = 0$ for $k \ge 5$ and $H_k(\RP^2
  \times \RP^2) = 0$ even for $k \ge 4$.

  Next we use the Mayer-Vietoris for the cofiber sequence we just
  established. Since $\tilde{H}_k(\Sigma^4 \RP^2) \neq 0$ only for
  $k=5$, we learn that the map $X \to \Ecom SU(2)$ is an integral
  homology isomorphism except possibly in degrees $4$ and $5$. The
  portion of the sequence for those degrees is as follows:
  \[ 0 \to H_5 X \to H_5 \Ecom SU(2) \to H_5\Sigma^4\RP^2 \to H_4X
    \to H_4\Ecom SU(2) \to 0.\]

  Since we found that $H_5X = 0$, the homomorphism
  $H_5 \Ecom SU(2) \to H_5\Sigma^4\RP^2$ is injective. Now, we
  calculated the integral cohomology of $\Ecom SU(2)$ previously in
  Corollary \ref{cor:cohEcomSU2} and found it was
  $(\Z, 0, 0, \Z, 0, \Z/2)$, so by the universal coefficient theorem,
  the homology is $(\Z, 0, 0, \Z, \Z/2)$. Thus both the domain and
  codomain of the injective homomorphism
  $H_5 \Ecom SU(2) \to H_5\Sigma^4\RP^2$ are $\Z/2$, which forces it
  to be an isomorphism. This in turn implies that
  $H_4X \to H_4 \Ecom SU(2)$ is also an isomorphism. We conclude that
  $X$ is an integral homology $4$-sphere.

  Finally, to show $X$ is simply-connected, we use the Seifert-van
  Kampen theorem. Since $S^2 \times S^2$ is the universal cover of
  $S^2 \times_W S^2$, we have $\pi_1(S^2 \times_W S^2) \cong W$. So the
  Seifert-van Kampen theorem says that we have the following pushout
  in the category of groups:
  \[\xymatrix{
      \Z/2 \times \Z/2 \ar[r]^-{\sigma_{|*}} \ar[d]_{\pi_2} & \Z/2
      \ar[d] \\ 
      \Z/2 \ar[r] & \pi_1(X). \\
    }\]

  It's not too hard to check that $\sigma_{|*}$ is the projection
  $\pi_1$. Indeed, we are looking for the effect of the map
  $\sigma_|$, given by $\sigma_|(gN, sZ) = [gT, sgT]$ on the
  fundamental group. Since $G/T \to G/N$ is the universal cover of
  $G/N$, a generator of $\pi_1(G/N)$ is a loop that lifts in $G/T$ to
  path connecting $T$ to $n_0T$, where $n_0T$ is a generator of the Weyl
  group $N/T$. Under $\sigma_|$ this loop will go to one that lifts in
  $G/T \times G/T$ to a path connecting $(T,T)$ with
  $(n_0T, n_0T) = n_0T \cdot (T,T)$ and thus $\sigma_{|*}$ is non-zero on
  the generator of $\pi_1(G/N)$. For the second factor note that the
  restriction of $\sigma_|$ to $S/Z$ is given by $sZ \mapsto [IT, sT]$
  and so factors through the map $G/T \to G/T \times_W G/T$, $gT
  \mapsto [IT, gT]$. Since $G/T \cong S^2$, we see that the
  restriction $\sigma_{|S/Z}$ is null-homotopic.
  
  Knowing the two homomorphisms, it is straightforward to compute that
  $\pi_1(X) \cong 1$, which concludes the proof.
\end{proof}

\subsection{Writing $\Ecom SU(2)$ as a homotopy cofiber}

To obtain the full pushout computing $\Ecom SU(2)$ from this, we need
to attach a 3-cell to $\RP^2$ to obtain $\RP^3$. Consider the
following diagram, a span of spans:

\begin{equation}
  \xymatrix{
    \ast & \RP^2 \ar[l] \ar[r]^{\id} & \RP^2 \\
    S^2 \ar[u] \ar[d]_{p} &
    \RP^2 \times S^2 \ar[l]_-{\pi_2} \ar[d]^{\id \times p}
    \ar[r]^{\pi_1} \ar[u] & \RP^2 \ar[u]_{\id} \ar[d]^{\Delta_{/W}} \\
    \RP^2 & \RP^2 \times \RP^2  \ar[l]_-{\pi_2} \ar[r]^{\sigma_|} &
    S^2 \times_W S^2. \\
  }
  \label{dgr:ecomsu2=hocofib}
\end{equation}

Here $p \colon S^2 \to \RP^2$ is the canonical quotient map
$S \to S/Z$, and $\Delta_{/W}$ denotes the map induced on $W$-orbits
by the diagonal $\Delta \colon S^2 \to S^2 \times S^2$. All squares
except for the bottom right one commute. Perhaps surprisingly, the
bottom right square commutes up to homotopy, as we shall prove next.

\begin{lemma}\label{lem:sqcommhtpy}
  The bottom right square of (\ref{dgr:ecomsu2=hocofib}) commutes up
  to homotopy.
\end{lemma}

\begin{proof}
  The two composite maps $\RP^2 \times S^2 \to S^2 \times_W S^2$, or
  better, the two maps $G/N \times S \to G/T \times_W G/T$, are given
  by $(gN, s) \mapsto [gT, gT]$ for the top route, and
  $(gN, s) \mapsto [gT, sgT]$ for the bottom route. Now, the inclusion
  $S \hookrightarrow G$ is null-homotopic, and picking a homotopy
  $H \colon S \times [0,1] \to G$ with $H(s,0) \equiv I$ and
  $H(s,1) = s$, we can easily define a homotopy
  $G/N \times S \times [0,1] \to S^2 \times_W S^2$, given by
  $(gN, s, t) \mapsto [gT, H(s,t) g T]$, between those two maps.
\end{proof}

Analyzing the diagram (\ref{dgr:ecomsu2=hocofib}) we will show the
following.

\begin{lemma}
  There is some map $\theta \colon \Sigma^3 \RP^2 \to S^4$ such that
  $\Ecom SU(2) \simeq \hocofib(\theta)$.
\end{lemma}

\begin{proof}
  We can compute the homotopy colimit of (\ref{dgr:ecomsu2=hocofib})
  in two different ways. Comparing the two will establish the result.

  \smallskip
  
    \paragraph*{Starting by columns:}
    \begin{itemize}
    \item The homotopy pushout of the left column is $\RP^3$, since
      $p$ is the attaching map of the $3$-cell of $\RP^3$.
    \item The middle column is just the image of the left column under
      the functor $\RP^2 \times (-)$, which preserves homotopy
      colimits, so the homotopy pushout of the middle column is
      $\RP^2 \times \RP^3$ and the induced map between the homotopy
      pushouts, $\RP^3 \leftarrow \RP^2 \times \RP^3$, is just the
      projection.
    \item Because one of the maps in the rightmost column is the
      identity, the homotopy pushout is the third space, namely,
      $S^2 \times_W S^2$.
    \end{itemize}

  So taking pushouts by columns gives a span of the form
  \[\RP^3 \xleftarrow{\pi_2} \RP^2 \times \RP^3 \to S^2 \times_W
    S^2.\] Once we show that the right-pointing map is homotopic to
  the map $\sigma$ in the pushout (\ref{dgr:ecomsu2}) for
  $\Ecom SU(2)$, we will have proven that the homotopy colimit of
  (\ref{dgr:ecomsu2=hocofib}) is $\Ecom SU(2)$.

  To prove that however, it is not enough to know Lemma
  \ref{lem:sqcommhtpy} as is: not every homotopy making the lower right
  square commute will induce the correct map between the homotopy
  pushouts of the middle and right column! We instead need to make a
  good choice of the null-homotopy $H$ of $S \hookrightarrow G$ used
  in the proof of Lemma \ref{lem:sqcommhtpy}. Notice that the space of
  such null-homotopies is $\Omega^3 S^3$, so there are a $\Z$'s worth
  of non-homotopic choices. Here's how we pick the right one: given a
  null-homotopy $H \colon S \times [0,1] \to G$ it induces a map
  $\bar{H} \colon C \to G/Z$, where $C = \hocofib(S \to S/Z)$. (Explicitly,
  say $C$ is given by $C = (S \times [0,1] \sqcup S/Z)/\!\!\sim$,
  where $\sim$ is the equivalence relation generated by
  $(s,0) \sim (s'',0)$ and $(s,1) \sim sZ$, then
  $\bar{H}(s,t) = H(s,t)Z$ and $\bar{H}(sZ) = sZ$.) Since
  $C \simeq \RP^3 \cong G/Z$ we can choose $H$ so that the induced
  $\bar{H}$ is an equivalence.
  
  To see that such a choice of $H$ induces the correct homotopy class
  of map, we factor the natural transformation between the last two
  columns of (\ref{dgr:ecomsu2=hocofib}) as follows:

  \[\xymatrix{
      \RP^2  \ar[r]^{\id}  & \RP^2 \ar[r]^{\id} & \RP^2 \\
      \RP^2 \times S^2 \ar[d]^{\id \times p} \ar[r]^{\pi_1} \ar[u] &
      \RP^2 \ar[r]^{\id} \ar[d]^{i} \ar[u]^{\id} &
      \RP^2 \ar[u]_{\id} \ar[d]^{\Delta_{/W}} \\
      \RP^2 \times \RP^2  \ar[r]^{\id \times j}  & \RP^2 \times \RP^3
      \ar[r]^{\sigma} & S^2 \times_W S^2, \\
    }\]
  where $i$ is the inclusion at the basepoint
  $IZ \in G/Z \cong \RP^3$. Notice that all squares here commute
  except for the lower left one which commutes up to homotopy,
  specifically up to the homotopy $\id_{\RP^2} \times (\pi \circ H)$,
  where $\pi \colon G \to G/Z$ is the canonical projection.

  In this new diagram take pushouts by columns to get
  \[\RP^2 \times C \xrightarrow{\id \times \bar{H}} \RP^2 \times \RP^3
  \xrightarrow{\sigma} S^2 \times_W S^2,\] whose composite is indeed
  equivalent to $\sigma$ as desired.

  \smallskip
  
    \paragraph*{Starting by rows:}
    \begin{itemize}
    \item Again because one of the maps in the top row is the
      identity, the homotopy pushout of the top row is the third
      space, that is, contractible.
  
    \item The homotopy pushout of the middle row is one of the
      standard ways of describing the join
      $S^2 * \RP^2 \simeq \Sigma^3 \RP^2$.
  
    \item The bottom row is the span appearing in Lemma \ref{lem:s4},
      whose homotopy pushout was shown to be $S^4$.
    \end{itemize}

  So taking pushouts by columns produces a span
  $\ast \leftarrow \Sigma^3 \RP^2 \xrightarrow{\theta} S^4$, showing
  that the homotopy colimit of (\ref{dgr:ecomsu2=hocofib}) is
  equivalent to the homotopy cofiber of whatever map $\theta$ happens
  to be induced between the pushouts of the middle and bottom row.
\end{proof}

\subsection{The two candidates for $\Ecom SU(2)$}

Unfortunately, the description we have of the map $\theta$ is very
indirect, as it is depends both on the homotopy chosen in Lemma
\ref{lem:sqcommhtpy} and the homotopy equivalence that was established
by calculational means in Lemma \ref{lem:s4}. Let's see what the
possibilities for $\theta$ are, that is, let's determine the group of
pointed homotopy classes of maps $[\Sigma^3 \RP^2, S^4]$. We have the
cofiber sequence $S^1 \xrightarrow{2} S^1 \to \RP^2$ where $2$ denotes
a degree $2$ map (more generally, we will use an integer to denote a
map of that degree on a sphere). Suspending thrice produces the
cofiber sequence $S^4 \xrightarrow{2} S^4 \to \Sigma^3 \RP^2$ and then
applying the pointed mapping space functor $\Map_*(-, S^4)$ turns it
into a fiber sequence:
\[\Map_*(\Sigma^3 \RP^2, S^4) \to \Map_*(S^4, S^4) \xrightarrow{2}
  \Map_*(S^4, S^4).\]
The long exact sequence of homotopy groups of that
fiber sequences reads in part:
\[ \pi_5(S^4) \xrightarrow{2} \pi_5(S^4) \xrightarrow{\partial}
  [\Sigma^3 \RP^2, S^4] \to \pi_4(S^4) \xrightarrow{2} \pi_4(S^4).\]

Since $\pi_5(S^4) \cong \Z/2$ and $\pi_4(S^4) \cong \Z$, we obtain
that $\partial$ is an isomorphism and thus
$[\Sigma^3 \RP^2, S^4] \cong \Z/2$. So there are exactly two possible
homotopy classes for the map $\theta$, and thus two candidates for
$\Ecom SU(2)$. To decide which one is correct, we'll need an alternate
description of the candidates, which we will give in terms of
$\vartheta := \partial^{-1}(\theta) \in \pi_5(S^4)$. Recall that
$\partial$ can be described in terms of the Puppe sequence
\[ S^4 \xrightarrow{2} S^4 \to \Sigma^3\RP^2 \xrightarrow{\delta} S^5
  \xrightarrow{-2} S^5 \to \cdots, \]
namely, $\partial = (-) \circ \delta$, so that $\theta = \vartheta
\circ \delta$.

Consider then the following diagram all of whose squares are homotopy
pushouts:

\[\xymatrix{
    \Sigma^3 \RP^2 \ar[r]^{\delta} \ar[d] & S^5 \ar[r]^{\vartheta}
    \ar[d]^{-2} & S^4 \ar[d] \\
    \ast \ar[r] & S^5 \ar[r] \ar[d]^{-1} & \hocofib(\theta)
    \ar[d]^{\simeq} \\
    & S^5 \ar[r] & \mathrm{pushout(2,\vartheta)}
  }\]
Since $-1$ is an equivalence, so is the map labeled $\simeq$.

This alternate description of $\Ecom SU(2)$ as the pushout of
$S^5 \xleftarrow{2} S^5 \xrightarrow{\vartheta} S^4$ is the one we'll
need later. But we can also make the CW-structure more explicit by
using the standard fact that the homotopy pushout square of pointed
spaces on the left can be rewritten as the homotopy pushout square on
the right (where $\nabla \colon A \vee A \to A$ is the fold map, given by
the identity on each wedge summand):

\[ \xymatrix{
    A \ar[r]^{b} \ar[d]_{c} & B \ar[d] \\
    C \ar[r] & D
  }
  \quad \quad
  \xymatrix{
    A \vee A \ar[d]_{\nabla} \ar[r]^{b \vee c} & B \vee C \ar[d] \\
    A \ar[r] & D
  }\]

Using this we can produce the following diagram both of whose squares
are homotopy pushouts:

\[\xymatrix{
    S^5 \ar[r]^-{\rho} \ar[d] & S^5 \vee S^5 \ar[r]^{-2 \vee \vartheta}
  \ar[d]_{\nabla} & S^5 \vee S^4 \ar[d] \\
  \ast \ar[r] & S^5 \ar[r] & \hocofib(\theta)
}\]

Here $\rho$ is a variant of the pinch map which is degree $-1$ onto
the first wedge summand, and degree $1$ onto the second. The fact that
the outer rectangle is a homotopy pushout shows that
$\hocofib(\theta)$ is homotopy equivalent to the CW-complex with
$5$-skeleton $S^5 \vee S^4$ with a $6$-cell attached along the map
$S^5 \xrightarrow{\mathrm{pinch}} S^5 \vee S^5 \xrightarrow{2 \vee
  \vartheta} S^5 \vee S^4$.

Using any of the three descriptions of the candidates it is easy to
see that if $\theta$, or equivalently $\vartheta$, is null-homotopic,
we get $\Ecom SU(2) \simeq \Sigma^4 \RP^2 \vee S^4$, which is what we
will show happens.

\subsection{Resolving the ambiguity}

To decide which of the two candidates for $\Ecom SU(2)$ is the correct
one, we will compute the action of the Steenrod algebra on the mod 2
cohomology; which is why in this section all cohomology groups
mentioned are with $\F_2$ coefficients. Notice first that both
candidates have the same cohomology (even integrally): the description
of the CW-structure shows the two possible attaching maps of the
$6$-cell, $S^5 \to S^5 \vee S^4$, differ only in how they map to the
$S^4$ wedge summand, which is in the $4$-skeleton of the CW-complex
and thus invisible to the cellular chain complex.

The candidate with $\theta$ null-homotopic, namely $\Sigma^4 \RP^2
\vee S^4$, is easy: the only non-zero Steenrod square is $\Sq^1 \colon H^5
\to H^6$, which comes from $\Sq^1 \colon H^1(\RP^2) \to  H^2(\RP^2)$. In
particular, $\Sq^2$ is zero for this candidate.

We will now show that for the other candidate, the one with non-null
$\theta$, the Steenrod square $\Sq^2 \colon H^4 \to H^6$ is non-zero. For
this we will use the description as the pushout of
$S^5 \xleftarrow{2} S^5 \xrightarrow{\vartheta} S^4$. Recall that the
non-zero element of $\pi_5(S^4)$ is given by the double suspension of
the Hopf fibration $\eta \colon S^3 \to S^2$, so this candidate is actually
$\Sigma^2 Y$, where $Y$ is the homotopy pushout of
$S^3 \xleftarrow{2} S^3 \xrightarrow{\eta} S^2$. It is enough to show
that $\Sq^2 \colon H^2(Y) \to H^4(Y)$ is non-zero. Consider the following
diagram both of whose squares are homotopy pushouts:

\[\xymatrix{
    S^3 \ar[d]_{\eta} \ar[r]^{2} & S^3 \ar[d] \ar[r] & \ast \ar[d] \\
    S^2 \ar[r] & Y \ar[r] & \CP^2 \\
  }\]

The rightmost space on the bottom is $\CP^2$, because $\eta$ is the
attaching map for the $4$-cell in the standard CW-structure of
$\CP^2$. Using the Mayer-Vietoris sequence first for the left-hand
square, we see that the map $S^3 \to Y$ is an isomorphism on $H^3$;
and then for the right-hand square, we learn that the map
$Y \to \CP^2$ is an isomorphism on $H^4$. Since $\Sq^2$ is non-zero
for $\CP^2$, we conclude it also non-zero for $Y$.

So now all that is left to do is calculate $\Sq^2$ on $H^4(\Ecom
SU(2))$: if it is zero then $\Ecom SU(2) \simeq \Sigma^4 \RP^2 \vee
S^4$; if it is non-zero then $\Ecom SU(2) \simeq \Sigma^2 Y$.

To compute $\Sq^2$ for $\Ecom SU(2)$, we'll go back to the homotopy
pushout square (\ref{dgr:ecomsu2}) and use the fact that Steenrod
squares commute with the connecting homomorphism in the Mayer-Vietoris
sequence. A small portion of the Mayer-Vietoris sequence reads:
\[ H^3(\RP^3) \oplus H^3(S^2 \times_W S^2) \to H^3(\RP^2 \times \RP^3)
  \xrightarrow{\delta} H^4(\Ecom SU(2)). \]

We can compute $H^*(S^2 \times_W S^2)$ by using the Serre spectral
sequence for the fiber sequence
\[S^2 \to S^2 \times_W S^2 \to S^2/W \cong \RP^2.\]
The fundamental group of the base, $\pi_1(\RP^2) \cong \Z/2$, acts by
antipodes on the fiber $S^2$, but with $\F_2$ coefficients the
corresponding action on $H^*(S^2)$ is trivial! The $E_2$-page of the
spectral sequence is thus $E_2^{p,q} = H^p(\RP^2) \otimes H^q(S^2)$.
There is no room for differentials and thus, as graded vector spaces,
$H^*(S^2 \times_W S^2) \cong H^*(\RP^2) \otimes H^*(S^2)$. In
particular, $H^3(S^2 \times_W S^2) \cong \F_2$.

Now, since both $H^3(\RP^3)$ and $H^3(S^2 \times_W S^2)$ are
$1$-dimensional over $\F_2$, and $H^3(\RP^2 \times \RP^3)$ is
$3$-dimensional, there must be some element
$x \in H^3(\RP^2 \times \RP^3)$ outside the kernel of the connecting
homomorphism $\delta$. This means that $\delta x$ is the sole non-zero
class in $H^4(\Ecom SU(2))$, and we have
$\Sq^2(\delta x) = \delta(\Sq^2 x)$ which is zero because it can
easily be checked that $\Sq^2$ is zero on all of
$H^3(\RP^2 \times \RP^3)$.

We conclude that $\Sq^2$ is zero on $H^4(\Ecom SU(2))$, and by the
discussion above, this proves the main result on $\Ecom SU(2)$ stated
in the introduction.

\renewcommand{\introthmref}{thm:ecomsu2}
\begin{introthm*}
There is a homotopy equivalence $\Ecom SU(2) \simeq S^4 \vee \Sigma^4 \RP^2$.
\end{introthm*}

\begin{remark}\label{rem:pi_nBcomSU2}
  As in Remark \ref{rem:pi_nBcomO2}, we can use our calculation of $\Ecom SU(2)$ to compute the first few homotopy groups of $\Bcom SU(2)$. The Hilton-Milnor theorem \cite[Section XI.6]{Whitehead}, readily shows that $\Omega \Ecom SU(2) \simeq \Omega(S^4 \times \Sigma^4 \RP^2 \times \Sigma^7 \RP^2 \times A)$, where the space $A$ is 10-connected. Many homotopy groups of $\Sigma^k \RP^2$ have been computed in \cite{JieWu}. Using those we obtain the following table (recall that $\Bcom SU(2)$ is 3-connected):

  \begin{center}
    \def\arraystretch{1.2}
    \begin{minipage}[t]{0.4\linewidth}
      \begin{tabular}[t]{ll}
        $n$ & $\pi_n(\Bcom SU(2))$ \\ \hline
        $4$ & $\Z^2$     \\
        $5$ & $(\Z/2)^3$ \\
        $6$ & $(\Z/2)^3$ \\
        $7$ & $\Z \oplus \Z/4 \oplus (\Z/{12})^2$ \\
      \end{tabular}
    \end{minipage}\quad
    \begin{minipage}[t]{0.4\linewidth}
      \begin{tabular}[t]{ll}
        $n$ & $\pi_n(\Bcom SU(2))$ \\ \hline
        $8$ & $(\Z/2)^6$ \\ 
        $9$ & $(\Z/2)^6$ \\ 
        $10$& $\Z/12 \oplus (\Z/24)^2$ \\
      \end{tabular}
    \end{minipage}
  \end{center}
\end{remark}

\begin{remark}
  The calculation of $\Ecom SU(2)$ yields $\Ecom SO(4)_\one$ as a corollary. According to Lemma \ref{lem:covhom}, if $\tilde{G} \to G$ is a homomorphism of compact connected Lie groups which is also a covering map, then $\Ecom \tilde{G}_\one \simeq \Ecom G_\one$. Since there is a double cover $SU(2) \times SU(2) \to SO(4)$, we see that \[\Ecom SO(4)_\one \simeq \Ecom \left(SU(2) \times SU(2)\right) \cong \Ecom SU(2) \times \Ecom SU(2) \simeq (S^4 \vee \Sigma^4 \RP^2)^2.\]
\end{remark}

\appendix\section{Computing kernels of ring homomorphisms in {\textsc{Singular}}}\label{app:singular}

We used the \textsc{Singular} computer algebra system \cite{Singular}
several times to compute generators for the kernel of a ring homomorphism
needed in a proof. This appendix shows those calculations.

\subsection{The kernel needed in the proof of Theorem \ref{thm:cru2}} \label{sec:singular1}

In the proof of Theorem \ref{thm:cru2} we needed to know the kernel of the
ring homomorphism
$f \colon \Z[c_1, c_2, y_1, y_2] \to \Z[t,\tilde{c}_1,\tilde{c}_2,c,U]/(t U, t \tilde{c}_1, t \tilde{c}_2,
t c, c^2, U^2, 2 U, c U, \tilde{c}_1 U)$ given by $f(c_1)=2t+\tilde{c}_1$,
$f(c_2)=t^2+\tilde{c}_2$, $f(y_1)=2c$ and $f(y_2)=\tilde{c}_1c$.

First we define the polynomial rings $P:=\Z[c_1,c_2,y_1,y_2]$ and $F:=\Z[t,U,a_1,a_2, c]$, the ideal of $F$ given by $I:=(t U, t a_1, t a_2, t c, c^2, U^2, 2 U, c U, a_1 U)$ and the quotient ring $R:=F/I$.

\begin{verbatim}
ring P = ZZ, (c1, c2, y1, y2), dp;
ring F = ZZ, (t, u, a1, a2, c), dp;
ideal I = t*U, t*a1, t*a2, t*c, c^2, U^2, 2*u, c*U, a1*U;
qring R = std(I);
\end{verbatim}

In \textsc{Singular} there is always a ``current ring'' and one defines homomorphisms \emph{into} it by specifying the domain and the list of images of the chosen generators for said domain:

\begin{verbatim}
map f = P, 2*t+a1, t^2+a2, 2*c, a1*c;
\end{verbatim}

To ask for the kernel of the homomorphism, one sets the current to be the domain of the homomorphism and uses the \texttt{kernel} function:

\begin{verbatim}
setring P;
kernel(R, f);
\end{verbatim}

\textsc{Singular} responds with a list of generators for the ideal $\ker(f)$:

\begin{verbatim}
_[1]=y2^2
_[2]=y1*y2
_[3]=y1^2
_[4]=c1*y1-2*y2
\end{verbatim}

\subsection{The kernel needed in the proof of Proposition \ref{prop:cohbcomo2part1}} \label{sec:singular2}

In the proof of Proposition \ref{prop:cohbcomo2part1} we considered the ring homomorphism
\[
f: \Z[W_1,W_2,p_1,r,b_1,b_2,b_3]\to \Z[a,t,\tilde{W}_1,\tilde{W}_2, x,\tilde{p}_1]/J
\]
where $J$ is the ideal generated by $at$, $a\tilde{W}_1$, $a\tilde{W}_2$, $ax$, $a\tilde{p}_1$, $2\tilde{W}_1$, $2\tilde{W}_2$, $2x$, $2\tilde{p}_1$, $t^2$, $t\tilde{W}_1$, $tx$, $x^2$, $\tilde{W}_2^2+\tilde{p}_1 \tilde{W}_1$ and $f$ is given by $f(W_1)=\tilde{W}_1$, $f(W_2)=\tilde{W}_2$, $f(p_1)=a^2+\tilde{p}_1$, $f(r)=2a$, $f(b_1)=t\tilde{W}_2$, $f(b_2)=x\tilde{W}_1$ and $f(b_3)=x\tilde{W}_2$.

The codomain of $f$ is a presentation for the cohomology ring $H^\ast(BSO(2)\vee (O(2)\times_{N(D_4)}BD_4);\Z)$. We needed to know the kernel of $f$, which in \textsc{Singular} can be computed as follows:

\begin{verbatim}
ring P = ZZ, (W1, W2, p1, r, b1, b2, b3), dp;
ring F = ZZ, (a,t,W1,W2,x,p1), dp;
ideal J = a*t, a*W1, a*W2, a*x, a*p1, 2*W1, 2*W2, 2*x, 2*p1, t^2, t*W1,
   t*x, x^2, W2^2+p1*W1;
qring R = std(J);
map f = P, W1, W2, a^2+p1, 2*a, t*W2, x*W1, x*W2;
setring P;
kernel(R, f);
\end{verbatim}

\textsc{Singular}'s response (in two columns to save space) is:\bigskip

\begin{minipage}{0.3\linewidth}
\begin{verbatim}
_[1]=2*b3
_[2]=2*b2
_[3]=2*b1
_[4]=2*W2
_[5]=2*W1
_[6]=b3^2
_[7]=b2*b3
_[8]=b1*b3
_[9]=r*b3
_[10]=b2^2
_[11]=b1*b2
\end{verbatim}
\end{minipage}
\begin{minipage}{0.3\linewidth}
\begin{verbatim}
_[12]=r*b2
_[13]=p1*b2+W2*b3
_[14]=W2*b2-W1*b3
_[15]=b1^2
_[16]=r*b1
_[17]=W2*b1
_[18]=W1*b1
_[19]=r^2-4*p1
_[20]=W2*r
_[21]=W1*r
_[22]=W2^2+W1*p1
\end{verbatim}
\end{minipage}

\subsection{The kernel needed in the proof of Proposition \ref{prop:icrso3}}\label{sec:singular3}

In the proof of Proposition \ref{prop:icrso3} we needed to know the kernel of the ring homomorphism
\[ \Z[p_1, w, y_1] \to \Z[b,U,x,y]/(2U,2x,U^2,xy,Ux,Uy,y^2,x^2-bU) \]
given by $p_1 \mapsto b$, $w \mapsto x$ and $y_1 \mapsto y$.

The calculation in {\sc Singular} is as follows:

\begin{verbatim}
ring P = ZZ, (p1, w, y1), dp;
ring F = ZZ, (b, U, x, y), dp;
ideal I = 2*U, 2*x, U^2, x*y, U*x, U*y, y^2, x^2-b*U;
qring R = std(I);
map f = P, b, x, y;
setring P;
kernel(R, f)
\end{verbatim}

{\sc Singular} confirms the expected answer:

\begin{verbatim}
_[1]=2*w
_[2]=y1^2
_[3]=w*y1
_[4]=w^3
\end{verbatim}

\section{Cohomology of $\RP^2$ with various local coefficient systems}

Let $\tau\in \Z/2$ denote the nontrivial element and let $\Z[\Z/2]=\{m+\tau n\,|\, m,n\in \Z,\, \tau^2=1\}$ be the regular $\Z/2$-representation. Denote by $\Z_{\omega}$ the sign representation and by $\Z$ and $\Z/2$ the respective trivial representations. There is a commutative diagram of representations
\[
\xymatrix{
\Z[\Z/2]\ar[r]^-{\varepsilon'} \ar[d]^-{\varepsilon}	& \Z_{\omega} \ar[d]^-{R} \\
\Z \ar[r]^-{R}									& \Z/2
}
\]
where $\varepsilon(m+\tau n)=m+n$, $\varepsilon'(m+\tau n)=m-n$ and $R$ is reduction modulo 2. This induces a commutative diagram of graded groups
\begin{equation} \label{dgr:coh}
\xymatrix{
H^\ast(\RP^2;\Z[\Z/2]) \ar[r]^-{\varepsilon_\ast'} \ar[d]^-{\varepsilon_\ast}	& H^\ast(\RP^2;\Z_{\omega}) \ar[d]^-{R_\ast}\\
H^\ast(\RP^2;\Z) \ar[r]^-{R_\ast}										& H^\ast(\RP^2;\Z/2)\, .
}
\end{equation}

\begin{lemma}\label{lem:rp2coh}
In the respective degrees diagram (\ref{dgr:coh}) reads
\[
H^0: \xymatrix{
\Z \ar[r] \ar[d]^-{2}	& 0 \ar[d] \\
\Z \ar[r]^-{R}			& \Z/2
}\quad\quad H^1:
\xymatrix{
0 \ar[r] \ar@{=}[d]	& \Z/2 \ar@{=}[d] \\
0 \ar[r]			& \Z/2
}\quad\quad H^2:
\xymatrix{
\Z \ar@{=}[r] \ar[d]^-{R}	& \Z \ar[d]^-{R} \\
\Z/2 \ar@{=}[r]			& \Z/2.
}
\]
\end{lemma}
\begin{proof}
Recall two facts: For any $\Z/2$-representation $V$ we have Poincar\'e duality in the form $H^\ast(\RP^2;V)\cong H_{2-\ast}(\RP^2;\Z_{\omega}\otimes V)$, and $H^\ast(\RP^2;\Z[\Z/2])\cong H^\ast(S^2;\Z)$ is the cohomology of the universal cover. Since $\Z_{\omega}\otimes \Z_{\omega}\cong\Z$, this determines the top groups in each of the three diagrams. The bottom groups and arrows in all three diagrams are clear. 
\begin{itemize}
\item \emph{Degree 0}: Under the identification $H^\ast(\RP^2;\Z[\Z/2])\cong H^\ast(S^2;\Z)$, the map $\varepsilon_\ast$ is the transfer for the covering $S^2\rightarrow \RP^2$. In degree $0$ the transfer is multiplication by the degree of the covering, i.e. multiplication by $2$.\bigskip

\item \emph{Degree 1}: We apply Poincar\'e duality to the map $R_\ast$ and note that
\[
(R: \Z_{\omega}\rightarrow \Z/2)\otimes \Z_{\omega} \cong (R: \Z\rightarrow \Z/2)\, .
\]
Since $R_\ast: H_{1}(\RP^2;\Z)\rightarrow H_1(\RP^2;\Z/2)$ is an isomorphism, the second diagram follows.\bigskip

\item \emph{Degree 2}: We apply Poincar\'e duality to the map $\varepsilon'_\ast$ and note that
\[
(\varepsilon': \Z[\Z/2]\rightarrow \Z_{\omega})\otimes \Z_{\omega} \cong (\varepsilon: \Z[\Z/2]\rightarrow \Z)\, ,
\]
since $\Z[\Z/2]\otimes \Z_{\omega}\cong \Z[\Z/2]$ via $m+\tau n\mapsto m-\tau n$. Therefore, the top arrow in the third diagram is $\varepsilon_\ast: H_0(\RP^2;\Z[\Z/2])\rightarrow H_0(\RP^2;\Z)$. It can be identified with the map of coinvariant modules $\Z[\Z/2]_{\Z/2}\rightarrow \Z$ induced by $\varepsilon$, which is an isomorphism. The right arrow follows again by applying Poincar\'e duality to the right arrow $R_\ast$ in diagram (\ref{dgr:coh}). The left arrow is then determined by commutativity of the diagram.
\end{itemize}
This finishes the proof.
\end{proof}

\end{document}